\documentclass[hyperref={colorlinks,citecolor=red,urlcolor=black,bookmarks=false,hypertexnames=true}]{article}

\usepackage{arxiv}

\usepackage[utf8]{inputenc} 
\usepackage[T1]{fontenc}    
\usepackage{url}            
\usepackage{booktabs}       
\usepackage{amsfonts}       
\usepackage{nicefrac}       
\usepackage{microtype}      
\usepackage{lipsum}         
\usepackage{graphicx}
\usepackage{subcaption}
\usepackage{doi}

\usepackage{colortbl}
\usepackage{listings}
\usepackage{bm}
\usepackage{physics}

\usepackage{amsmath,accents}
\usepackage{cleveref}       

\usepackage[ruled,vlined]{algorithm2e}

\usepackage{appendix}
\usepackage{cite}

\usepackage{multirow}
\usepackage{array}

\usepackage{chngcntr}
\usepackage{epigraph}

\SetKwComment{Comment}{/* }{ */}

\newcommand{\ttl}{%
	\begin{minipage}{0.1\textwidth}
		{\includegraphics[height=2cm]{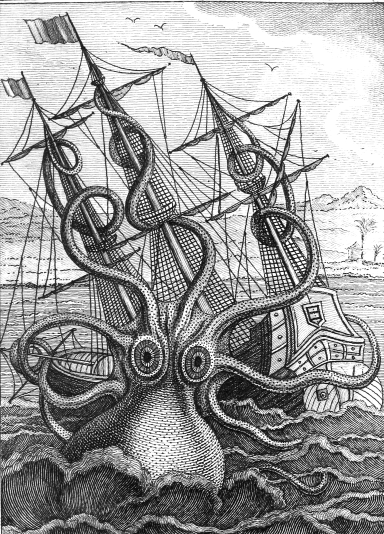}}    
	\end{minipage}
	\begin{minipage}{0.85\textwidth}
		\centering
		Are Two Hidden Layers Still Enough for the Physics--Informed Neural Networks?
	\end{minipage}   
}

\title{\ttl}

\date{}

\author{ \href{https://orcid.org/0000-0002-4930-1846}{\includegraphics[scale=0.06]{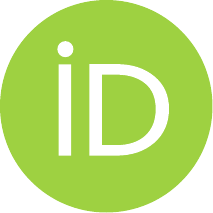}\hspace{1mm}Vasiliy A. Es'kin}\thanks{Corresponding author: Vasiliy Alekseevich Es’kin (\href{vasiliy.eskin@gmail.com}{vasiliy.eskin@gmail.com})} \\
	Department of Radiophysics, University of Nizhny Novgorod\\
	 Nizhny Novgorod, Russia, 603950\\
	 and\\
	 Huawei Nizhny Novgorod Research Center\\
	 Nizhny Novgorod, Russia\\
	\href{vasiliy.eskin@gmail.com}{\texttt{vasiliy.eskin@gmail.com}} \\
	\And
	\href{https://orcid.org/0000-0001-5210-8281}{\includegraphics[scale=0.06]{orcid.pdf}\hspace{1mm}Alexey O. Malkhanov} \\
	Huawei Nizhny Novgorod Research Center\\
	Nizhny Novgorod, Russia\\	\texttt{alexey.malkhanov@gmail.com} \\
	\And
	\href{https://orcid.org/0000-0002-0454-5249}{\includegraphics[scale=0.06]{orcid.pdf}\hspace{1mm}Mikhail E. Smorkalov} \\
	Skolkovo Institute of Science and Technology\\
	Moscow, Russia\\
	and\\
	Huawei Nizhny Novgorod Research Center\\
	Nizhny Novgorod, Russia\\	\texttt{smorkalovme@gmail.com} \\
}

\renewcommand{\headeright}{}
\renewcommand{\undertitle}{}
\renewcommand{\shorttitle}{A preprint}

\renewcommand{\vec}{\bf}
\newcommand{\rasymbol}[1]{\accentset{\circ}{#1}}

\usepackage{tikz}
\newcommand\tikznode[3][]%
{\tikz[remember picture,baseline=(#2.base)]
	\node[minimum size=0pt,inner sep=0pt,#1](#2){#3};%
}

\hypersetup{
pdftitle={A template for the arxiv style},
pdfsubject={math.na, cs.AI, physics.comp-ph},
pdfauthor={Vasiliy A. Es'kin, Alexey O. Malkhanov, Mikhail E. Smorkalov},
pdfkeywords={Deep Learning, Physics-informed Neural Networks, Partial differential equations, Predictive modeling, Computational physics, Nonlinear dynamics},
}

\begin{document}
\maketitle

\epigraph{It is not the gods who burn pots}{An ancient Greek proverb (attributed to $\Sigma{o}\lambda\omega\nu$)}


\begin{abstract}
The article discusses the development of various methods and techniques for initializing and training neural networks with a single hidden layer, as well as training a separable physics-informed neural network consisting of neural networks with a single hidden layer to solve physical problems described by ordinary differential equations (ODEs) and partial differential equations (PDEs). A method for strictly deterministic initialization of a neural network with one hidden layer for solving physical problems described by an ODE is proposed. Modifications to existing methods for weighting the loss function ($\delta$-causal training and gradient normalization) are given, as well as new methods developed for training strictly deterministic-initialized neural networks to solve ODEs (detaching, additional weighting based on the second derivative, predicted solution-based weighting, relative residuals). An algorithm for physics-informed data-driven initialization of a neural network with one hidden layer is proposed. A neural network with pronounced generalizing properties is presented, whose generalizing abilities of which can be precisely controlled by adjusting network parameters. A metric for measuring the generalization of such neural network has been introduced. A gradient-free neuron-by-neuron (NbN) fitting method has been developed for adjusting the parameters of a single-hidden-layer neural network, which does not require the use of an optimizer or solver for its implementation. The proposed methods have been extended to 2D problems using the separable physics-informed neural networks (SPINN) approach. Numerous experiments have been carried out to develop the above methods and approaches. Experiments on physical problems, such as solving various ODEs and PDEs, have demonstrated that these methods for initializing and training neural networks with one or two hidden layers (SPINN) achieve competitive accuracy and, in some cases, state-of-the-art results.

\end{abstract}


\usetikzlibrary {arrows.meta,bending,positioning}

\section{Introduction}

The last decades have been notable by significant achievements in the field of machine learning for the science and engineering~\cite{Ronneberger2021,bodnar2024foundationmodelearth,Raissi2019,Wang2022,Wang2022_2,Lu_2021,Krinitskiy2022,Fanaskov2022,Raissi2017,eskin2023optimal}. Notable examples include AlphaFold~\cite{Ronneberger2021}, which has been used to predict protein structures and earned its creators the Nobel prize, and Aurora~\cite{bodnar2024foundationmodelearth}, which is used for weather forecasting. These examples demonstrate the integration of physical laws and principles into machine learning models, leading to the development of physics-informed machine learning techniques. Recall that auxiliary systems for solving scientific problems can be categorized into those that are trained using previously acquired data (from experiments or calculations) and make predictions based on this information~\cite{Lu_2021, Krinitskiy2022, Fanaskov2022}, those that apply previously established physical principles~\cite{Raissi2019, Wang2022, Wang2022_2}, and those that utilize both approaches~\cite{Raissi2017}. Previous century studies explored the initial ideas for constraining neural networks with physical laws in~\cite{Psichogios1992AHN} and~\cite{Lagaris_1998}. More recent research in the area was realized a few years ago with physics-informed neural networks (PINNs) using modern computational tools~\cite{Raissi2019}. In the PINN approach, a neural network is trained to approximate the dependences of physical values on spatial and temporal variables for a given physical problem, described by a set of physical equations, together with additional constraints such as initial and boundary conditions. In some cases, data from numerical simulations or experiments can also be used to help train the network. This method has been used in recent years to solve a wide range of problems described by ordinary differential equations~\cite{Raissi2017,Raissi2019}, integro-differential equations~\cite{Yuan2022}, nonlinear partial differential equations~\cite{Raissi2017,Raissi2019,Kharazmi2021}, PDE with noisy data~\cite{Yang2021}, etc.~\cite{Cuomo2022}, related to various fields of knowledge such as thermodynamics~\cite{Patel2022,Kharazmi2021}, hydrodynamics~\cite{Cai2021}, mechanics~\cite{Moseley2020}, electrodynamics~\cite{Kharazmi2021}, geophysics~\cite{He2020}, systems biology~\cite{daryakenari2023aiaristotle}, finance~\cite{Yuan2022}, etc.

Despite the success of researches for using PINN approach for scientific problems, most current works employ small neural networks, with only a few layers, or neural networks with specific architectures that have a large number of layers~\cite{wang2024,jiang2024}. Such choices of neural networks are associated with the problem of vanishing gradients in deep PINNs, which reduces the representativeness of PINNs and limits their effectiveness, and for the mentioned networks this problem is relatively mitigated. At the same time, neural networks with one or two hidden layers are devoid of this drawback. In our opinion, the potential of these neural networks has not yet been fully explored despite the fact that they constitute a relatively simple and straightforward subject of study.

It follows from the universal approximation theorem that any continuous function can be arbitrarily closely approximated by a multi-layer perceptron with only one hidden layer and a finite number of neurons~\cite{HORNIK1989359,Cybenko1989ApproximationBS,Yarotsky2016ErrorBF}. Currently, two approaches to this issue are presented in the literature: shallow neural networks, which have a single hidden layer and can still produce any non-linear continuous function; deep neural networks (DNN), that use more than two layers of neural networks to model complicated relationships of entities and values. Some researchers argue that two hidden layers are sufficient for approximating any complex function~\cite{Arnold1957,Kolmogorov1956,Kolmogorov1957,Maiorov1999LowerBF,GULIYEV2018262}, while others claim that three hidden layers are sufficient and necessary~\cite{SHEN201974,Shen2020NeuralNA}. In our work, we adhere to the opinion of Kolmogorov, Arnold, Guliyev, and Ismailov~\cite{Arnold1957,Kolmogorov1956,Kolmogorov1957,Maiorov1999LowerBF,GULIYEV2018262}, and for the PINN approach, we limit the consideration to neural networks with one hidden layer for ODEs and two layers for PDEs or parametrized ODEs. By a two-layer network, we mean several single-hidden-layer neural networks that form a larger neural network through multiplication, within the framework of the separable PINN (SPINN) approach.


This work covers all aspects of PINN for solving physical problems: the architecture of a neural network, the initialization of its trainable parameters, the weighting of elements of the loss function, and the training procedure. The development of methods and techniques related to these aspects of PINNs is presented in this paper in the form of an iterative progression from simple to more complex, through a series of cycles in the following procedure:
\begin{center}
	\begin{tikzpicture}[outer sep=auto]
		\draw[line width=5pt]
		(-1.2,0)  node (f0) {}
		(0,0)  node (f) {Suggestions}
		(3.8,0)  node (d) {Numerical Experiments}
		(8.8,0)  node (s) {Discussion of disadvantages}
		(10.8,0)  node (s0) {};
		
		\draw[->] (f) -- (d);
		\draw[->] (d) -- (s);
		\draw[overlay, arrows = {->[]}] ++(s0.east) -- +(0.2,0) |- +(0,-0.5) -| (f0.west) -- +(0.3,0);  
	\end{tikzpicture}
\end{center}
This approach allows for a systematic and iterative process of refining and improving the methods and techniques used in PINN. We will briefly list the contributions made in the paper:
\begin{enumerate}
	\item A method of strictly deterministic initialization of a neural network with one hidden layer is proposed for solving physical problems described by the ODE;
	\item Modifications of the existing methods of weighting the loss function are given ($\delta$-causal training, gradient normalization) and new weighting methods are developed for training strictly deterministic initialized neural networks for solving the ODE (detaching, additional weighting based on second derivative, predicted solution based weighting (PSBW), relative residuals (RR));
	\item A algorithm of physics-informed data-driven (PIDD) initialization of a neural network with one hidden layer is proposed;
	\item A neural network with pronounced generalizing properties is presented, the generalizing properties of which can be precisely controlled by changing network parameters;
	\item A measure of generalization of the above-mentioned neural network has been introduced;
	\item A gradient-free neuron-by-neuron (NbN) fitting method has been developed to adjust the parameters of a neural network with one hidden layer, which does not require an optimizer or solver in its application;
	\item The proposed methods are extended to 2D problems using the separable PINN (SPINN) approach;
	\item Numerous experiments have been carried out to develop the above methods and approaches.
\end{enumerate}
Our approaches for initializing and training neural networks with one or two hidden layers open new frontiers for traditional scientific research and contribute to address challenges in machine learning for science, such as robustness, interpretability, and generalization. 

The paper is structured as follows. In section 2, the methods of the strictly deterministic initialization of the trainable parameters of neural networks for solving ODE problems are suggested and described. Section 3 includes the descriptions of modifications to training and a new approach to training neural networks for the solving of ODE problems. Section 4 presents physics-informed data-driven initialization and neuron-by-neuron training of neural networks. In the same section, there is a neural network, which possesses pronounced generalization properties. In Section 5, the architecture of the neural network, initialization, and training procedures for the solving of PDE are given. Finally, in Section 6 concluding remarks are given.

\section{Strictly deterministic initialization of the weights and biases}

In this section, we consider a neural network with a single hidden layer. We are inspired by the simple Euler method for numerical integration of ordinary differential equations (ODEs) when developing methods for initializing the weights and biases of this type of neural network.

Consider nonlinear ordinary differential equations, which depend on coordinates $x$, and in general, take the following form
\begin{equation}\label{eq1}
	\frac{\partial \vec u}{\partial x} + \mathcal{N} [{\vec u}, x] = 0,\quad x\in [0, X]
\end{equation}
under the initial conditions
\begin{equation}
	{\vec u}(0) = {\vec g}, \label{eq2}
\end{equation}
where ${\vec u}(x)$ denotes the latent solution that is governed by the ODE system of equations (\ref{eq1}), consisting of $n$ components ${\vec u} = (u_1, u_2,\dots,u_n)$, $\mathcal{N}$ is a nonlinear differential operator, ${\vec {g}}$ is the initial distribution of ${\vec u}$.

According to the PINN approach~\cite{Raissi2019} we approximate the unknown solution ${\vec u}(x)$ using neural networks ${\vec u}_{{\bm \theta}}(x)$, each of which the component ${u}_{{\bm \theta};l}(x)$ is a separate neural network with a single hidden layer, as follows:
\begin{equation}\label{eq3}
	{u}_{{\bm \theta};l}(x) = \sum\limits_{k = 0}^{N-1} {W}^{(2)}_k\sigma\left({W}^{(1)}_k x + {b}^{(1)}_k\right) + {b}^{(2)}_0,
\end{equation}
where ${W}^{(j)}_k$ is $k$th weight of the $j$th layer, $b^{(j)}_k$ is $k$th the bias of the $j$th layer, ${\bm \theta}$ denote all trainable parameters ({weights} and biases) of the neural network ${u}_{{\bm \theta};l}$, $\sigma$ is an activation function, and $N$ is number neurons in the hidden layer. 

Finding optimal parameters is an optimization problem, which requires {the definition} of a loss function such that its minimum gives the solution of the ODE. The physical-informed model is trained by minimizing a composite loss function which consists of the local residuals of the differential equation over the problem domain and its initial conditions as shown below:
\begin{equation}\label{eq4}
	\mathcal{L}({\bm \theta}) = \lambda_{ic} \mathcal{L}_{ic}({\bm \theta}) + \lambda_{r} \mathcal{L}_{r}({\bm \theta}),	  
\end{equation}
where
\begin{eqnarray}
	&& \mathcal{L}_{ic} \left( {\bm \theta} \right)  =  \left| {\vec u}_{{\bm \theta}}\left(0\right) - {\vec g}  \right|^{2},\label{eq5}\\
	&& \mathcal{L}_{r}\left(\bm \theta \right) = \frac{1}{N_r} \sum_{i=1}^{N_{r}} \left| \mathcal{R}\left[{\vec u}_{{\bm \theta}} \right] \left({x}_{i}^{(r)} \right)\right|^{2},\label{eq6}\\
	&& \mathcal{R}\left[{\vec u} \right]:= \frac{\partial \vec u}{\partial x} + \mathcal{N}\left[ {\vec u}\left(x \right), x  \right].\label{eq7}
\end{eqnarray}
Here $\left\{x_{i}^{(r)}\right\}^{N_{r}}_{i=1}$ is a set of points sampled from the ODE domain. All required gradients w.r.t. input variable $x$ and network parameters $\bm \theta$ can be efficiently computed via automatic differentiation~\cite{Griewank2008} with algorithmic accuracy, which is defined by the precision of the computation system. Note, the hyperparameters $\lambda_{ic}$ and $\lambda_{r}$ allow for separate tuning of the learning rate for each loss term in order to enhance the convergence of the model~\cite{Wang2021,WANG2022110768}.

The optimization problem can be defined as follows
\begin{equation}
	{\bm \theta}^* = \underset{{\bm \theta}}{\arg}\min \mathcal{L}({\bm \theta}),\label{eq8}
\end{equation}
where ${\bm \theta}^*$ are optimal parameters of the neural network which minimize the discrepancy between the exact unknown solution ${\vec u}$ and the approximate one ${\vec u}_{\bm \theta^*}$.

\begin{figure}[th!]
	\centering
	\includegraphics[width=0.9\textwidth]{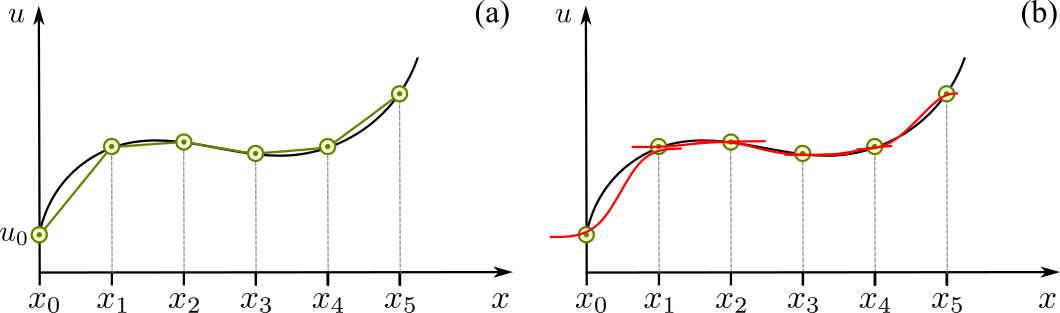}
	\caption{Approximations of the graph of a given function $u(x)$ by a polyline curve (a) and segments of sigmoids (b).}
	\label{fig1}
\end{figure}

{\bf Suggestion 1: Activation function.} We take the sigmoid function as the activation function, since this allows us to naturally approximate the piecewise linear function obtained by the Euler integration method (see the explanatory Figure~\ref{fig1}), and since it has been proven that a neural network with such an activation function is an universal approximator~\cite{Cybenko1989ApproximationBS}.

\begin{figure}[th!]
	\centering
	\includegraphics[width=0.9\textwidth]{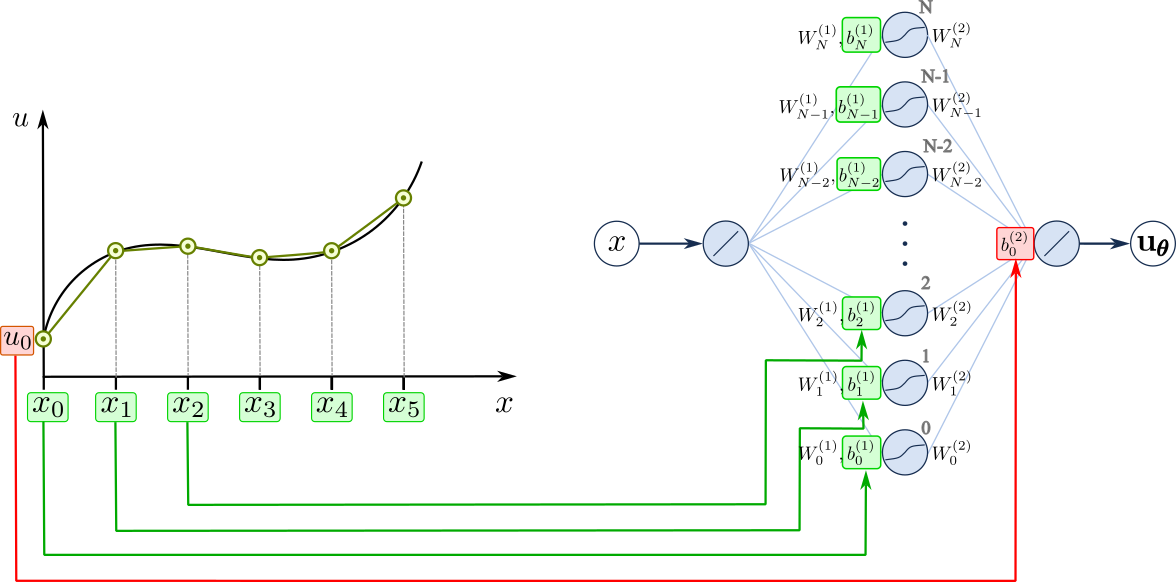}
	\caption{Explanatory figure for Algorithm~\ref{alg1}.}
	\label{fig2}
\end{figure}

{\bf Suggestion 2: Initialization.} We initialize the parameters of a neural network as follows (see the explanatory Figure~\ref{fig2}):
\begin{enumerate}
	\item Weights of the hidden layer are taken $1$ (${W}^{(1)}_k = 1$), as suggested in~\cite{GULIYEV2018262}.
	\item We divide the interval $[0, X]$ into $N$ identical segments. Then we assign each bias $b^{(1)}_k$ a value equal to the coordinate $-x_k-\Delta x$ of the beginning of $k$th interval ($b^{(1)}_k=-x_k-\Delta x$, where $\Delta x = X / N$ and $x_0=0$). In this case of choosing weights, the $k$th neuron of the hidden layer shows the largest response (the derivative of the output value of the neuron is maximal) only if the input coordinate $x$ is close to the value $x_{k+1}$.
	\item The weights and biases of the hidden layer are fixed as it was suggested in~\cite{GULIYEV2018262}. Note, the values of centres of the neurons of physics-informed radial basis network (PIRBN)~\cite{Bai_2023} are assigned in a similar way.
	\item The bias of the output layer is taken equal to initial value $u_0$ ($b^{(2)}_0=u_0, u_0 = u(0)$).
	\item The parameters ${W}^{(2)}_k$ are initialized with the Glorot scheme~\cite{Glorot10a}.
\end{enumerate}
As a result, we have the following neural network initialization Algorithm~\ref{alg1}.

\begin{algorithm}[h!]
	\caption{Initialization of physics-informed neural networks}\label{alg1}
	\KwData{---}
	\KwResult{Initialized neural network ${u}_{{\bm \theta};l}(x)$ of PINN ${\vec u}_{\bm \theta}$, which consists of $N$ neurons on hidden layer}
	$\Delta x \gets X / N$\;
	\For{$k=0,\dots, N-1$}{
		${W}^{(1)}_k \gets 1$\;
		${b}^{(1)}_k \gets - (k+1) \Delta x$\;
	}
	${W}^{(2)}_k$ are initialized with the Glorot scheme~\cite{Glorot10a}\;
	${b}^{(2)}_0 \gets {u_l}(0)$.
\end{algorithm}
Note that the global error in the Euler method is of the order $X/N$. For this reason, we can expect the same order of prediction error for a neural network initialized in a given way.

\subsection{Numerical experiments}
For our experiments we used Pytorch~\cite{Paszke2019} version 2.1.2 and the training was carried out on a node with a GPU with characteristics similar to industry-leading GPU and CPU.


The networks had been trained via stochastic gradient descent using the Adam optimizer~\cite{Kingma2014} from the \verb*|pytorch| and LBFGS from the python library \verb*|pytorch-minimize|.

In the examples given in this subsection, we investigate how neural network parameters and training parameters affect the accuracy of predicting solutions when using the neural network initialization Algorithm~\ref{alg1}. The relative ${\mathbb L}_2$ error is used as a measure of accuracy, which is defined in the application~\ref{appA}. In all experiments, a uniform distribution of coordinates for collocation points is taken.

\subsubsection{Example 1. Harmonic Oscillator}
Consider the harmonic oscillator, which is governed by a system of two first-order equations for the $t\in [0, T]$
\begin{eqnarray}
	&& \frac{d u_1}{d t} = u_2, \quad t \in [0,T],\label{eq9}\\
	&& \frac{d u_2}{d t} = - \frac{k}{m} u_1,\label{eq10}\\
	&& u_1(0) = 1,\quad u_2(0) = 0 \label{eq11},
\end{eqnarray}
where $m$ is the mass and $k$ is the spring force constant. We used the following parameters: $T=100$, $m=1$, $k=1$. The exact analytical solution of this problem is $u_1^{\text{ref}} = \cos\left(t\right)$, $u_2^{\text{ref}} = - \sin\left(t\right)$.
The latent variables $u_1$ and $u_2$ are represented by neural networks $u_{\bm \theta; 1}$ and ${u}_{\bm \theta;2}$ described above. Consider how the accuracy of the solution provided by a neural network depends on the number of neurons $N$ in the hidden layer, the number of collocation points $N_t$ and the upper limit of the interval under consideration $T$. 

In our experiments, we used two-stage and three-stage learning. Two-stage training consisted of 1000 epochs of optimization by Adam optimizer with learning rate $10^{-3}$ and 5 epochs of optimization by LBFGS optimizer (maximal number of iterations per optimization step is 5000) under frozen parameters of the hidden layer of neural network. Three-stage training is the two-stage training with an additional stage of 5 epochs of optimization by LBFGS optimizer (maximum number is previous) under unfrozen parameters of the hidden layer of neural network. Such an additional training procedure with unfrozen parameters should lead to an increase in the accuracy of the training result. In our calculations the weights of loss (\ref{eq4}) were taken as $\lambda_{ic} =1$ and $\lambda_{r} = 10$. We used 5000 uniform distributed collocation points for calculating the relative error trained model.


\begin{figure}[t!]
	\centering
	\begin{subfigure}[t]{0.4\textwidth}
		\centering
		\includegraphics[width=\linewidth]{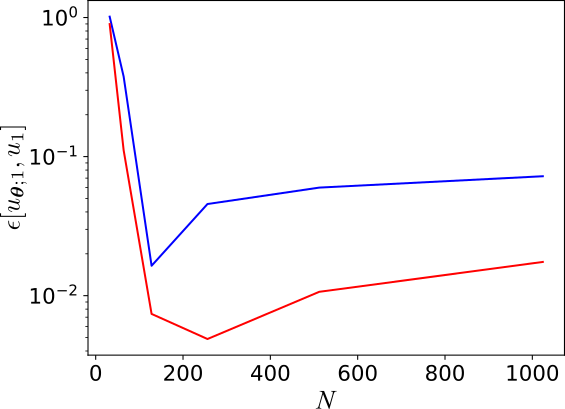}
		\caption{$u_1$}
	\end{subfigure}\hspace{5mm}
	\begin{subfigure}[t]{0.4\textwidth}
		\centering
		\includegraphics[width=\linewidth]{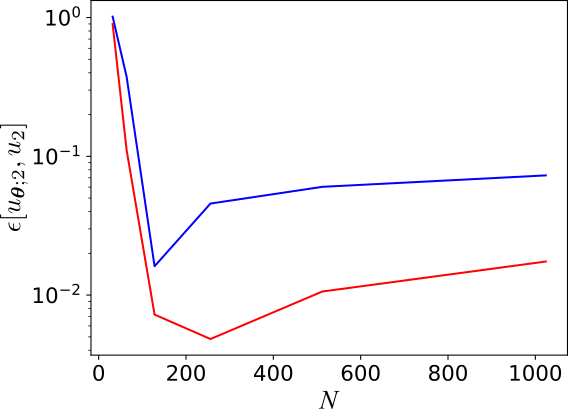}
		\caption{$u_2$}
	\end{subfigure}
	\caption{Harmonic Oscillator. (a) and (b) are dependences of relative ${\mathbb L}_2$ error $\epsilon[u_{\bm \theta; 1}, u_{1}]$ and $\epsilon[u_{\bm \theta; 2}, u_{2}]$, respectively, on the number of neurons $N$. Blue and red solid lines correspond to two-stage and three-stage training. Number of collocation points is $N_t=2048$, maximum of time is $T=100$.}
	\label{fig3}
\end{figure}

Figure~\ref{fig3} shows dependences of relative ${\mathbb L}_2$ error $\epsilon[u_{\bm \theta; 1}, u_{1}]$ and $\epsilon[u_{\bm \theta; 2}, u_{2}]$, respectively, on the number of neurons $N$ under two-stage and three-stage training (definition of $\epsilon$ is given in appendix~\ref{appA}). The figures show that the minimum error value is achieved with 128 neurons for two-stage learning and 256 neurons for three-stage learning. The relative error obtained in three-stage learning is much smaller than the error obtained in two-stage learning.

\begin{figure}[t!]
	\centering
	\begin{subfigure}[t]{0.4\textwidth}
		\centering
		\includegraphics[width=\linewidth]{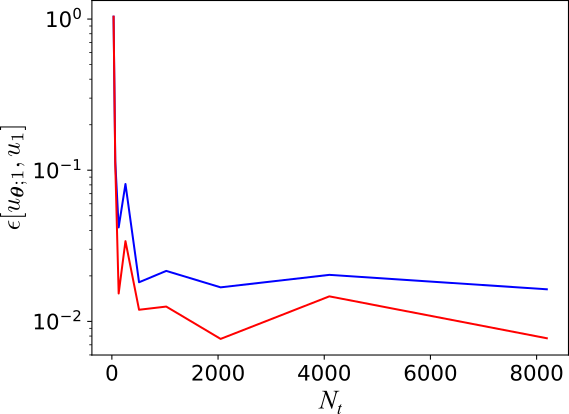}
		\caption{$u_1$}
	\end{subfigure}\hspace{5mm}
	\begin{subfigure}[t]{0.4\textwidth}
		\centering
		\includegraphics[width=\linewidth]{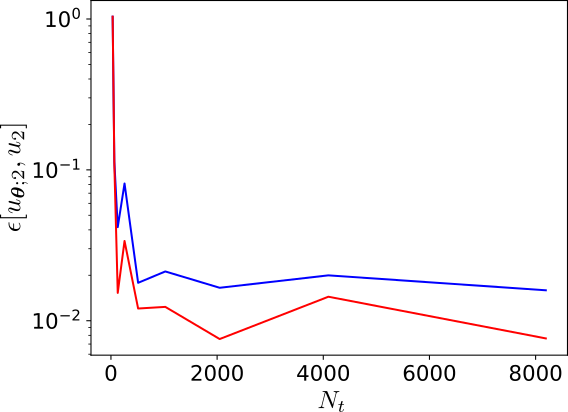}
		\caption{$u_2$}
	\end{subfigure}
	\caption{Harmonic Oscillator. (a) and (b) are dependences of relative ${\mathbb L}_2$ error $\epsilon[u_{\bm \theta; 1}, u_{1}]$ and $\epsilon[u_{\bm \theta; 2}, u_{2}]$, respectively, on the number of collocation points $N_t$. Blue and red solid lines correspond to two-stage and three-stage training. Number of neurons is $N=128$, maximum of time is $T=100$.}
	\label{fig4}
\end{figure}

Figure~\ref{fig4} shows dependences of relative ${\mathbb L}_2$ error for $u_{\bm \theta; 1}$ and ${u}_{\bm \theta;2}$, respectively, on the number of collocation points $N_t$ under two-stage and three-stage training. The figures show that the minimum error value is achieved with 2048 collocation points for both scenarios of trainings. The relative error obtained in three-stage learning is much less than the error obtained in two-stage training.

\begin{figure}[t!]
	\centering
	\begin{subfigure}[t]{0.4\textwidth}
		\centering
		\includegraphics[width=\linewidth]{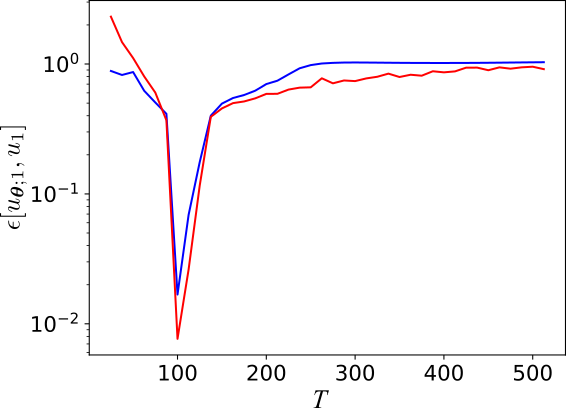}
		\caption{$u_1$}
	\end{subfigure}\hspace{5mm}
	\begin{subfigure}[t]{0.4\textwidth}
		\centering
		\includegraphics[width=\linewidth]{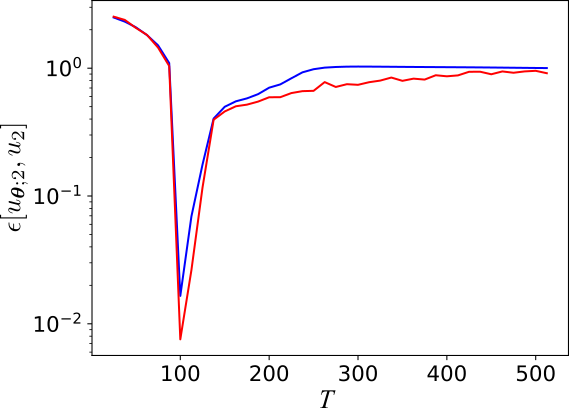}
		\caption{$u_2$}
	\end{subfigure}
	\caption{Harmonic Oscillator. (a) and (b) are dependences of relative ${\mathbb L}_2$ error $\epsilon[u_{\bm \theta; 1}, u_{1}]$ and $\epsilon[u_{\bm \theta; 2}, u_{2}]$, respectively, on the boundary $T$. Blue and red solid lines correspond to two-stage and three-stage training. Number of neurons is $N=128$,  number of collocation points is $N_t=2048$.}
	\label{fig5}
\end{figure}

Figure~\ref{fig5} shows dependences of relative ${\mathbb L}_2$ error for $u_{\bm \theta; 1}$ and ${u}_{\bm \theta;2}$, respectively, on the boundary $T$ under two-stage and three-stage trainings. The figures show that the minimum error value is achieved with value $T=100$ for both trainings. The relative error obtained in three-stage learning is much less than the error obtained in two-stage learning.

\begin{figure}[t!]
	\centering
	\begin{subfigure}[t]{0.4\textwidth}
		\centering
		\includegraphics[width=\linewidth]{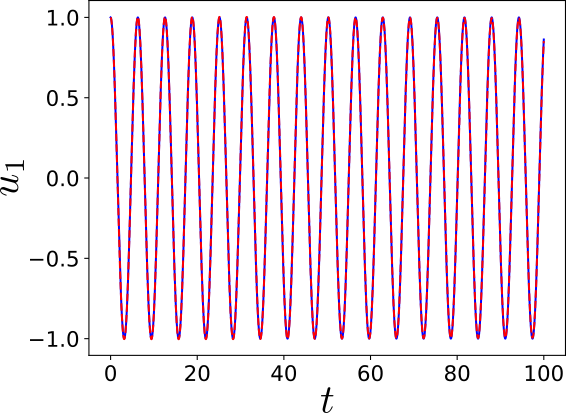}
		\caption{$u_1$}
	\end{subfigure}\hspace{5mm}
	\begin{subfigure}[t]{0.4\textwidth}
		\centering
		\includegraphics[width=\linewidth]{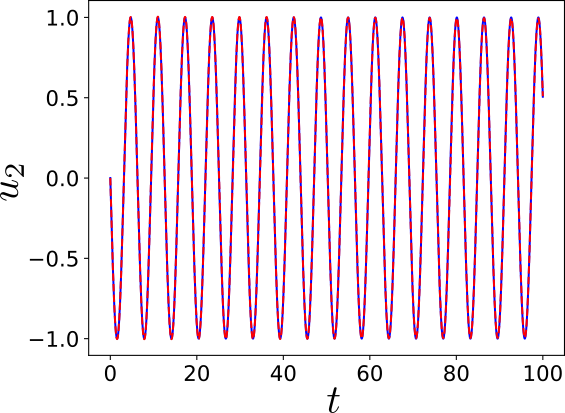}
		\caption{$u_2$}
	\end{subfigure}\vspace{5mm}
	\begin{subfigure}[t]{0.4\textwidth}
		\centering
		\includegraphics[width=\linewidth]{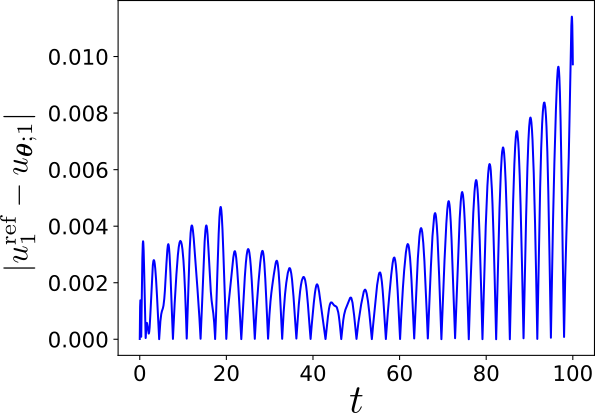}
		\caption{$u_1$}
	\end{subfigure}\hspace{5mm}
	\begin{subfigure}[t]{0.4\textwidth}
		\centering
		\includegraphics[width=\linewidth]{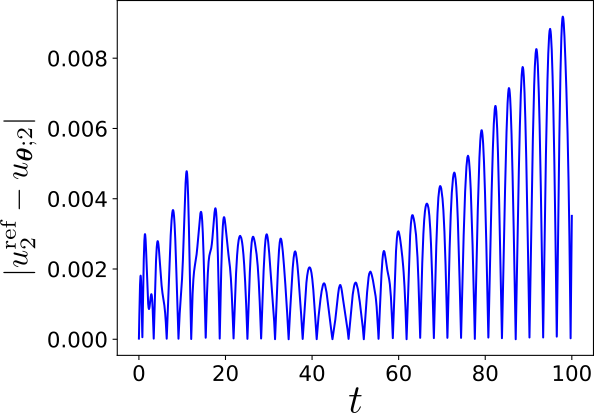}
		\caption{$u_2$}
	\end{subfigure}
	\caption{Harmonic Oscillator. (a) and (b) are comparisons of the predicted (red dash lines) and reference solutions (blue solid lines) corresponding to $u_1$ and $u_2$, respectively. (c) and (b) are absolute errors $|u_1^{\text{ref}} - u_{\bm \theta; 1}|$ and $|u_2^{\text{ref}} - {u}_{\bm \theta;2}|$, respectively. Number of neurons is $N=128$, the number of collocation points is $N_t = 2048$, maximum of time is $T=100$.}
	\label{fig6}
\end{figure}
Results of three-stage training for a neural network with $N=128$, number of collocation points $N_t=2048$ and $T=100$ are shown on Figure~\ref{fig6}. The relative ${\mathbb L}_2$ errors are $\epsilon[{u}_{\bm \theta;1}, {u_1^{\text{ref}}}] = 4.71{\times} 10^{-3}$ and $\epsilon[{u}_{\bm \theta;2}, {u_2^{\text{ref}}}] = 4.61{\times} 10^{-3}$. The training was carried out on CPU and took 9459 iterations of LBFGS optimizer and time $87.38$ seconds.

\begin{figure}[t!]
	\centering
	\begin{subfigure}[t]{0.2\textwidth}
		\centering
		\includegraphics[width=\linewidth]{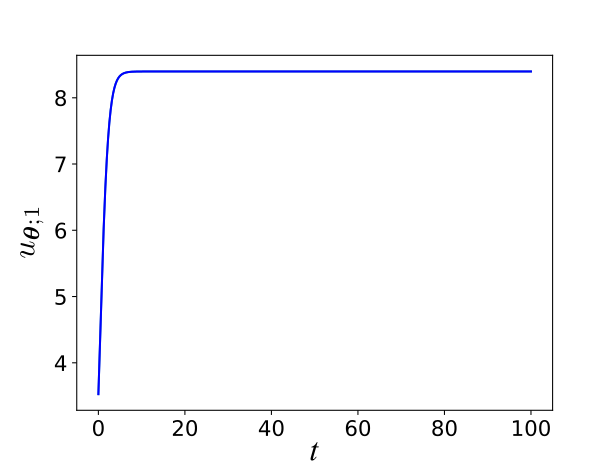}
		\caption{1 neurons}
	\end{subfigure}
	\begin{subfigure}[t]{0.2\textwidth}
		\centering
		\includegraphics[width=\linewidth]{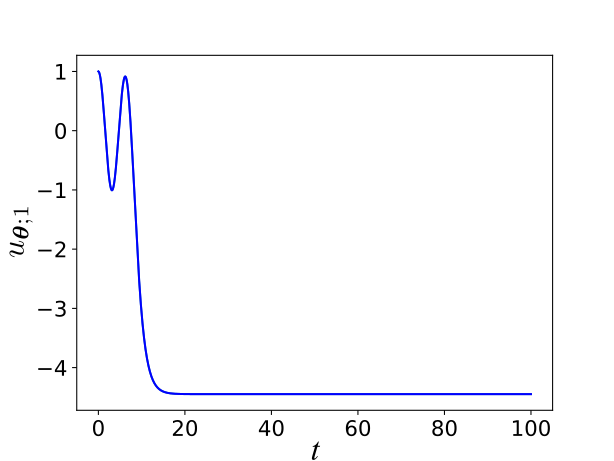}
		\caption{11 neurons}
	\end{subfigure}
	\begin{subfigure}[t]{0.2\textwidth}
		\centering
		\includegraphics[width=\linewidth]{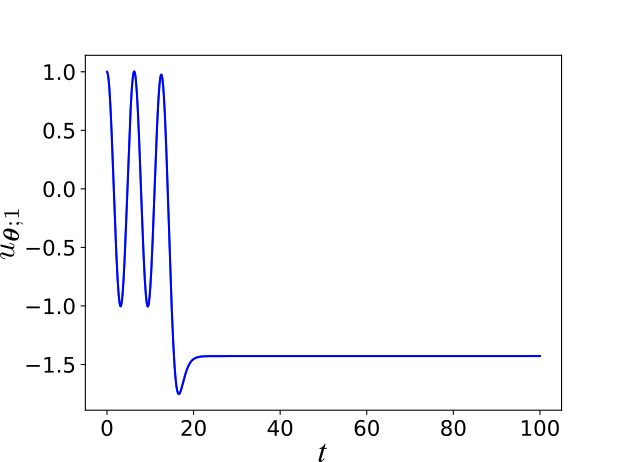}
		\caption{21 neurons}
	\end{subfigure}
	\begin{subfigure}[t]{0.2\textwidth}
		\centering
		\includegraphics[width=\linewidth]{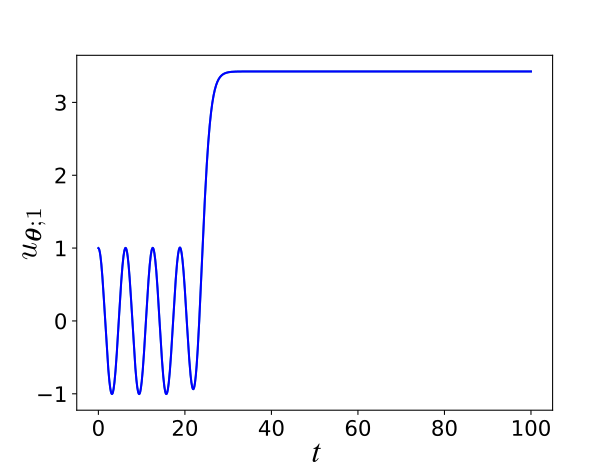}
		\caption{31 neurons}
	\end{subfigure}
	\begin{subfigure}[t]{0.2\textwidth}
		\centering
		\includegraphics[width=\linewidth]{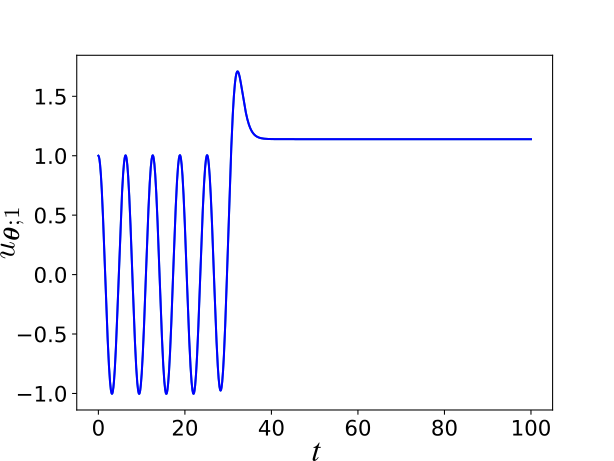}
		\caption{41 neurons}
	\end{subfigure}
	\begin{subfigure}[t]{0.2\textwidth}
		\centering
		\includegraphics[width=\linewidth]{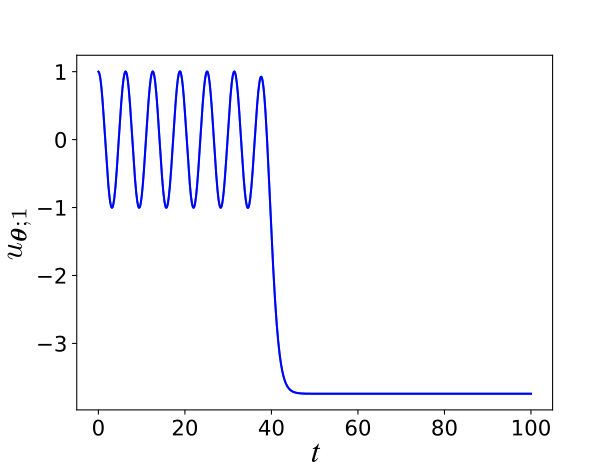}
		\caption{51 neurons}
	\end{subfigure}
	\begin{subfigure}[t]{0.2\textwidth}
		\centering
		\includegraphics[width=\linewidth]{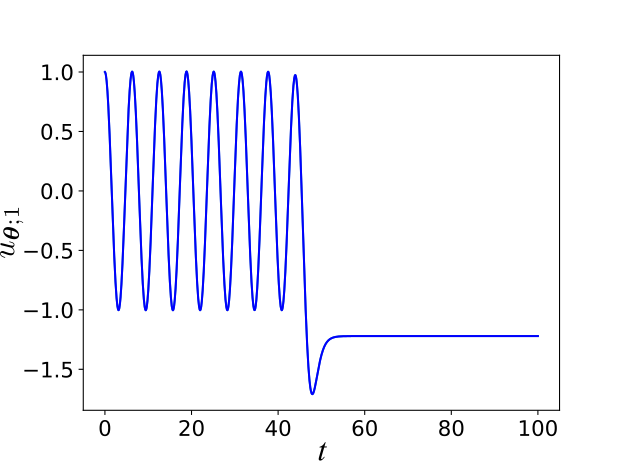}
		\caption{61 neurons}
	\end{subfigure}
	\begin{subfigure}[t]{0.2\textwidth}
		\centering
		\includegraphics[width=\linewidth]{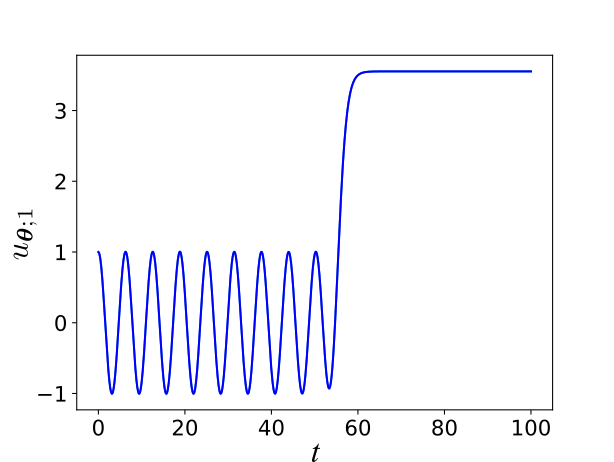}
		\caption{71 neurons}
	\end{subfigure}
	\begin{subfigure}[t]{0.2\textwidth}
		\centering
		\includegraphics[width=\linewidth]{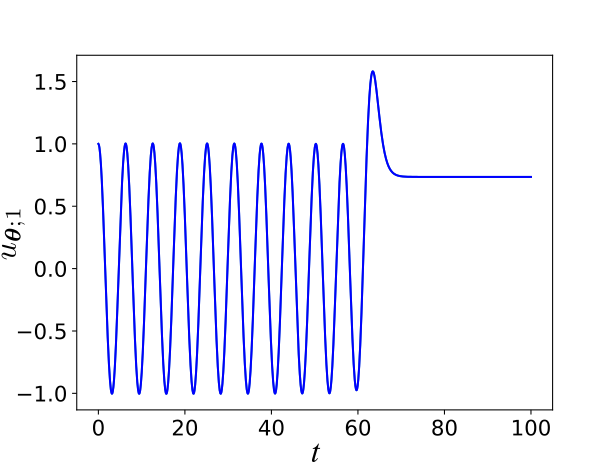}
		\caption{81 neurons}
	\end{subfigure}
	\begin{subfigure}[t]{0.2\textwidth}
		\centering
		\includegraphics[width=\linewidth]{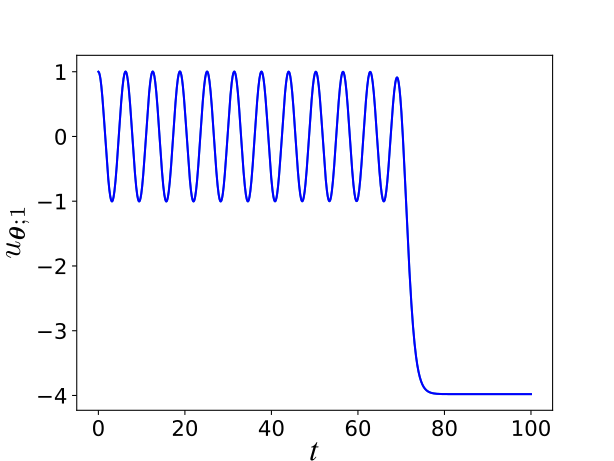}
		\caption{91 neurons}
	\end{subfigure}
	\begin{subfigure}[t]{0.2\textwidth}
		\centering
		\includegraphics[width=\linewidth]{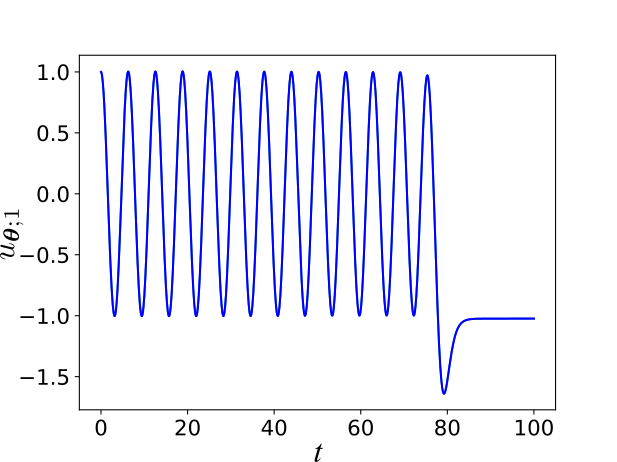}
		\caption{101 neurons}
	\end{subfigure}
	\begin{subfigure}[t]{0.2\textwidth}
		\centering
		\includegraphics[width=\linewidth]{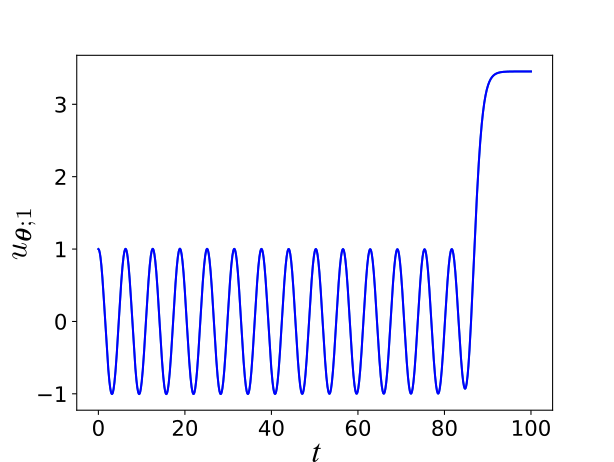}
		\caption{111 neurons}
	\end{subfigure}
	\caption{Harmonic Oscillator. Contribution to the solution $u_{\bm \theta; 1}$ of the first 1 (a), 11 (b), 21 (c), 31 (d), 41 (e), 51 (f), 61 (g), 71 (h), 81 (i), 91 (j), 101 (k), 111 (l) neurons.}
	\label{fig7}
\end{figure}
Contribution to the solution $u_{\bm \theta; 1}$ of the first neurons are presented on Figure~\ref{fig7}. As can be seen from these figures, with an increase in the number of neurons in the trained neural network, the interval of $t$ where the solution is close to the truth solution increases.

\subsubsection{Example 2. The electric potential of the charged sphere}
Consider the electric potential of the charged sphere, which is governed by the equation for the $r\in [0, R]$
\begin{eqnarray}
	&& \frac{d u}{d r} = - \frac{1}{r^2}, \quad r \in [a,R],\label{eq12}\\
	&& u(a) = 1 / a, \label{eq13}
\end{eqnarray}
where $a$ is the radius of the sphere. We used the following parameters: $a=0.01$ and $R=10$. The exact analytical solution of this problem is $u^{\text{ref}} = 1 / r$.
The latent variable $u$ is represented by the neural network $u_{\bm \theta}$ described above.

Training and evaluating parameters are the same as in the previous example. In our calculations the weights of loss (\ref{eq4}) were taken as $\lambda_{ic} = 10$ and $\lambda_{r} = 10$ (it is the result of empirical choice).

\begin{figure}[t!]
	\centering
	\begin{subfigure}[t]{0.3\textwidth}
		\centering
		\includegraphics[width=\linewidth]{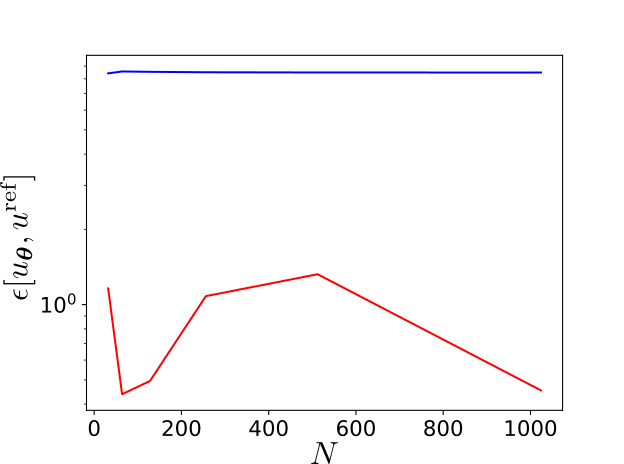}
		\caption{}
	\end{subfigure}
	\begin{subfigure}[t]{0.3\textwidth}
		\centering
		\includegraphics[width=\linewidth]{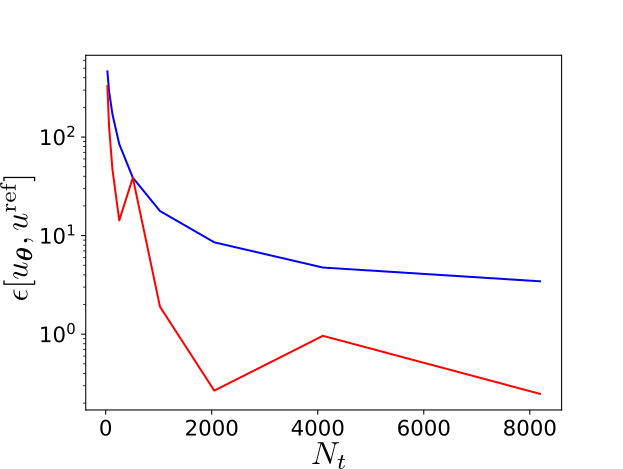}
		\caption{}
	\end{subfigure}
	\caption{Electric Potential. (a) and (b) are dependences of relative ${\mathbb L}_2$ error for $u_{\bm \theta}$ on the number of neurons $N$ and $N_t$, respectively. Blue and red solid lines correspond to two-stage and three-stage training. Number of collocation points is $N_t=2048$ for (a) and number of neurons is $N=128$ for the (b), maximum of time is $R=10$.}
	\label{fig8}
\end{figure}

Figure~\ref{fig8} shows dependences of relative ${\mathbb L}_2$ error for $u_{\bm \theta}$ on the number of neurons $N$ ((a) panel) and number of collocation points $N_t$ ((b) panel) under two-stage and three-stage training, respectively. The figures show that the minimum error value is achieved with 64 neurons for three-stage learning and at 2048 collocation points. The relative error obtained in three-stage learning is much smaller than the error obtained in two-stage learning.

\begin{figure}[t!]
	\centering
	\begin{subfigure}[t]{0.3\textwidth}
		\centering
		\includegraphics[width=\linewidth]{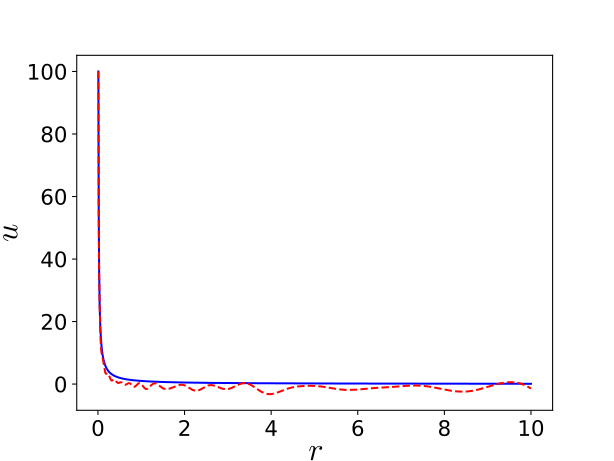}
		\caption{}
	\end{subfigure}
	\begin{subfigure}[t]{0.3\textwidth}
		\centering
		\includegraphics[width=\linewidth]{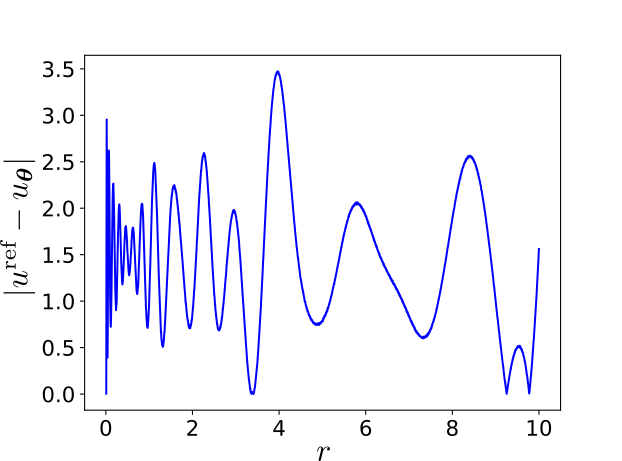}
		\caption{}
	\end{subfigure}
	\caption{Electric Potential. (a) is a comparison of the predicted (red dash lines) and reference solutions (blue solid lines) corresponding to $u$. (b) is absolute errors $|u^{\text{ref}} - u_{\bm \theta}|$. Number of neurons is $N=128$, the number of collocation points $N_t = 2048$, maximum of radius is $R=10$.}
	\label{fig9}
\end{figure}
Results of three-stage training for a neural network with $N=128$, number of collocation points $N_t=2048$ and $R=10$ are shown on Figure~\ref{fig9}. The relative ${\mathbb L}_2$ errors is $\epsilon[u_{\bm \theta}, {u^{\text{ref}}}] = 4.81{\times} 10^{-1}$. The training was carried out on CPU and took $978$ iterations of LBFGS optimizer and time $8.377$ seconds.

\subsubsection{Discussion}
The PINN solutions of the presented problems demonstrate non-monotonic dependences of accuracy on the number of neurons of the hidden layer and on the number of collocation points in the training set. This behaviour on the one hand is explained by overlapping areas of the greatest change in the output of neurons of the hidden layer due to the same values of the weights on the hidden layer, equal to 1. On the other hand, the weights of the output layer are chosen randomly. Let's attempt to address both of these drawbacks.

We determine the weights of the hidden layer as ${W}^{(1)}_k = 2 \Delta \zeta / \Delta x$. With such a choice of these weights in the range $[-\Delta \zeta,\Delta \zeta]$ the derivative of the sigmoid of the neuron must change slightly and the sigmoid must approximate by a straight line in this interval. The choice of $\Delta\zeta$ allows for some arbitrariness (in our experiments, the best results were achieved at $\Delta\zeta=0.7$). These weights act as scaling coordinate factors for $x$.  With this set of weights, the biases of the hidden layer must be converted to $b^{(1)}_k=-2(x_k+\Delta x) \Delta \zeta / \Delta x$ (it means $b^{(1)}_k=- 2 (k+1) \Delta \zeta$). In this case of weights and biases of the hidden layer, the output of $k$th neuron of the hidden layer is proportional to $\sigma\left(\dfrac{2 \Delta \zeta}{\Delta x} \left[x -(x_k+\Delta x)\right]\right)$. Therefore, such representation of biases preserves the biases relative to the coordinates $x$ for the given $k$th neuron.

Find the derivative of the component of solution $u_{\bm \theta;l}(x)$  at the point with the coordinate $x_m$ ($x_m = m \Delta x$ is the boundary of the cell for which $m$ neurons are responsible):
\begin{equation}\label{eq14}
	\left.\frac{\partial u_{\bm \theta;l}(x)}{\partial x}\right\vert_{x=x_m} = \sum\limits_{k = 0}^{N-1} \frac{2\Delta \zeta}{\Delta x}{W}^{(2)}_k\sigma'\left(\frac{2\Delta \zeta}{\Delta x} \left[x_m -(x_k+\Delta x)\right]\right).
\end{equation}
Here, $\sigma'(x)$ means a derivate on the argument of function $\sigma(x)$. Pay attention $\sigma'(x)$ has maximum at $x=0$, $\sigma'(0) = 0.25$ and $x_{k+1} -(x_k+\Delta x) = 0$.We assume that the main contribution to the sum on the right-hand side of Eqn.~(\ref{eq14})  is provided by a term containing $\sigma'(0)$. As result we have following approximation of derivative of $u_{\bm \theta;l}(x)$
\begin{equation}\label{eq15}
	\left.\frac{\partial u_{\bm \theta;l}(x)}{\partial x}\right\vert_{x=x_m} \approx \frac{\Delta \zeta}{2 \Delta x}{W}^{(2)}_{m-1}.
\end{equation}
On the other hand, derivative of $u_{\bm \theta;l}(x)$ must satisfy the Equations (\ref{eq1}). Using (\ref{eq1}) and (\ref{eq15} we have the following equation for the weights ${W}^{(2)}_{m}$
\begin{equation}\label{eq16}
	{W}^{(2)}_{m} =  - 2 \Delta x \mathcal{N}_l [{\vec u}, x_{m+1}] / \Delta \zeta.
\end{equation}
Here values ${\vec u}$ are taken from the initial conditions (\ref{eq2}), $x_{N} = x_{N-1} + \Delta x$ and $\mathcal{N}_l [{\vec u}, x_{m+1}]$ is the $l$ component of the result of operator action $\mathcal{N} [{\vec u}, x_{m+1}]$.

\begin{figure}[th!]
	\centering
	\includegraphics[width=0.9\textwidth]{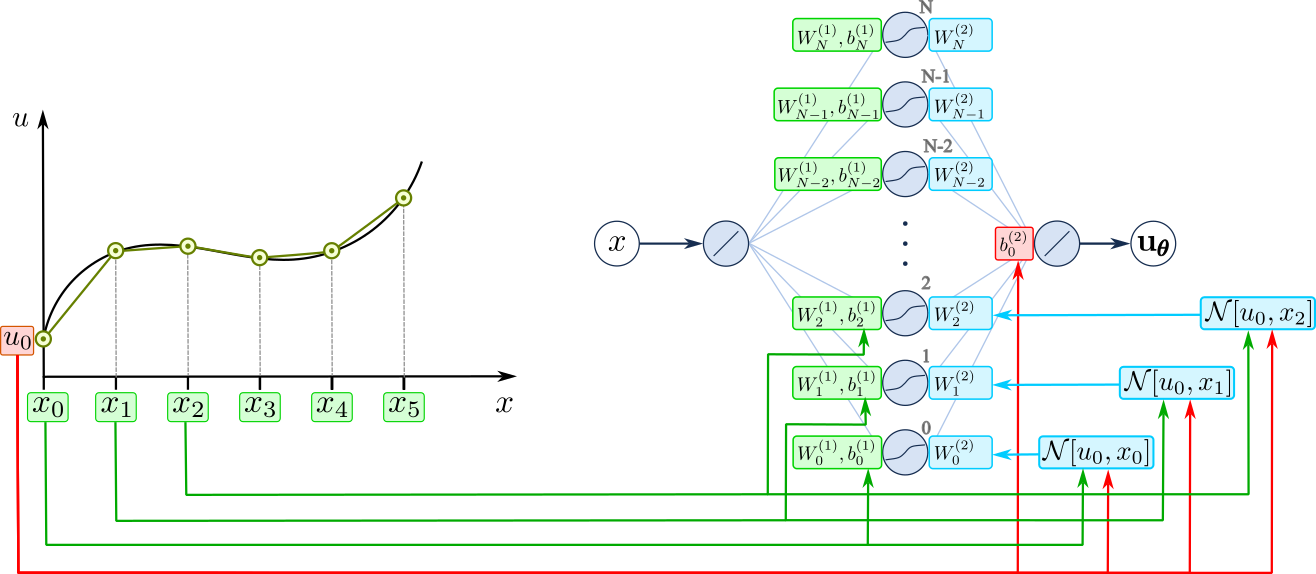}
	\caption{Explanatory figure for strictly deterministic initialization of PINN (Algorithm~\ref{alg2}).}
	\label{fig10}
\end{figure}

{\bf Suggestion 3: Strictly deterministic initialization.} Based on the above reasoning we initialize the parameters of neural network $u_{\bm \theta;l}(x)$ as follows (see the explanatory  Figure~\ref{fig10}):
\begin{enumerate}
	\item Weights of the hidden layer are taken $2 \Delta \zeta / \Delta x$ (${W}^{(1)}_k = 2 \Delta \zeta / \Delta x$) and biases $b^{(1)}_k=-2 (k+1)\Delta \zeta$ for the $k$th neuron.
	\item The bias of the output layer is taken equal to the initial value $u_0$ ($b^{(2)}_0=u_0$).
	\item The parameters of output layer are initialized with ${W}^{(2)}_{k} =  - 2 \Delta x \mathcal{N}_l [{\vec u}(0), x_{k+1}] / \Delta \zeta$.
\end{enumerate}
As a result, we have the following neural network initialization Algorithm~\ref{alg2}.
\begin{algorithm}[h!]
	\caption{Strictly deterministic initialization of physics-informed neural network}\label{alg2}
	\KwData{---}
	\KwResult{Initialized neural network ${u}_{{\bm \theta};l}(x)$ of PINN ${\vec u}_{\bm \theta}$, which consists of $N$ neurons on hidden layer}
	$\Delta x \gets X / N$\;	
	\For{$k=0,\dots, N-1$}{
		${W}^{(1)}_k \gets 2 \Delta \zeta / \Delta x$\;
		${b}^{(1)}_k \gets -2 (k+1) \Delta \zeta$\;
		${W}^{(2)}_k \gets - 2 \Delta x \mathcal{N}_l [{\vec u}(0), x_{k+1}] / \Delta \zeta$\;}
	${b}^{(2)}_0 \gets {u}_l(0)$.
\end{algorithm}

\subsection{Numerical experiments}
In the examples  given in this subsection, we investigate the dependences of the accuracy of predicting solutions when using the neural network initialization Algorithm~\ref{alg2} on neural network parameters and training parameters.

\subsubsection{Example 3. Harmonic Oscillator}

In our experiments, we used only three-stage learning, which is described above, because it previously provided a more accurate solution than the two-stage one. In our calculations the weights of loss (\ref{eq4}) were taken as $\lambda_{ic} =10^3$ and $\lambda_{r} = 10 \lambda_{ic}$. Note that a large absolute value of $\lambda_{ic}$ contributes to improving training outcomes with the help of an LBFGS optimizer. We used 5000  uniform distributed collocation points for calculating the relative error trained model.

\begin{figure}[t!]
	\centering
	\begin{subfigure}[t]{0.3\textwidth}
		\centering
		\includegraphics[width=\linewidth]{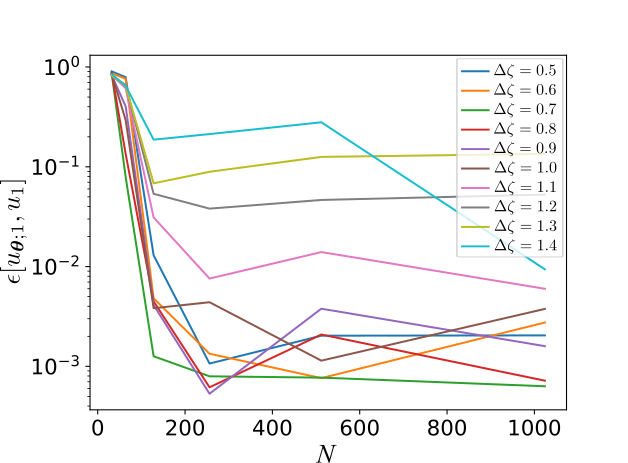}
		\caption{$u_1$}
	\end{subfigure}
	\begin{subfigure}[t]{0.3\textwidth}
		\centering
		\includegraphics[width=\linewidth]{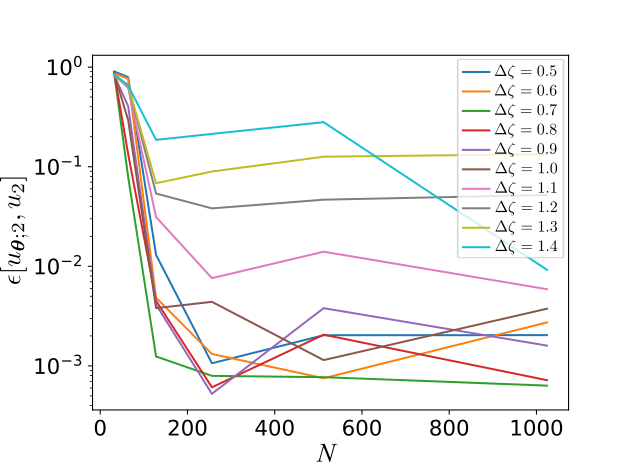}
		\caption{$u_2$}
	\end{subfigure}
	\caption{Harmonic Oscillator. (a) and (b) are dependences of relative ${\mathbb L}_2$ error for $u_{\bm \theta; 1}$ and ${u}_{\bm \theta;2}$, respectively, on the number of neurons $N$ for the different value of $\Delta \zeta$. Number of collocation points is $N_t=2048$, maximum of time is $T=100$.}
	\label{fig11}
\end{figure}

Figure~\ref{fig11} shows dependences of relative ${\mathbb L}_2$ error for $u_{\bm \theta; 1}$ and ${u}_{\bm \theta;2}$, respectively, on the number of neurons $N$ under three-stage training for the different value of $\Delta \zeta$. The figures show that the minimum stable errors are achieved for the $\Delta \zeta = 0.7$.

\begin{figure}[t!]
	\centering
	\begin{subfigure}[t]{0.3\textwidth}
		\centering
		\includegraphics[width=\linewidth]{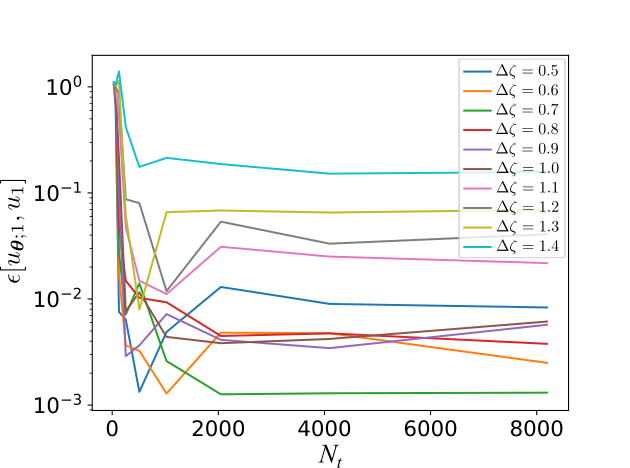}
		\caption{$u_1$}
	\end{subfigure}
	\begin{subfigure}[t]{0.3\textwidth}
		\centering
		\includegraphics[width=\linewidth]{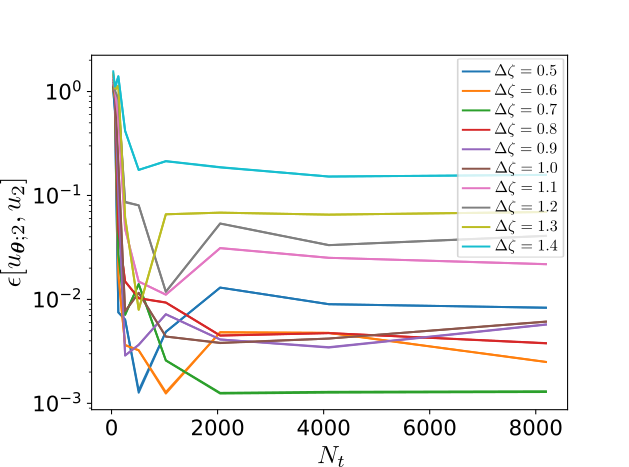}
		\caption{$u_2$}
	\end{subfigure}
	\caption{Harmonic Oscillator. (a) and (b) are dependences of relative ${\mathbb L}_2$ error for $u_{\bm \theta; 1}$ and ${u}_{\bm \theta;2}$, respectively, on the number of collocation points $N_t$ for the different value of $\Delta \zeta$. Number of neurons is $N=128$, maximum of time is $T=100$.}
	\label{fig12}
\end{figure}

Figure~\ref{fig12} shows dependences of relative ${\mathbb L}_2$ error for $u_{\bm \theta; 1}$ and ${u}_{\bm \theta;2}$, respectively, on the number of collocation points $N_t$ under three-stage training for the different value of $\Delta \zeta$. The figures show that the minimum stable errors are achieved for the $\Delta \zeta = 0.7$.

\begin{figure}[t!]
	\centering
	\begin{subfigure}[t]{0.3\textwidth}
		\centering
		\includegraphics[width=\linewidth]{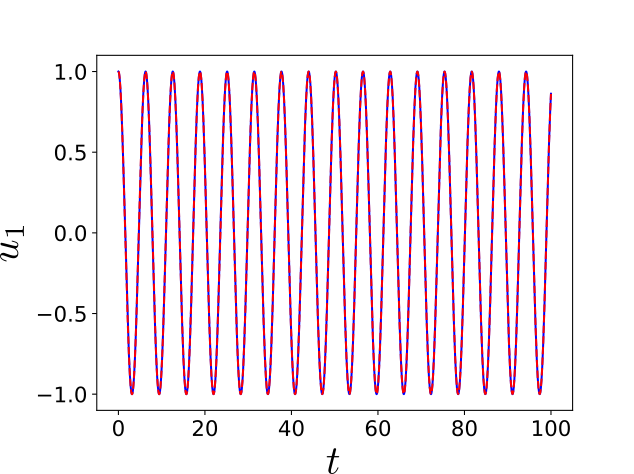}
		\caption{$u_1$}
	\end{subfigure}
	\begin{subfigure}[t]{0.3\textwidth}
		\centering
		\includegraphics[width=\linewidth]{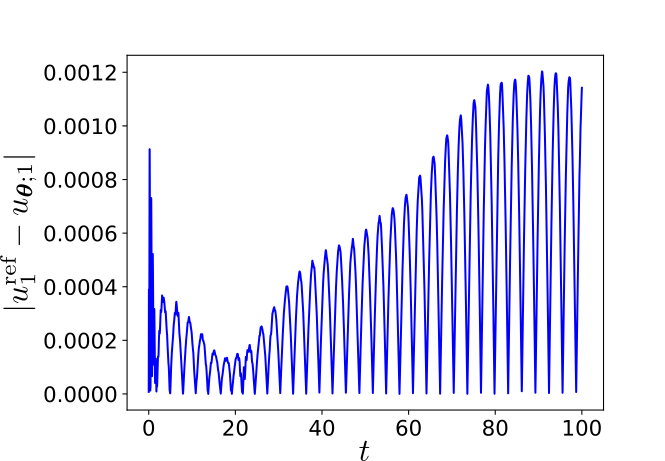}
		\caption{$u_2$}
	\end{subfigure}
	\begin{subfigure}[t]{0.3\textwidth}
		\centering
		\includegraphics[width=\linewidth]{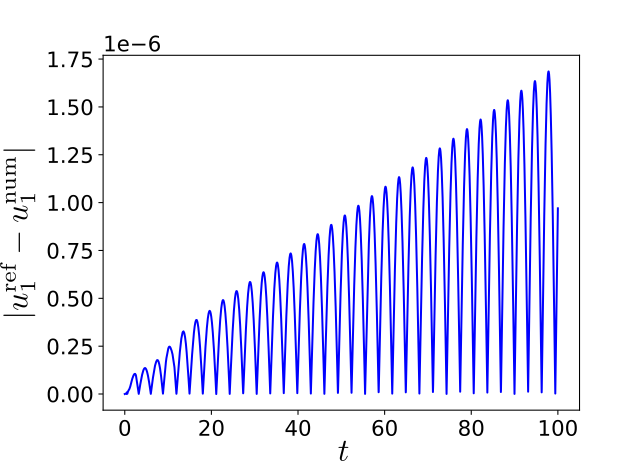}
		\caption{$u_1$}
	\end{subfigure}
	\begin{subfigure}[t]{0.3\textwidth}
		\centering
		\includegraphics[width=\linewidth]{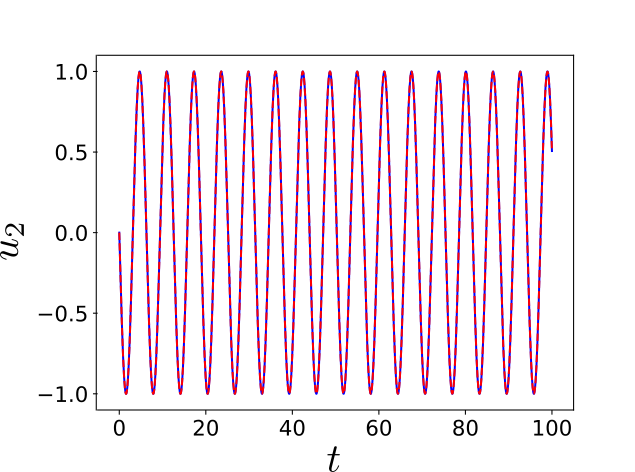}
		\caption{$u_2$}
	\end{subfigure}
	\begin{subfigure}[t]{0.3\textwidth}
		\centering
		\includegraphics[width=\linewidth]{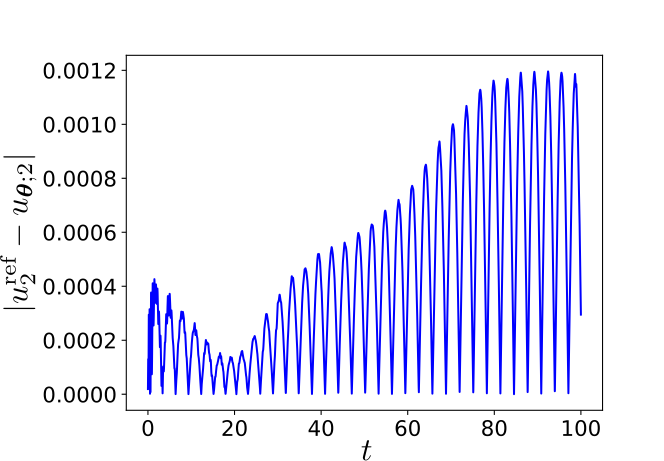}
		\caption{$u_1$}
	\end{subfigure}
	\begin{subfigure}[t]{0.3\textwidth}
		\centering
		\includegraphics[width=\linewidth]{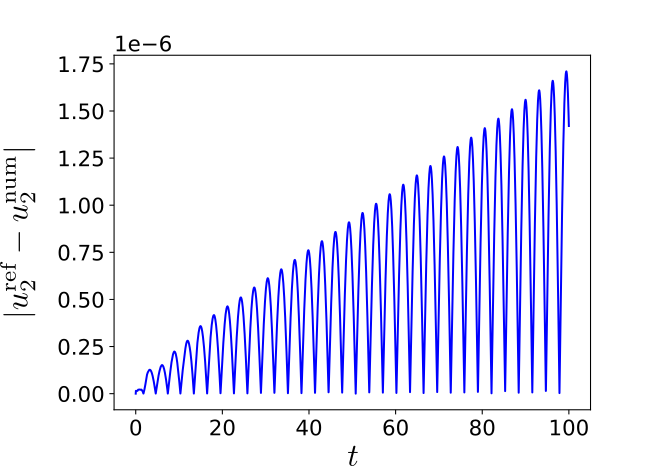}
		\caption{$u_2$}
	\end{subfigure}
	\caption{Harmonic Oscillator. (a) and (d) are comparisons of the predicted (red dash lines) and reference solutions (blue solid lines) corresponding to $u_1$ and $u_2$, respectively. (b) and (e) are absolute errors $|u_1^{\text{ref}} - u_{\bm \theta; 1}|$ and $|u_2^{\text{ref}} - {u}_{\bm \theta;2}|$, respectively. Number of neurons is $N=512$, the number of collocation points $N_t = 2048$, maximum of time is $T=100$, $\Delta \zeta = 0.7$. Panels (c) and (f) correspond to the panels (b) and (e) were obtained by numeric ODE solver ($|u_1^{\text{ref}} - u_1^{\text{num}}|$ and $|u_2^{\text{ref}} - u_1^{\text{num}}|$).}
	\label{fig13}
\end{figure}
Results of three-stage training for a neural network with $N=512$, number of collocation points $N_t=2048$, $T=100$, and $\Delta \zeta = 0.7$ are shown on Figure~\ref{fig13}. The relative ${\mathbb L}_2$ errors are $\epsilon[{u}_{\bm \theta;1}, {u_1^{\text{ref}}}] = 7.71{\times} 10^{-4}$ and $\epsilon[{u}_{\bm \theta;2}, {u_2^{\text{ref}}}] = 7.71{\times} 10^{-4}$. The training was carried out on CPU and took 8751 iterations of LBFGS optimizer and time $200.64$ seconds. It can be seen, the accuracy of the solution given by the neural network initialized by Algorithm~\ref{alg2} has increased by an order of magnitude compared to the initialization Algorithm~\ref{alg1} (compare with the results of Example~1).

\subsubsection{Example 4. The electric potential of the charged sphere}
In this our experiments, we used three-stage learning. In out calculations the weights of loss (\ref{eq4}) were taken as $\lambda_{ic} = 10$ and $\lambda_{r} = 10$.

\begin{figure}[t!]
	\centering
	\begin{subfigure}[t]{0.3\textwidth}
		\centering
		\includegraphics[width=\linewidth]{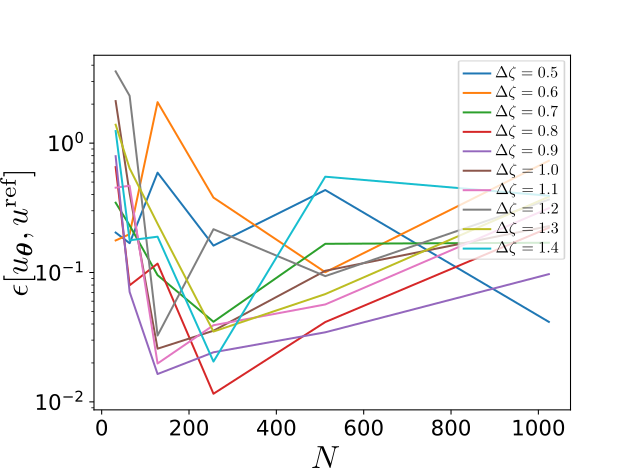}
		\caption{}
	\end{subfigure}
	\begin{subfigure}[t]{0.3\textwidth}
		\centering
		\includegraphics[width=\linewidth]{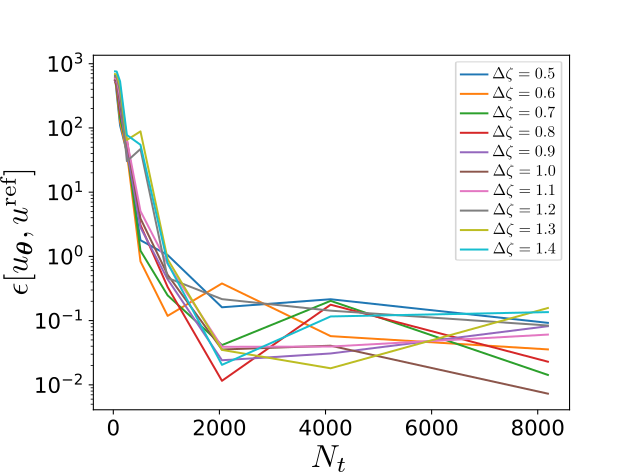}
		\caption{}
	\end{subfigure}
	\caption{Electric Potential. (a) and (b) are dependences of relative ${\mathbb L}_2$ error for $u_{\bm \theta}$ on the number of neurons $N$ and $N_t$, respectively, for the different value of $\Delta \zeta$. Number of collocation points is $N_t=2048$ for (a) and number of neurons is $N=128$ for the (b), maximum of distance is $R=10$.}
	\label{fig14}
\end{figure}
Figure~\ref{fig14} shows dependences of relative ${\mathbb L}_2$ error for $u_{\bm \theta}$ on the number of neurons $N$ ((a) panel) and number of collocation points $N_t$ ((b) panel) under three-stage training for the different value of $\Delta \zeta$.

\begin{figure}[t!]
	\centering
	\begin{subfigure}[t]{0.3\textwidth}
		\centering
		\includegraphics[width=\linewidth]{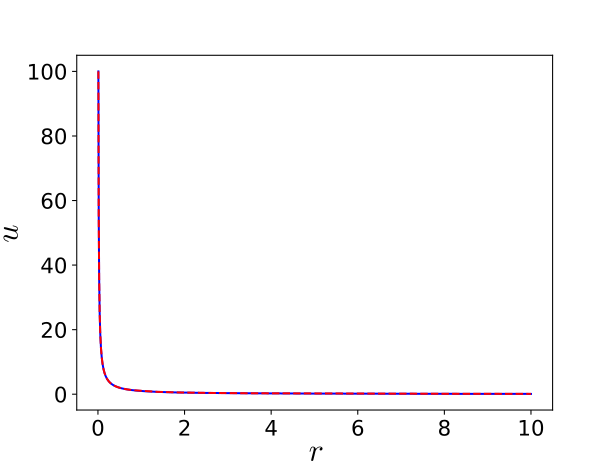}
		\caption{}
	\end{subfigure}
	\begin{subfigure}[t]{0.3\textwidth}
		\centering
		\includegraphics[width=\linewidth]{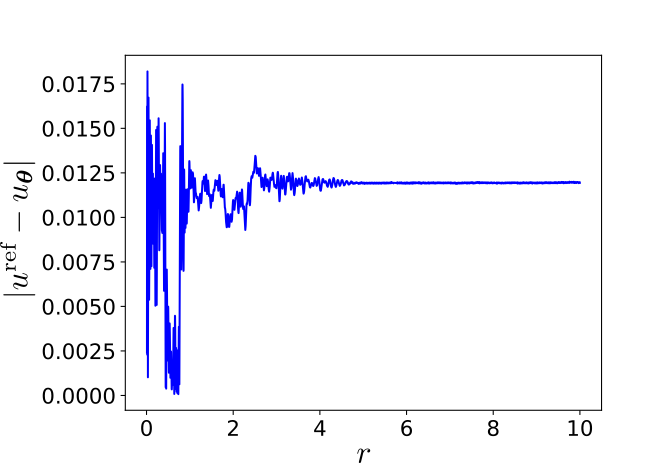}
		\caption{}
	\end{subfigure}
	\begin{subfigure}[t]{0.3\textwidth}
		\centering
		\includegraphics[width=\linewidth]{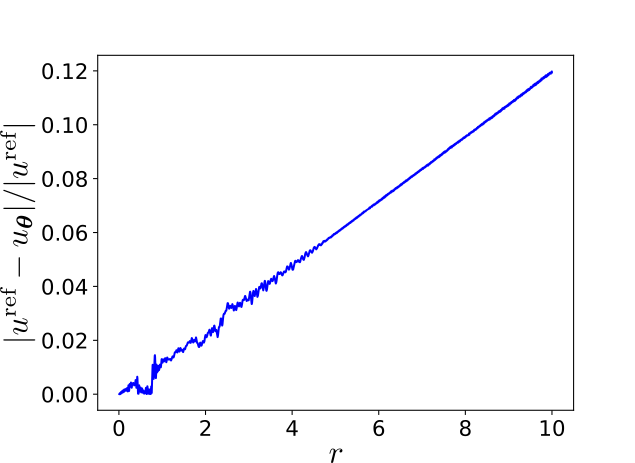}
		\caption{}
	\end{subfigure}
	\begin{subfigure}[t]{0.3\textwidth}
		\centering
		\includegraphics[width=\linewidth]{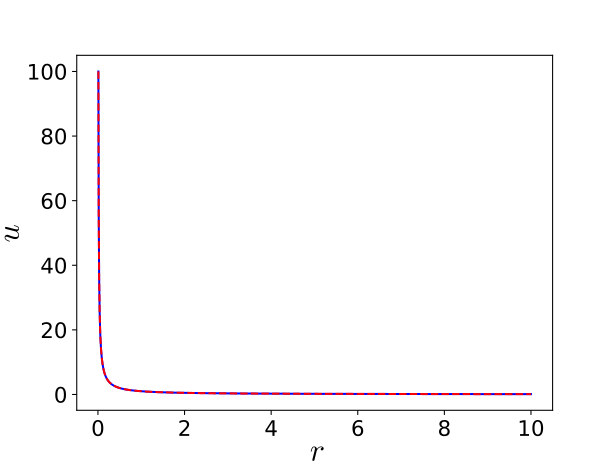}
		\caption{}
	\end{subfigure}
	\begin{subfigure}[t]{0.3\textwidth}
		\centering
		\includegraphics[width=\linewidth]{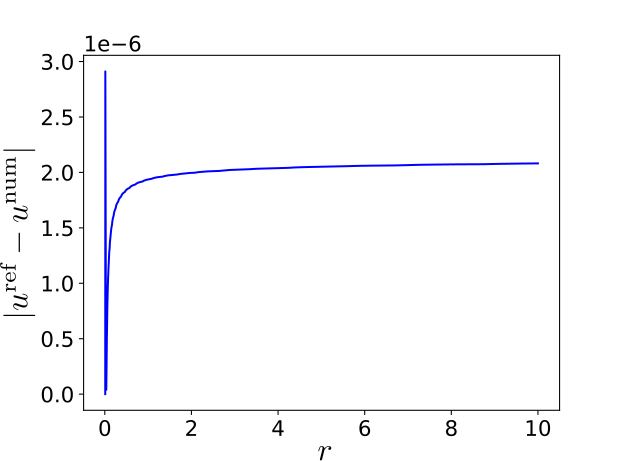}
		\caption{}
	\end{subfigure}
	\begin{subfigure}[t]{0.3\textwidth}
		\centering
		\includegraphics[width=\linewidth]{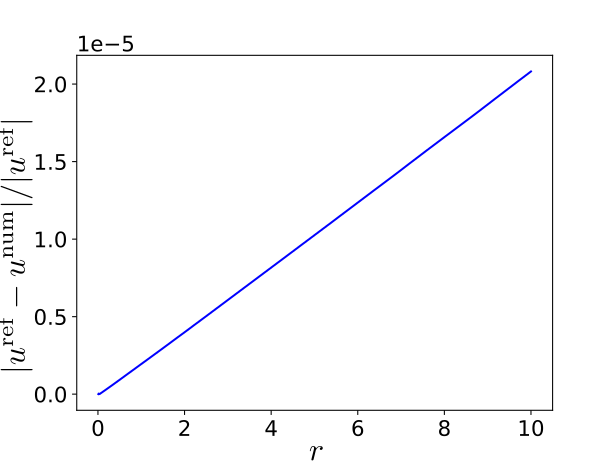}
		\caption{}
	\end{subfigure}
	\caption{Electric Potential. (a) is a comparison of the predicted (red dash lines) and reference solutions (blue solid lines) corresponding to $u$. (b) and (c) are absolute errors $|u^{\text{ref}} - u_{\bm \theta}|$ and relative error $|u^{\text{ref}} - u_{\bm \theta}|/|u^{\text{ref}}|$. Panels (d), (e), (f) correspond to the panels (a), (b), (c) were obtained by ODE solver. Number of neurons is $N=512$, the number of collocation points $N_t = 8000$, maximum of radius is $R=10$.}
	\label{fig15}
\end{figure}
Results of three-stage training for neural network with $N=512$, number of collocation points $N_t=8000$, $R=10$ and $\Delta \zeta=1$ are shown on the Figure~\ref{fig15}(a)--(c). The relative ${\mathbb L}_2$ errors is $\epsilon[u_{\bm \theta}, {u^{\text{ref}}}] = 3.49{\times} 10^{-3}$. The training was carried out on CPU and took $1903$ iterations of LBFGS optimizer and time $108$ seconds.  Panels (d), (e), (f) correspond to the panels (a), (b), (c) were obtained by \verb*|odeint| solver of \verb*|scipy.integrate| library. It is noteworthy that the dependences of absolute error and relative error on the argument $r$ of the solution provided by the neural network qualitatively repeat similar dependences for a deterministic solver. 


\subsubsection{Example 5. Relativistic slingshot}\label{ex5}

Consider the problem of a source for single circularly polarized attosecond x-ray pulses~\cite{Wang2013}. The problem is formulated by following a system of ordinary differential equations
\begin{align}
	& \frac{dh}{d\xi} = E_x - \varepsilon \frac{u_{\perp}^{2}}{1+u_{\perp}^{2}},   \notag\\
	&\frac{dx}{d\xi} = \frac{1+u_{\perp}^{2}-h^{2}}{2 h^{2}},\notag\\
	& \frac{dy}{d\xi} = \frac{u_y}{h},\notag\\
	& \frac{dz}{d\xi} = \frac{u_z}{h},\notag\\
	& u_{y} = a_{y,L} - \varepsilon y,\notag\\
	& u_{z}=a_{z,L}-\varepsilon z,\label{eq17}		
\end{align}
where $h=\gamma - u_{x}$, $\gamma^{2} = 1 + u_{\perp}^{2} + u_{x}^{2}$ is the relativistic gamma factor, $u_{x,y,z}$ are the space components of the four-velocity, $u_{\perp}^{2} = u_{y}^{2}+ u_{z}^{2}$, $\xi = t - x$, $\varepsilon = 2 \pi e^{2}n'l'/m\omega_{L}c$ ($l'$ is foil thickness), $a_{y,L}$ and $a_{z,L}$ are the $y$ and $z$ components of the laser vector potential ${\vec A}$. It is assumed that the electron is initially at $x=0$, $y=0$ and $z=0$, $E_{x} = \varepsilon {\rm th}\left({x}/ {4l'}\right)$. For the calculations we take $\varepsilon=4\pi$,  $l'=0.01 \lambda_{L}$ ($\lambda_{L} = 2 \pi {c}/{\omega_{L}}$).The period of impulse is $T_L = 1$ ($\omega_{L} = 2\pi / T_{L}$) and duration of impulse is  $T = 4 T_{L}$. In our calculation the value $l'$ was given $l'=0.01 2 \pi$. The amplitude of the vector potential components of the initiating pulse is described by
\begin{align}
	& a_{y} = a_0 \sin\left( 2 \pi t / T_L \right) \sin^2 \left( \pi t / T \right),\quad t \in [0,T],   \notag\\
	& a_{z} = a_{0} \cos\left( 2 \pi {t}/{T_{L}} \right) \sin^2 \left( \pi t / T \right), \quad t\in [0,T].
\end{align}

Rewrite the equations (\ref{eq17}) in the absolute coordinate $t$. First note
\begin{equation}
	{\frac{\partial }{\partial \xi} } = \frac{\partial }{\partial t}\frac{\partial t}{\partial \xi} = \left( 1-\frac{\partial x}{\partial t}  \right)^{-1}\textit{}\frac{\partial }{\partial t} .      \notag
\end{equation}
From second equation of the (\ref{eq17}) we have
\begin{equation}
	\frac{\partial x}{\partial t} = \frac{1 + u_{\perp}^{2} - h^{2}}{2 h^{2}}\left( 1-\frac{\partial x}{\partial t}  \right),    \notag
\end{equation}
\begin{equation}
	\frac{\partial x}{\partial t}  = \frac{b}{1 + b}, \quad b = \frac{1 + u_{\perp}^{2} - h^{2}}{2 h^{2}}.   \notag
\end{equation}
As a results we have equations
\begin{align}
	& \frac{\partial h}{\partial t} = \left(E_x - \varepsilon \frac{u_{\perp}^{2}}{1+u_{\perp}^{2}}\right) \frac{1}{1 + b},   \notag\\
	&\frac{\partial x}{\partial t} = \frac{b}{1 + b},\notag\\
	& \frac{\partial y}{\partial t} = \frac{u_y}{h (1+b)},\notag\\
	& \frac{\partial z}{\partial t} = \frac{u_z}{h (1+b)}.\label{eq19}		
\end{align}
At the $t=0$ the value $h=1$, because $u_{x,y,z} =0$, and $\gamma = 1$.
Normalize the equations with replacing $\tilde{h} = h / \varepsilon$, $\tilde{y} = \varepsilon y$ and $\tilde{z} = \varepsilon z$:
\begin{align}
	& \frac{\partial \tilde{h}}{\partial t} = \left(\tilde{E}_x - \frac{u_{\perp}^{2}}{1+u_{\perp}^{2}}\right) \frac{1}{1 + b},   \notag\\
	&\frac{\partial x}{\partial t} = \frac{b}{1 + b},\notag\\
	& \frac{\partial \tilde{y}}{\partial t} = \frac{u_y}{\tilde{h} (1+b)},\notag\\
	& \frac{\partial \tilde{z}}{\partial t} = \frac{u_z}{\tilde{h} (1+b)}.\label{eq20}		
\end{align}
Here $\tilde{E}_{x} = {\rm th}\left({x}/ {4l'}\right)$, $u_{y} = a_{y,L} - \tilde{y}$, $u_{z}=a_{z,L}- \tilde{z}$, $b =({1 + u_{\perp}^{2} - (\varepsilon \tilde{h})^{2}})/({2 (\varepsilon \tilde{h})^{2}})$.

\begin{figure}[t!]
	\centering
	\begin{subfigure}[t]{0.24\textwidth}
		\centering
		\includegraphics[width=\linewidth]{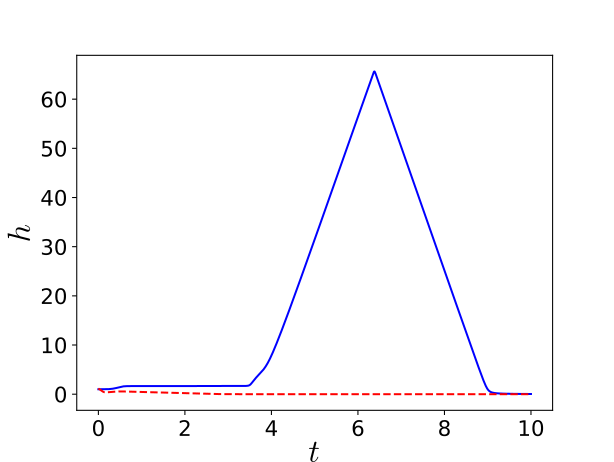}
		\caption{}
	\end{subfigure}
	\begin{subfigure}[t]{0.24\textwidth}
		\centering
		\includegraphics[width=\linewidth]{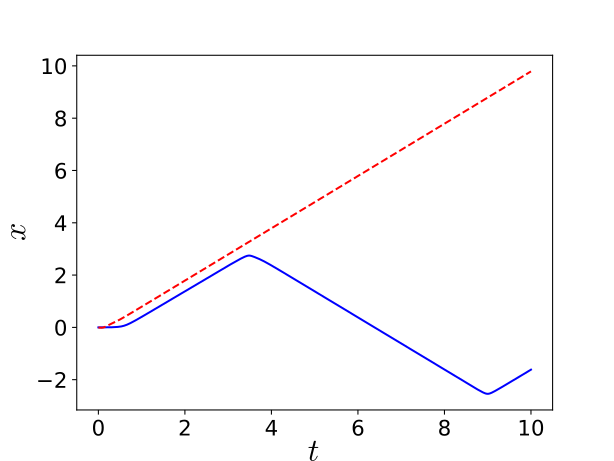}
		\caption{}
	\end{subfigure}
	\begin{subfigure}[t]{0.24\textwidth}
		\centering
		\includegraphics[width=\linewidth]{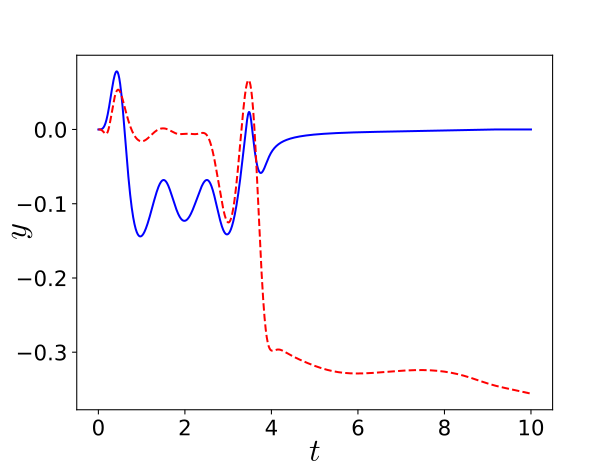}
		\caption{}
	\end{subfigure}
	\begin{subfigure}[t]{0.24\textwidth}
		\centering
		\includegraphics[width=\linewidth]{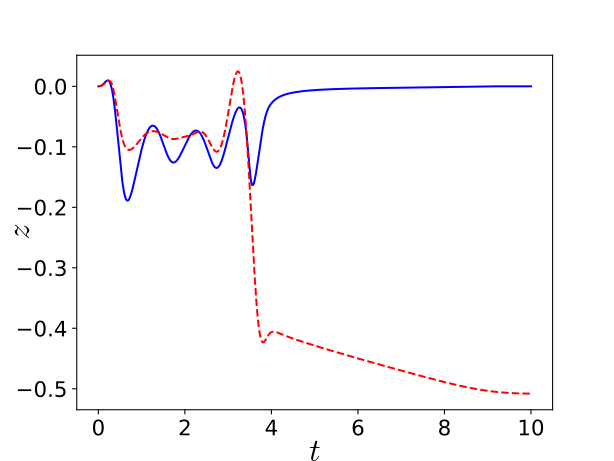}
		\caption{}
	\end{subfigure}
	\caption{Relativistic slingshot. (a), (b), (c) and (d) are comparisons of the predicted (red dash lines) and reference solutions (blue solid lines) corresponding to $h$, $x$, $y$ and $z$, respectively. Number of neurons is $N=128$, the number of collocation points $N_t = 2048$, maximum of time is $T=10$.}
	\label{fig16}
\end{figure}

Results of three-stage training for neural networks with $N=128$, number of uniformly distributed collocation points $N_t=2048$, $T=10$ and $\Delta \zeta = 0.7$ are shown on Figure~\ref{fig16}. The reference solutions are obtained by using the \verb*|odeint| solver of \verb*|scipy.integrate| library. The relative ${\mathbb L}_2$ errors are $\epsilon[{h}_{\bm \theta}, {h^{\text{ref}}}] = 1.00$, $\epsilon[{x}_{\bm \theta}, {x^{\text{ref}}}] = 3.8$, $\epsilon[{y}_{\bm \theta}, {y^{\text{ref}}}] = 4.5$, and $\epsilon[{z}_{\bm \theta}, {z^{\text{ref}}}] = 5.65$. The training was carried out on CPU and took 3928 iterations of LBFGS optimizer and time is $82$ seconds.

\subsubsection{Discussion}
As has been demonstrated in a number of examples, neural networks with the proposed strictly deterministic initialization of weights and biases, described by algorithm~\ref{alg2}, provide more accurate solutions than neural networks with random initialization of output layer parameters, described by algorithm~\ref{alg1}. In the case of Example 5, such initialization of weights alone is not enough to obtain a high-precision solution to this nonlinear problem. In order to improve the accuracy of solutions provided by these strictly deterministic initialized neural networks, we modify the training procedure by using both best-established approaches and methods that are not covered in the literature on this subject.

\newpage
\section{Training of strictly deterministic initialized neural networks for solving ODE}
{\bf Suggestion 4:} To improve the accuracy of the ODE solution by the neural network, we use the modifications of the training methods listed below.

\subsection{Detaching the right-hand sides of the equations from the calculation graph}
Take a look at the Euler method again: based on the value of the function $u(x_i)$ at this point $x_i$, the derivative $\dfrac{\partial u}{\partial x}(x_i)$ is calculated, based on which the value of the function at the new point $u(x_{i + 1})$ is obtained. It can be seen, the value of the derivative $\dfrac{\partial u}{\partial x}(x_i)$, calculated by the function at a given point $x_i$, is an isolated and independent value from the value of the function $u(x_{i + 1})$ at the next point $x_{i + 1}$. When training a neural network based on the equation (\ref{eq1}), we actually teach it to correctly predict the {\it derivative} of $u_{\bm \theta} (x)$ at a given point $x$. For more effective training, we do the same as in the Euler method --- isolate the value of the derivative of the function given by the neural network from the value of the function:
\begin{equation}\label{eq21}
	\frac{\partial{\vec u}_{{\bm \theta}}}{\partial x} + \rasymbol{\mathcal{N}} [{\vec u_{\bm \theta}}, x] = 0.
\end{equation}
Here and below Ra symbol over values means calculated values (${\mathcal{N}} [{\vec u}, x]$) are detached from the calculation graph. This operation corresponds to \verb*|detach| in PyTorch, and \verb*|lax.stop_gradient| in JAX. This approach can play a significant role in accelerating the convergence of learning in the case of the operator ${\mathcal{N}} [{\vec u}, x]$ is a nonlinear operator.

\subsection{Respecting Causality}
Decompose the loss function  (\ref{eq4}) into terms corresponding to different components of the desired function ${\vec u}_{{\bm \theta}}$ in the following form
\begin{equation}\label{eq22}
	\mathcal{L}({\bm \theta}) = \sum\limits_{k = 1}^K \left[\lambda_{ic}^{k} \mathcal{L}_{ic}^{k}({\bm \theta}) + \lambda^{k}_{r} \mathcal{L}^{k}_{r}({\bm \theta})\right],	  
\end{equation}
where
\begin{eqnarray}
	&& \mathcal{L}^{k}_{ic} \left( {\bm \theta} \right)  = \left| {u_{\bm \theta;k}}\left(0\right) - {g}_{k}  \right|^{2},\label{eq23}\\
	&& \mathcal{L}^{k}_{r}\left(\bm \theta \right) = \frac{1}{N_r} \sum_{i=1}^{N_{r}} \mathcal{L}_r^{{(\rm t),k}} \left[\left({x}_{i}^{(r)} \right), {\bm \theta} \right] ,\label{eq24}\\
	&& \mathcal{L}_r^{{(\rm t),k}} \left( x, {\bm \theta} \right) = \left| \mathcal{R}^{k}\left[{\vec u}_{{\bm \theta}} \right] \left({x}\right)\right|^{2} \\
	&& \mathcal{R}^{k}\left[{\vec u} \right]:= \frac{\partial {u}_{k}}{\partial {x}} + \mathcal{N}_k\left[ {\vec u}\left(x \right), x  \right].\label{eq25}
\end{eqnarray}
Here $u_{k}$ and $g_{k}$ are $k$th components of vectors $\vec u$ and $\vec g$, respectively. The original problem described by a differential equation (\ref{eq1}) accompanied with the initial condition (\ref{eq2}) can be reformulated to the problem described by the differential equation only. The procedure of reformulation is given in the~\cite{eskin2023optimal}. According to this method the weights $\lambda_{ic}^{k}$ and $\lambda^{k}_{r}$ are related to each other by means of the relation 
\begin{equation}
	\lambda_{ic}^{k} = {\beta^k} \lambda^{k}_{r},\label{eq26}
\end{equation}
where $\beta^k = \max\left(\dfrac{N_r}{X}, \dfrac{1}{\sigma^{[0,5]}}, \dfrac{1}{\sigma^{[2]}} \right)$, $\sigma^{[0,5]}$ and $\sigma^{[2]}$ are argument of solution at which the initial value $u^k(0)$ of the field decreases and increases twice, respectively.

To impose the causal structure within the optimization process and acceleration of learning, we use Dirac delta function causal training~\cite{eskin2023optimal}
\begin{align}
	& \mathcal{L}_{r}^{k}\left(\bm \theta \right) = \frac{1}{N_r}\sum_{i=1}^{N_{r}} w_i \mathcal{L}_{r}^{{(\rm t),k}} \left(x_{i}, {\bm \theta} \right),\label{eq28} \\
	& w_i = \underset{{k}}{\max} \left\{ \exp \left[- \varepsilon \left({s}_i^k\right)^2  \right]\right\}, \label{eq29}\\
	& {s}_i^k = \frac{1}{i} \sum_{m=1}^{i-1} \rasymbol{\mathcal{L}}_r^{{(\rm t),k}} \left(x_m, {\bm \theta} \right). \label{eq30}
\end{align}
In our calculation the initial value $\varepsilon$ is taken $\varepsilon = 10^{-8}$. The $\varepsilon$ is doubled under the satisfying condition $\underset{i}{\min} (\exp \left[ - \varepsilon s_i  \right]) > \delta_w$, $\delta_w = 0.99$. Note that we have modified the sum (\ref{eq28}) for the long time intervals (in our case is $\rasymbol{\mathcal{L}}_r^{{(\rm t),k}}$ instead $m \rasymbol{\mathcal{L}}_r^{{(\rm t),k}}$ as in paper~\cite{eskin2023optimal}). We mark this method in the tables as {\bf $\delta$-causal}.

\subsection{Gradient Normalization}
In order to balance the convergence rates of different loss terms such that their convergence rates are comparable to one another we use the following modification of Gradient Normalization method. We use two types of gradient normalization. For the first type (marked as {\bf GN1})  additional global weights are calculating as

\begin{eqnarray}
	&& \lambda_{ic}^{1} = \lambda_{ic}^{2} =...=\lambda_{ic}^{K}=\lambda_{ic}, \quad \lambda_{r}^{1} = \lambda_{r}^{2} =...=\lambda_{r}^{K}=\lambda_{r},\notag\\ 
	&& \hat{\lambda}_{ic} = \frac{\sum\limits_{k = 1}^K \left( \norm{ \nabla_{\bm \theta} \mathcal{L}_{ic}^{k}({\bm \theta})} + \norm{ \nabla_{\bm \theta} \mathcal{L}^{k}_{r}({\bm \theta})}\right)}{\sum\limits_{k = 1}^K\norm{ \nabla_{\bm \theta} \mathcal{L}_{ic}^{k}({\bm \theta})}},\notag\\
	&& \hat{\lambda}_{r} = \frac{\sum\limits_{k = 1}^K \left( \norm{ \nabla_{\bm \theta} \mathcal{L}_{ic}^{k}({\bm \theta})} + \norm{ \nabla_{\bm \theta} \mathcal{L}^{k}_{r}({\bm \theta})}\right)}{\sum\limits_{k = 1}^K \norm{ \nabla_{\bm \theta} \mathcal{L}_{r}^{k}({\bm \theta})}},\label{eq31}
\end{eqnarray}
where $\norm{\cdot}$ is the ${\mathbb L}_2$ norm. Then the global weights are updated by using a moving average of the form
\begin{align}
	& {\lambda}_{ic} = \alpha \rasymbol{{\lambda}}_{ic} + (1 - \alpha)\rasymbol{\hat{\lambda}}_{ic}\notag,\\
	& {\lambda}_{r} = \alpha \rasymbol{\lambda}_{r} + (1 - \alpha)\rasymbol{\hat{\lambda}}_{r},\label{eq32}
\end{align}
where the parameter $\alpha$ is equal $0.9$ in our calculation. Given that, our training pipeline integrates a self-adaptive learning rate annealing algorithm, which can automatically balance losses during training. With this choice of weights balance of loss terms during training  are performed automatically and the following relations are satisfied
\begin{eqnarray}
	\hat{\lambda}_{ic} \sum\limits_{k = 1}^K \norm{ \nabla_{\bm \theta} \mathcal{L}_{ic}^{k}({\bm \theta})}=\hat{\lambda}_{r} \sum\limits_{k = 1}^K \norm{ \nabla_{\bm \theta} \mathcal{L}_{r}^{k}({\bm \theta})}=\sum\limits_{k = 1}^K \left( \norm{ \nabla_{\bm \theta} \mathcal{L}_{ic}^{k}({\bm \theta})} + \norm{ \nabla_{\bm \theta} \mathcal{L}^{k}_{r}({\bm \theta})}\right).
\end{eqnarray}

For the second type (marked as {\bf GN2})  additional global weights are calculated as
\begin{eqnarray}
	&& \lambda_{ic}^{k} = {\beta^k} \lambda^{k}_{r},\notag\\
	&& \hat{\lambda}_{r}^{p} = \frac{\sum\limits_{k = 1}^K \left( {\beta^k} \norm{ \nabla_{\bm \theta} \mathcal{L}_{ic}^{k}({\bm \theta})} + \norm{ \nabla_{\bm \theta} \mathcal{L}^{k}_{r}({\bm \theta})}\right)}{2 K \left( {\beta^p} \norm{ \nabla_{\bm \theta} \mathcal{L}_{ic}^{p}({\bm \theta})} + \norm{ \nabla_{\bm \theta} \mathcal{L}^{p}_{r}({\bm \theta})}\right)},\notag\\
	&& {\lambda}^{p}_{r} = \alpha \rasymbol{\lambda}^{p}_{r} + (1 - \alpha)\rasymbol{\hat{\lambda}}^{p}_{r},\label{eq34}
\end{eqnarray}
and the following relations are satisfied
\begin{eqnarray}
	&&\hat{\lambda}^{1}_{r} \left( {\beta^1} \norm{ \nabla_{\bm \theta} \mathcal{L}_{ic}^{1}({\bm \theta})} + \norm{ \nabla_{\bm \theta} \mathcal{L}^{1}_{r}({\bm \theta})}\right)=...=\hat{\lambda}^{K}_{r} \left( {\beta^K} \norm{ \nabla_{\bm \theta} \mathcal{L}_{ic}^{K}({\bm \theta})} + \norm{ \nabla_{\bm \theta} \mathcal{L}^{K}_{r}({\bm \theta})}\right)\notag\\
	&& =\sum\limits_{k = 1}^K \left( {\beta^k}\norm{ \nabla_{\bm \theta} \mathcal{L}_{ic}^{k}({\bm \theta})} + \norm{ \nabla_{\bm \theta} \mathcal{L}^{k}_{r}({\bm \theta})}\right).
\end{eqnarray}

\subsection{Additional weighting based on second derivative}
In a number of problems, the hidden solution has a significant second derivative in a relatively narrow range of coordinate changes. An example of such a problem is Relativistic slingshot for which the dependences of components of solution on time are shown at Figure~{\ref{fig16}}. It can be seen, for some time intervals there are ``breaks'' of dependences. The times of such ``breaks'' are associated with a significant increase in the second derivatives. At these intervals, the neural network is poorly trained to predict solutions of the problem. To overcome this disadvantage, we propose to insert additional weights to the terms of the loss function, which represent the second derivatives of the governing equations. As a result the Equation (\ref{eq28}) has the following view
\begin{align}
	& \mathcal{L}_{r}^{k}\left(\bm \theta \right) = \frac{1}{N_r}\sum_{i=1}^{N_{r}} \left[1 + \rasymbol{\mathcal{M}}_{r}^{{(\rm t),k}} \left(x_{i}, {\bm \theta} \right)\right] w_i \mathcal{L}_{r}^{{(\rm t),k}} \left(x_{i}, {\bm \theta} \right),\label{eq36_2}
\end{align}
where
\begin{align}
	& \mathcal{M}_r^{{(\rm t),k}} \left( x, {\bm \theta} \right) = \left|\frac{\partial}{\partial x} \mathcal{R}^{k}\left[{\vec u}_{{\bm \theta}} \right] \left({x}\right)\right|^{2}. \label{eq36_3}
\end{align}
Note, we should take the value $w_i$ equal to $1$ in the case of using vanilla PINN training instead of Dirac delta function causal training. We mark this method in the tables as $D^2$.

\subsection{Numerical experiments}
We use the strictly deterministic initialization of physics-informed neural network (Algorithm~\ref{alg2}) for all examples in this subsection.

\subsubsection{Example 6. Harmonic Oscillator}\label{Ex6}

In this part of our experiments we used four-stage training which consists of 5000 epochs of optimization by Adam optimizer with learning rate $10^{-3}$ and 5 epochs of optimization by LBFGS optimizer (maximal number of iterations per optimization step is 5000) under frozen parameters of the hidden layer of neural network, then 5000 epochs of optimization by Adam optimizer with learning rate $10^{-6}$ and 5 epochs of optimization by LBFGS optimizer under unfrozen parameters of the hidden layer of neural network without detaching. In fact, we took the previously proven three-stage training procedure and supplemented it with another stage using Adam with a low learning rate. The parameters $\lambda_{ic}$ and  $\lambda_{r}$ are recalculated by Eqns.~(\ref{eq32}) and (\ref{eq34}) every 1000 epochs of the stage of Adam optimizations (in case of switching on gradient normalization approach). The initial values of lambdas were taken $\lambda^{k}_{r} = 10^5$ and $\lambda_{ic}^{k} = {\beta^k} \lambda^{k}_{r}$, ${\beta^k}=\dfrac{N_t}{T}$. Other parameters are $N=512$, number of collocation points $N_t=2048$, $T=100$ and $\Delta \zeta = 0.7$. Results are presented in the Table~{\ref{table1}}. The best accuracy is achieved by using $\delta$-causal approach and gradient normalization GN1.

We noticed that changing the weights $\lambda_{ic}$ and $\lambda_{r}$ can lead to better accuracy of the solution. For $\lambda_{ic} =10^3$ and $\lambda_{r} = 10 \lambda_{ic}$ the best results was archived which are $\epsilon[u_{\bm \theta;1}, {u_1^{\text{ref}}}] = 2.19{\times} 10^{-4}$ and $\epsilon[u_{\bm \theta;2}, {u_2^{\text{ref}}}] = 2.18{\times} 10^{-4}$. 

\begin{table*}[h!]
	\centering
	\begin{tabular}{cc|ccc}
		\toprule
		\multicolumn{2}{c}{\bf{Train Settings}} & \multicolumn{3}{c}{\bf{Performance}}\\
		\midrule
		\bf{$\delta$-causal} & \bf{Grad Norm}  &  \bf{Relative ${\mathbb L}_2$ errors ($u_1$, $u_2$)} & \bf{Run time (s)} & \bf{Iterations of LBFGS} \\
		\midrule
		{\bf --} & {\bf --} & {$9.5{\times} 10^{-4}$, $9.5{\times} 10^{-4}$} & {365} & {12118} \\
		{\bf +} & {\bf --} & {$5.5{\times} 10^{-4}$, $5.4{\times} 10^{-4}$} & 382 & 10710 \\
		{\bf --} & GN1 & {$7.7{\times} 10^{-4}$, $7.7{\times} 10^{-4}$} & 358 & 10536 \\
		{\bf +} & GN1 & {${\bf 4.3{\times} 10^{-4}}$, ${\bf 4.4{\times} 10^{-4}}$} & 479 & 14951 \\
		{\bf --} & GN2 & {$8.1{\times} 10^{-4}$, $8.1{\times} 10^{-4}$} & 382 & 12079 \\
		{\bf +} & GN2 & {$8.6{\times} 10^{-4}$, $8.5{\times} 10^{-4}$} & 491 & 15869 \\
		\bottomrule
	\end{tabular}
	\caption{Harmonic Oscillator. Relative ${\mathbb L}_2$ error, run time and iterations of LBFGS optimizer for illustrating the impact of individual components of the proposed methods of training}
	\label{table1}
\end{table*}

\subsubsection{Example 7. The electric potential of the charged plane}
Consider the electric potential of the charged plane, which is governed by the equation for the $x\in [0, X]$
\begin{eqnarray}
	&& \frac{d u}{d x} = \Omega, \quad x \in [0,X],\quad u(0) = 0,\label{eq26_1}
\end{eqnarray}
where $\Omega$ is surface charge density. We used the following parameters: $\Omega=10$ and $X=100$. The exact analytical solution of this problem is $u^{\text{ref}} = \Omega x$.
The latent variable $u$ is represented by neural networks $u_{\bm \theta}$ described above.

Training and evaluating parameters are the same as the previous example in subsection~\ref{Ex6}. In our calculations the weights of loss (\ref{eq4}) were taken as $\lambda_{r} = 10^5$ and $\lambda_{ic} = {\beta} \lambda_{r}$, ${\beta}=\dfrac{N_t}{X}$. The parameters $\lambda_{ic}$ and  $\lambda_{r}$ are recalculated by Eqns.~(\ref{eq32}) every 1000 epochs of stage of Adam optimizations (in case of switching on gradient normalization approach). Other parameters are $N=512$, a number of collocation points $N_t=2048$, $X=100$, $\Delta \zeta = 1.0$, and detaching is switched off. Results are presented in the Table~{\ref{table2}}. 

\begin{table*}[h!]
	\centering
	\begin{tabular}{cc|ccc}
		\toprule
		\multicolumn{2}{c}{\bf{Train Settings}} & \multicolumn{3}{c}{\bf{Performance}}\\
		\midrule
		\bf{$\delta$-causal} & \bf{Grad Norm}  &  \bf{Relative ${\mathbb L}_2$ errors ($u$)} & \bf{Run time (s)} & \bf{Iterations of LBFGS} \\
		\midrule
		{\bf --} & {\bf --} & {$5.0{\times} 10^{-6}$} & {365} & {2879} \\
		{\bf +} & {\bf --} & {${\bf 3.1{\times} 10^{-6}}$} & 117 & 2795 \\
		{\bf --} & GN1 & {$5.7{\times} 10^{-6}$} & 124 & 4408 \\
		{\bf +} & GN1 & {$3.1{\times} 10^{-5}$} & 127 & 3786 \\
		\bottomrule
	\end{tabular}
	\caption{The electric potential of the charged plane. Relative ${\mathbb L}_2$ error, run time and iterations of LBFGS optimizer for illustrating the impact of individual components of the proposed methods of training}
	\label{table2}
\end{table*}

\begin{figure}[t!]
	\centering
	\begin{subfigure}[t]{0.28\textwidth}
		\centering
		\includegraphics[width=\linewidth]{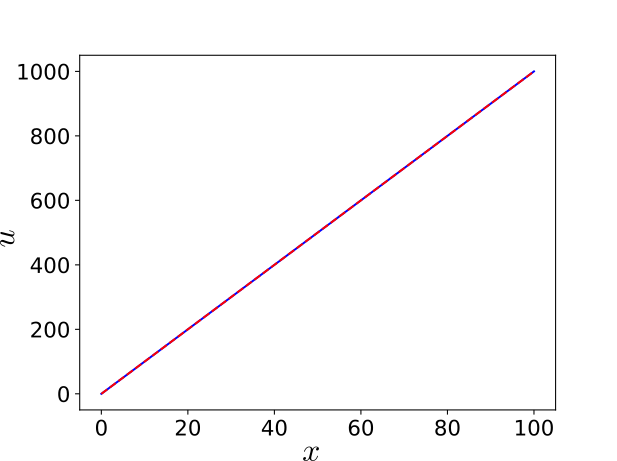}
		\caption{}
	\end{subfigure}
	\begin{subfigure}[t]{0.3\textwidth}
		\centering
		\includegraphics[width=\linewidth]{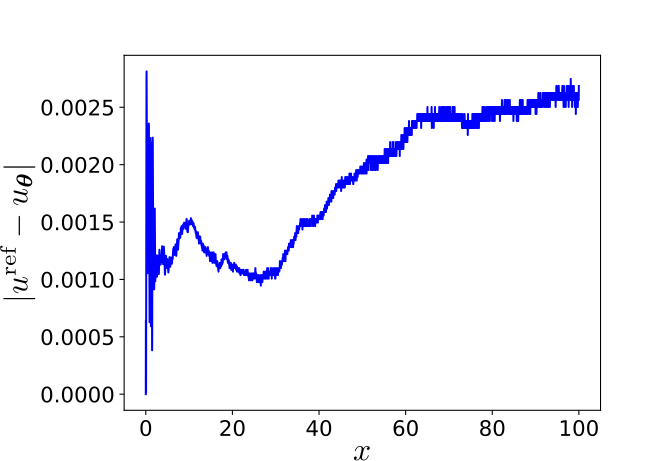}
		\caption{}
	\end{subfigure}
	\begin{subfigure}[t]{0.3\textwidth}
		\centering
		\includegraphics[width=\linewidth]{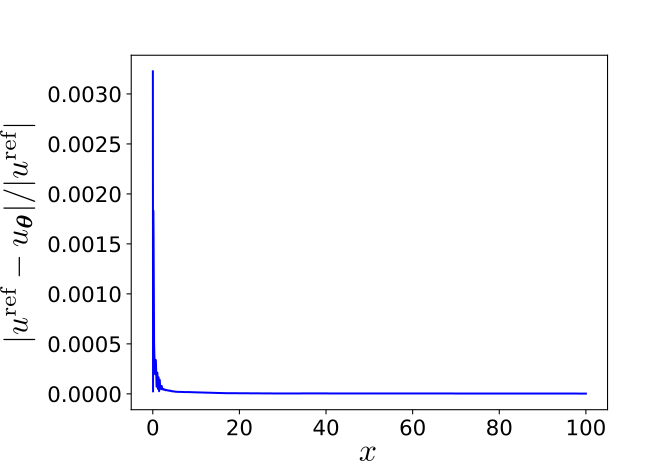}
		\caption{}
	\end{subfigure}
	\caption{Electric potential of the charged plane. (a) is a comparison of the predicted (red dash lines) and reference solutions (blue solid lines) corresponding to $u$. (b) and (c) are absolute errors $|u^{\text{ref}} - u_{\bm \theta}|$ and relative error $|u^{\text{ref}} - u_{\bm \theta}|/|u^{\text{ref}}|$. Number of neurons is $N=512$, the number of collocation points $N_t = 2024$, $X=100$.  The relative ${\mathbb L}_2$ errors is $\epsilon[u_{\bm \theta}, {u^{\text{ref}}}] = 3.38{\times} 10^{-6}$.}
	\label{fig17}
\end{figure}

Results of training for a neural network with $N=512$, number of collocation points $N_t=2048$ and $X=100$ are shown on Figure~\ref{fig17}. The relative ${\mathbb L}_2$ errors is $\epsilon[u_{\bm \theta}, {u^{\text{ref}}}] = 3.38{\times} 10^{-6}$. The training was carried out on CPU and took $3265$ iterations of LBFGS optimizer and time $122$ seconds.

\subsubsection{Example 8. The electric potential of the charged sphere}

In this part of our experiments, we used fifth-stage training. The initial training consists of 1000 epochs of optimization by Adam optimizer with learning rate $10^{-3}$ and 5 epochs of optimization by LBFGS optimizer (maximal number of iterations per optimization step is 5000) under frozen parameters of the hidden layer of neural network, then 5 epochs of optimization by LBFGS optimizer under unfrozen parameters of the hidden layer of neural network without $\delta$-causal and gradient normalization. The next two steps are 1000 epochs of optimization by Adam optimizer with learning rate $10^{-3}$ and 5 epochs of optimization by LBFGS optimizer (maximal number of iterations per optimization step is 5000) under the chosen additional approach. The parameters $\lambda_{ic}$ and  $\lambda_{r}$ are recalculated by Eqns.~(\ref{eq32}) every 500 epochs of stage of Adam optimizations (in case of switching on gradient normalization approach). The initial values of lambdas were taken $\lambda_{r} = 10$ and $\lambda_{ic} = 10$. Other parameters are $N=512$, number of collocation points $N_t=8000$, $R=10$ and $\Delta \zeta = 1.0$.  Results are presented in the Table~{\ref{table3}}.

\begin{table*}[h!]
	\centering
	\begin{tabular}{cc|ccc}
		\toprule
		\multicolumn{2}{c}{\bf{Train Settings}} & \multicolumn{3}{c}{\bf{Performance}}\\
		\midrule
		\bf{$\delta$-causal} & \bf{Grad Norm}  &  \bf{Relative ${\mathbb L}_2$ errors ($u_1$)} & \bf{Run time (s)} & \bf{Iterations of LBFGS} \\
		\midrule
		{\bf --} & {\bf --} & {$5.7{\times} 10^{-4}$} & {141} & {2053} \\
		{\bf +} & {\bf --} & {${\bf 5.5{\times} 10^{-4}}$} & 135 & 2089 \\
		{\bf --} & GN1 & {$5.8{\times} 10^{-4}$} & 131 & 2077 \\
		{\bf +} & GN1 & {$7.2{\times} 10^{-4}$} & 136 & 2085 \\
		\bottomrule
	\end{tabular}
	\caption{The electric potential of the charged sphere: Relative ${\mathbb L}_2$ error, run time and iterations of LBFGS optimizer for illustrating the impact of individual components of the proposed methods of training}
	\label{table3}
\end{table*}

\subsubsection{Example 9. Relativistic slingshot}
In this part of our experiments, we used three-stage training which is consisted of 5000 epochs of optimization by Adam optimizer with a learning rate $10^{-3}$ and 5 epochs of optimization by LBFGS optimizer (maximal number of iterations per optimization step is 5000) under frozen parameters of the hidden layer of neural network, then 5 epochs of optimization by LBFGS optimizer under unfrozen parameters of the hidden layer of neural network. The parameters $\lambda_{ic}$ and  $\lambda_{r}$ are recalculated by Eqns.~(\ref{eq32}) and (\ref{eq34}) every 1000 epochs of the stage of Adam optimization (in case of switching on gradient normalization approach). The initial values of lambdas were taken $\lambda^{k}_{r} = 1$ and $\lambda_{ic}^{k} = {\beta^k} \lambda^{k}_{r}$, ${\beta^k}=\dfrac{N_t}{T}$. Other parameters are $N=512$, number of collocation points $N_t=5048$, $T=10$ and $\Delta \zeta = 0.7$  Results are presented in the Table~{\ref{table4}}.

\begin{table*}[h!]
	\centering
	\begin{tabular}{cccc|ccc}
		\toprule
		\multicolumn{3}{c}{\bf{Train Settings}} & \multicolumn{3}{c}{\bf{Performance}}\\
		\midrule
		\bf{Detaching} & $\bf{D^2}$ & \bf{$\delta$-causal} & \bf{Grad}  &  \bf{Relative ${\mathbb L}_2$ error (h, x, y, z)} & \bf{Run time} & \bf{Iterations} \\
		&  &  & \bf{Norm}  &  & \bf{(s)} & \bf{of LBFGS} \\
		\midrule
		{\bf --} & {\bf --} & {\bf --} & {\bf --} & (1.0, 3.7, 6.2, 9.4) & 1399 & 21775 \\
		{\bf --} & {\bf --} & {\bf +} & {\bf --} & (1.0, 3.7, 6.2, 9.2) & 1584 & 21805 \\
		{\bf --} & {\bf +} & {\bf --} & {\bf --} & (1.0, 3.8, 7.6, 10.0) & 1112 & 149 \\
		{\bf --} & {\bf +} & {\bf +} & {\bf --} & (1.0, 3.8, 7.6, 10.0) & 1165 & 150 \\
		{\bf +} & {\bf --} & {\bf --} & {\bf --} & ($5.6{\times} 10^{-1}$, $3.2{\times} 10^{-1}$, $1.1{\times} 10^{-1}$, $2.4{\times} 10^{-1}$) & 354 & 250 \\
		{\bf +} & {\bf --} & {\bf +} & {\bf --} & ($5.5{\times} 10^{-1}$, $3.2{\times} 10^{-1}$, $1.1{\times} 10^{-1}$, $2.4{\times} 10^{-1}$) & 455 & 715 \\
		{\bf +} & {\bf +} & {\bf --} & {\bf --} &($6.0{\times} 10^{-2}$, $4.7{\times} 10^{-2}$, $1.5{\times} 10^{-1}$, $5.0{\times} 10^{-2}$) & 1135 & 290 \\
		{\bf +} & {\bf +} & {\bf +} & {\bf --} & ($6.2{\times} 10^{-2}$, $5.1{\times} 10^{-2}$, $1.5{\times} 10^{-1}$, $6.1{\times} 10^{-2}$) & 1207 & 549 \\
		{\bf --} & {\bf --} & {\bf --} & GN1 & (1.0, 3.8, 17.4, 19.4) & 1935 & 33476 \\
		{\bf --} & {\bf --} & {\bf +} & GN1 & (1.2, 3.8, 4.5, 2.4) & 528 & 278 \\
		{\bf --} & {\bf +} & {\bf --} & GN1 & (1.0, 3.8, 8.6, 10.1) & 1091 & 99 \\
		{\bf --} & {\bf +} & {\bf +} & GN1 & (1.0, 3.8, 8.8, 10.1) & 1179 & 187 \\
		{\bf +} & {\bf --} & {\bf --} & GN1 & ($5.0{\times} 10^{-1}$, $3.1{\times} 10^{-1}$, $3.1{\times} 10^{-1}$, $1.6{\times} 10^{-1}$) & 367 & 744 \\
		{\bf +} & {\bf --} & {\bf +} & GN1 & ($7.5{\times} 10^{-1}$, $5.0{\times} 10^{-1}$, $2.1{\times} 10^{-1}$, $2.5{\times} 10^{-1}$) & 437 & 356 \\
		{\bf +} & {\bf +} & {\bf --} & GN1 & (0.4, 0.3, 0.9, 1.6) & 1167 & 374 \\
		{\bf +} & {\bf +} & {\bf +} & {GN1} & ($\bf 4.8{\times} 10^{-3}$, $\bf 4.3{\times} 10^{-3}$, $\bf 3.1{\times} 10^{-2}$, $\bf 7.6{\times} 10^{-3}$) & 1225 & 560 \\
		{\bf --} & {\bf --} & {\bf --} & GN2 & (1.0, 3.7, 7.8, 8.9) & 1370 & 21887 \\
		{\bf --} & {\bf --} & {\bf +} & GN2 & (1.0, 3.7, 7.5, 8.6) & 1668 & 21857 \\
		{\bf --} & {\bf +} & {\bf --} & GN2 & (1.0, 3.8, 7.6, 10.0) & 1116 & 178 \\
		{\bf --} & {\bf +} & {\bf +} & GN2 & (1.0, 3.8, 7.6, 10.0) & 1186 & 181 \\
		{\bf +} & {\bf --} & {\bf --} & GN2 & ($3.5{\times} 10^{-1}$, $2.1{\times} 10^{-1}$, $1.6{\times} 10^{-1}$, $2.4{\times} 10^{-1}$) & 359 & 693 \\
		{\bf +} & {\bf --} & {\bf +} & GN2 & ($3.6{\times} 10^{-1}$, $2.2{\times} 10^{-1}$, $1.7{\times} 10^{-1}$, $2.4{\times} 10^{-1}$) & 431 & 281 \\
		{\bf +} & {\bf +} & {\bf --} & GN2 & ($4.2{\times} 10^{-2}$, $5.4{\times} 10^{-2}$, $8.4{\times} 10^{-2}$, $2.7{\times} 10^{-2}$) & 1078 & 245 \\
		{\bf +} & {\bf +} & {\bf +} & GN2 & ($4.0{\times} 10^{-2}$, $5.0{\times} 10^{-2}$, $1.6{\times} 10^{-1}$, $3.1{\times} 10^{-2}$) & 1050 & 394 \\
		
		\bottomrule
	\end{tabular}
	\caption{Relativistic slingshot: Relative ${\mathbb L}_2$ error, run time and iterations of LBFGS optimizer for illustrating the impact of individual components of the proposed methods of training}
	\label{table4}
\end{table*}

\begin{figure}[t!]
	\centering
	\begin{subfigure}[t]{0.24\textwidth}
		\centering
		\includegraphics[width=\linewidth]{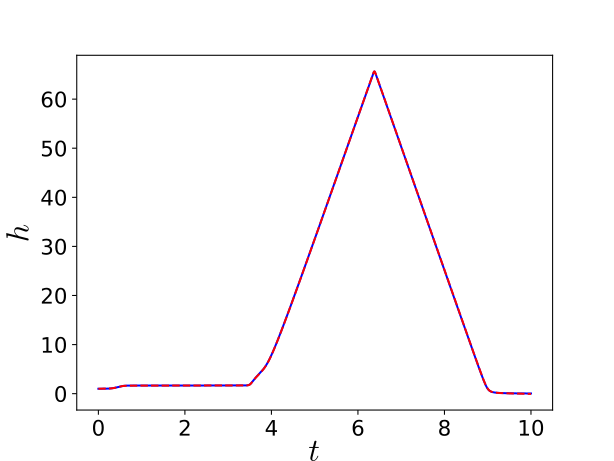}
		\caption{$h$}
	\end{subfigure}
	\begin{subfigure}[t]{0.24\textwidth}
		\centering
		\includegraphics[width=\linewidth]{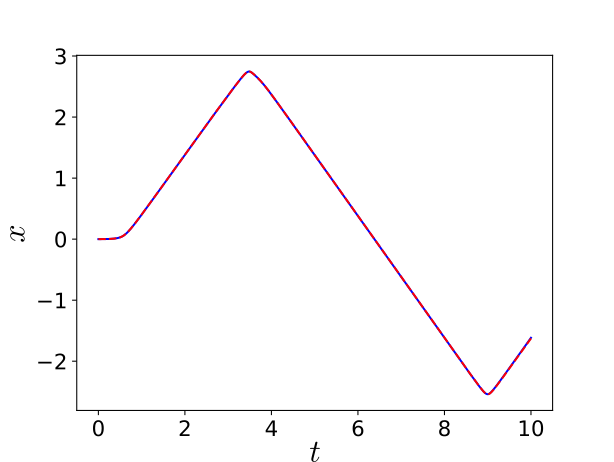}
		\caption{$x$}
	\end{subfigure}
	\begin{subfigure}[t]{0.24\textwidth}
		\centering
		\includegraphics[width=\linewidth]{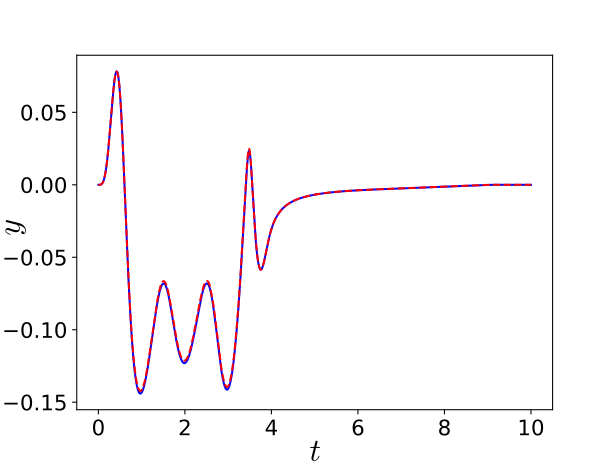}
		\caption{$y$}
	\end{subfigure}
	\begin{subfigure}[t]{0.24\textwidth}
		\centering
		\includegraphics[width=\linewidth]{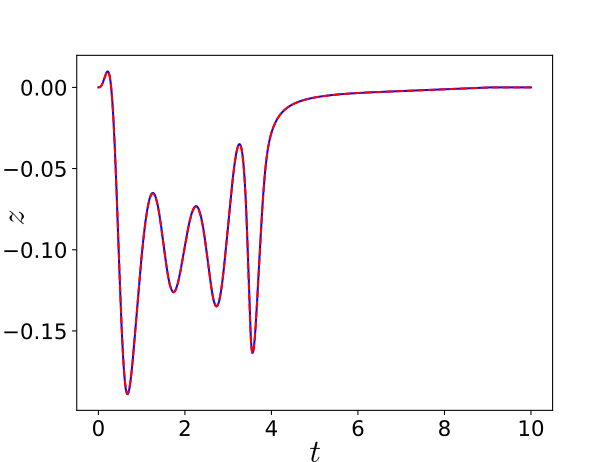}
		\caption{$z$}
	\end{subfigure}

	\begin{subfigure}[t]{0.24\textwidth}
		\centering
		\includegraphics[width=\linewidth]{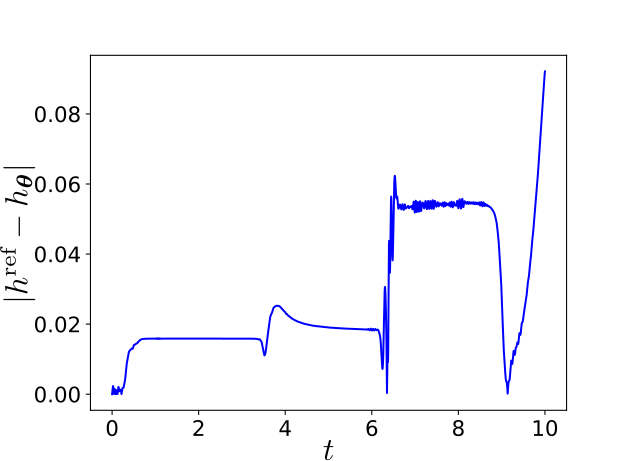}
		\caption{$h$}
	\end{subfigure}
	\begin{subfigure}[t]{0.24\textwidth}
		\centering
		\includegraphics[width=\linewidth]{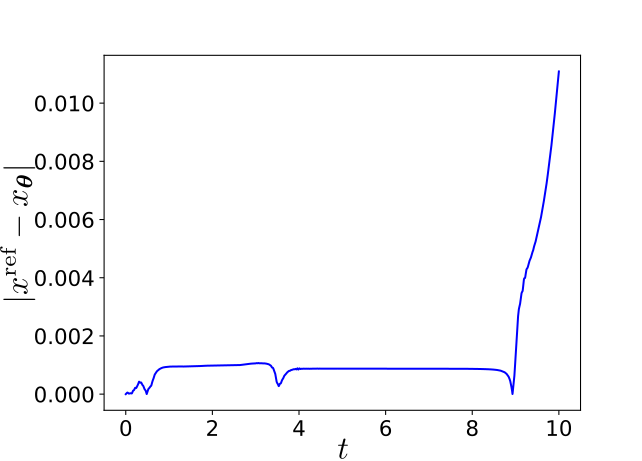}
		\caption{$x$}
	\end{subfigure}
	\begin{subfigure}[t]{0.24\textwidth}
		\centering
		\includegraphics[width=\linewidth]{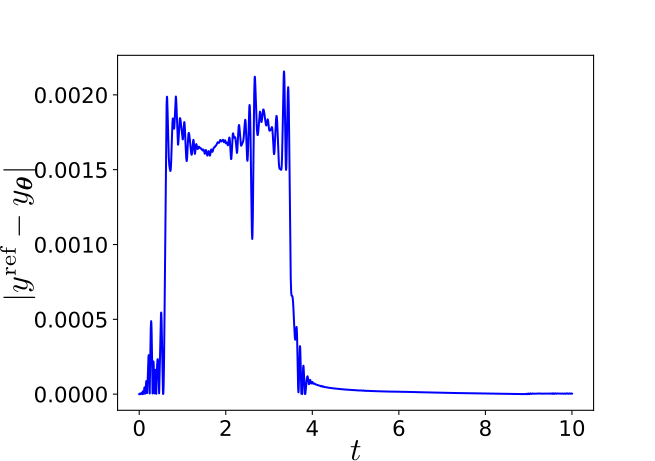}
		\caption{$y$}
	\end{subfigure}
	\begin{subfigure}[t]{0.24\textwidth}
		\centering
		\includegraphics[width=\linewidth]{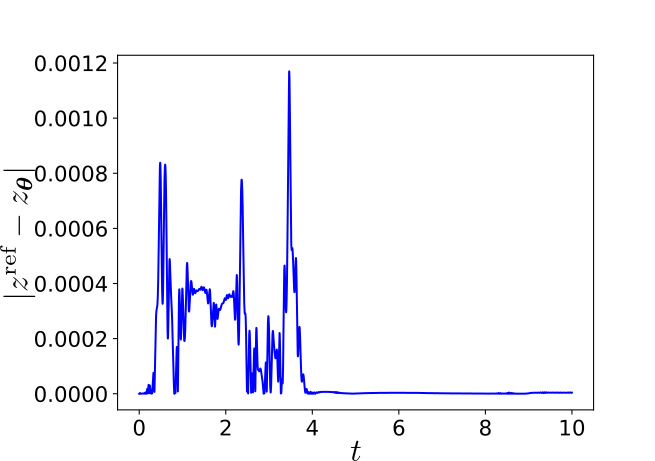}
		\caption{$z$}
	\end{subfigure}
	\caption{Relativistic slingshot. (a), (b), (c) and (d) are comparisons of the predicted (red dash lines) and reference solutions (blue solid lines) corresponding to $h$, $x$, $y$ and $z$, respectively. (e), (f), (g) and (h) are absolute errors $|h^{\text{ref}} - h_{\bm \theta}|$, $|x^{\text{ref}} - x_{\bm \theta}|$, $|y^{\text{ref}} - y_{\bm \theta}|$, and $|z^{\text{ref}} - z_{\bm \theta}|$, respectively. Number of neurons is $N=512$, the number of collocation points $N_t = 5048$, maximum of time is $T=10$.}
	\label{fig18}
\end{figure}

Results of four-stage training which consisted of 5000 epochs of optimization by Adam optimizer with learning rate $10^{-3}$ and 5 epochs of optimization by LBFGS optimizer (maximal number of iterations per optimization step is 5000) under frozen parameters of the hidden layer of neural network, then 5000 epochs of optimization by Adam optimizer with learning rate $10^{-6}$ and 5 epochs of optimization by LBFGS optimizer under unfrozen parameters of the hidden layer of neural network are shown on Figure~\ref{fig18}. The parameters are the number of neurons $N=512$, the number of collocation points $N_t=5048$, $T=10$, and $\Delta \zeta = 0.7$. The reference solutions are obtained by using the \verb*|odeint| solver of \verb*|scipy.integrate| library. The relative ${\mathbb L}_2$ errors are $\epsilon[{h}_{\bm \theta}, {h^{\text{ref}}}] = 1.2\times 10^{-3}$, $\epsilon[{x}_{\bm \theta}, {x^{\text{ref}}}] = 1.3\times 10^{-3}$, $\epsilon[{y}_{\bm \theta}, {y^{\text{ref}}}] = 1.6\times 10^{-2}$, and $\epsilon[{z}_{\bm \theta}, {z^{\text{ref}}}] = 3.2\times 10^{-3}$. The training was carried out on CPU and took 421 iterations of LBFGS optimizer and $1370$ seconds of calculation time.

\subsubsection{Example 10. Relativistic slingshot in interval $t\in[9,15]$}
\begin{figure}[t!]
	\centering
	\begin{subfigure}[t]{0.24\textwidth}
		\centering
		\includegraphics[width=\linewidth]{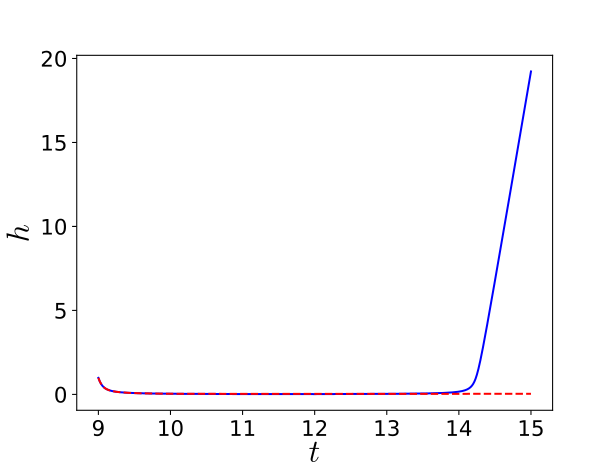}
		\caption{$h$}
	\end{subfigure}
	\begin{subfigure}[t]{0.24\textwidth}
		\centering
		\includegraphics[width=\linewidth]{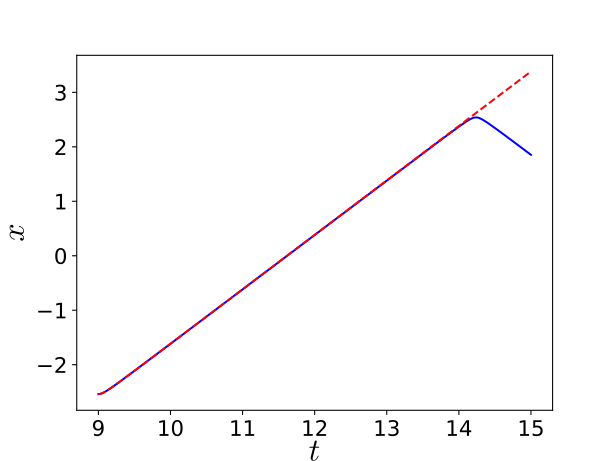}
		\caption{$x$}
	\end{subfigure}
	\begin{subfigure}[t]{0.24\textwidth}
		\centering
		\includegraphics[width=\linewidth]{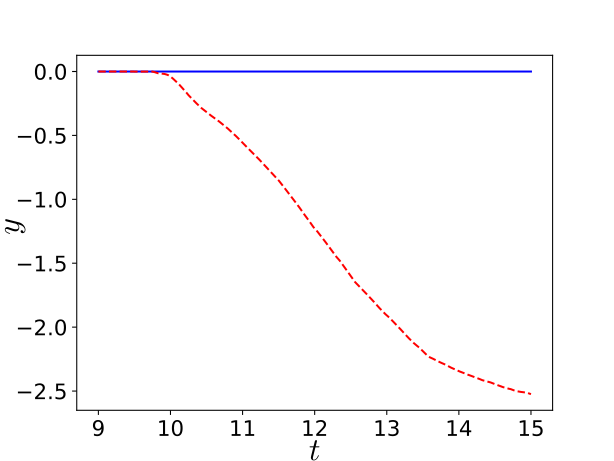}
		\caption{$y$}
	\end{subfigure}
	\begin{subfigure}[t]{0.24\textwidth}
		\centering
		\includegraphics[width=\linewidth]{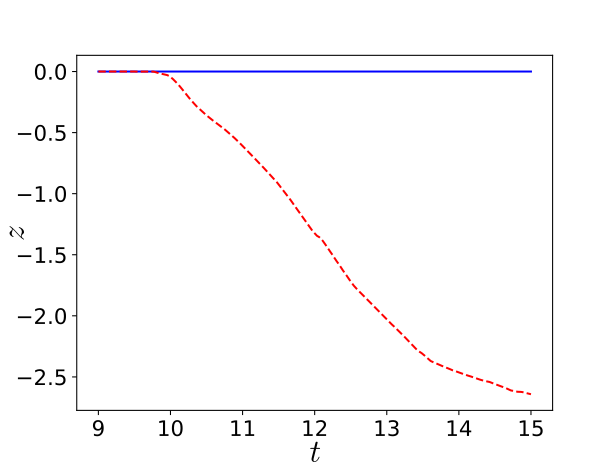}
		\caption{$z$}
	\end{subfigure}
	\caption{Relativistic slingshot. (a), (b), (c) and (d) are comparisons of the predicted (red dash lines) and reference solutions (blue solid lines) corresponding to $h$, $x$, $y$ and $z$, respectively. Number of neurons is $N=512$, the number of collocation points $N_t = 5048$, for $t\in [9,15]$.}
	\label{fig19}
\end{figure}

For the same hyperparameters that are for Figure~\ref{fig18}, but for the time interval $t\in [9,15]$ the results are shown on the Figure~\ref{fig19}. The relative ${\mathbb L}_2$ errors are $\epsilon[{h}_{\bm \theta}, {h^{\text{ref}}}] = 1.0$, $\epsilon[{x}_{\bm \theta}, {x^{\text{ref}}}] = 2.0\times 10^{-1}$, $\epsilon[{y}_{\bm \theta}, {y^{\text{ref}}}] = 335534$, and $\epsilon[{z}_{\bm \theta}, {z^{\text{ref}}}] = 387025$. The training was carried out on CPU and took 85 iterations of LBFGS optimizer and time is $617$ seconds.

\subsubsection{Example 11. The exponential growth equation}
Consider the exponential growth, which is governed by the equation for the $t\in [0, T]$
\begin{eqnarray}
	&& \frac{d u}{d t} = u, \quad t \in [0,T],\label{eq39}\\
	&& u(0) = 1. \label{eq40}
\end{eqnarray}
We used the following maximums of time $T=5$ and $T=10$. The exact analytical solution of this problem is $u^{\text{ref}} = \exp\left( t \right)$.
The latent variable $u$ is represented by neural networks $u_{\bm \theta}$ described above.

Training and evaluating parameters are the same as in the previous example. In out calculations the weights of loss (\ref{eq4}) were taken as $\lambda_{r} = 10$ and $\lambda_{ic} = \lambda_{r} N_t / T $.

\begin{figure}[t!]
	\centering
	\begin{subfigure}[t]{0.4\textwidth}
		\centering
		\includegraphics[width=\linewidth]{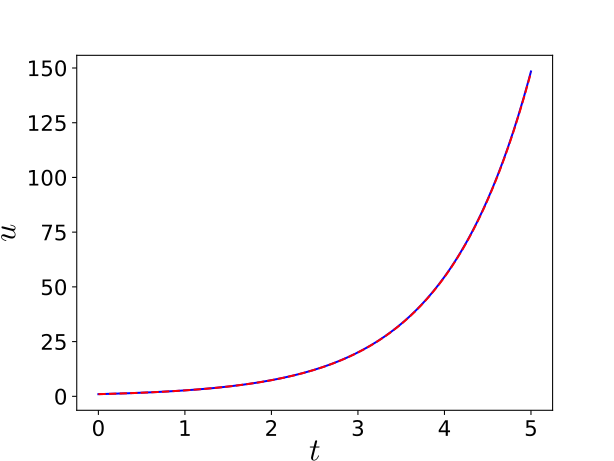}
		\caption{}
	\end{subfigure}
	\begin{subfigure}[t]{0.4\textwidth}
		\centering
		\includegraphics[width=\linewidth]{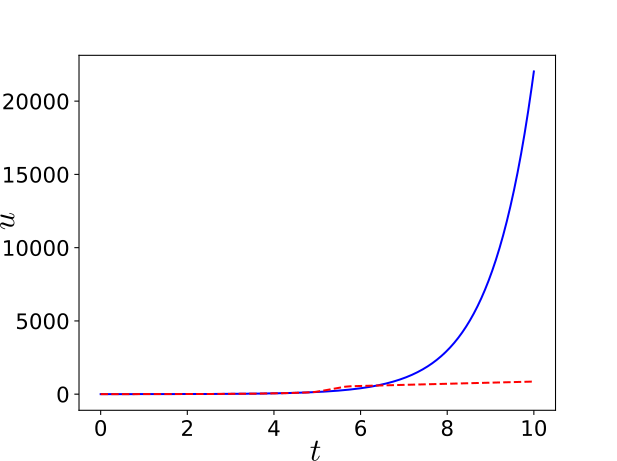}
		\caption{}
	\end{subfigure}
	\caption{The exponential growth. (a) and (b) are comparisons of the predicted (red dash lines) and reference solutions (blue solid lines) corresponding to $u$ for $T=5$ and $T=10$, respectively. Number of neurons is $N=512$, the number of collocation points $N_t = 5048$.}
	\label{fig20}
\end{figure}
Results of six-stage training for a neural network with $N=512$, number of collocation points $N_t=5048$, $T=5$ and $T=10$ are shown on Figure~\ref{fig20} (a) and (b), respectively. Training consisted of three stages for charged sphere calculation under the detach approach, and the same three stages under detach, $\delta$-causal, and gradient normalization GN1 approaches. The relative ${\mathbb L}_2$ errors is $\epsilon[u_{\bm \theta}, {u^{\text{ref}}}] = 2.9 {\times} 10^{-3}$ for (a) and ${\mathbb L}_2$ errors is $\epsilon[u_{\bm \theta}, {u^{\text{ref}}}] = 0.93$ for (b).

\subsection{Discussion}
As can be seen from the results of the previous subsection, the combination of the proposed learning methods significantly increases the accuracy of the solutions provided by neural networks. The last two examples (10 and 11) demonstrate that the proposed training methods are insufficient for that domain of the arguments of the solutions where the values of the solutions differ by orders of magnitude. To eliminate this disadvantage, we use very simple weighting schemes based on the decisions predicted by the neural network.

{\bf Suggestion 5: Predicted solution based weighting (PSBW).}

The residual loss term (\ref{eq6}) for every collocation points and initial condition term (\ref{eq5}) at $k$th iteration of training we are weighting by following scheme 
\begin{eqnarray}
	&& \mathcal{L}_{ic}^l \left( {\bm \theta} \right)  = \frac{1}{\left(v^{(ic)}_l + \epsilon_0\right)^2} \left| {u}_{{\bm \theta};l}\left(0\right) - {g}_l  \right|^{2},\label{eq41}\\
	&& \mathcal{L}_{r}^l\left(\bm \theta \right) = \frac{1}{N_r} \sum_{i=1}^{N_{r}} \frac{1}{\left(v^{(r)}_{l;i} + \epsilon_0\right)^2} \left| \mathcal{R}^l\left[{\vec u}_{{\bm \theta}} \right] \left({x}_{i}^{(r)} \right)\right|^{2}\label{eq42}.
\end{eqnarray}
Here
\begin{eqnarray}
	&& (v^{(ic)}_l)_k = \gamma(v^{(ic)}_l)_{k-1} + (1-\gamma ) \left|{u}_{{\bm \theta};l}\left(0\right)\right|,\label{eq43}\\
	&& \left(v^{(r)}_{l;i}\right)_{k} =  \gamma\left(v^{(r)}_{l;i}\right)_{k-1} + (1-\gamma ) \left|{u}_{{\bm \theta};l}\left({x}_{i}^{(r)} \right)\right| \label{eq44},
\end{eqnarray}
$\gamma$ is the decay parameter ($\gamma=0.9$ in our calculations), initial values are $(v^{(ic)}_l)_{0}={g}_l$ and  $\left(v^{(r)}_{l;i}\right)_{0}=1$, $\epsilon_0$ is parameter that excludes division by zero ($\epsilon=10^{-8}$ in our calculations), index $k$ means the the training iteration number. Due to this weighting, the elements of the loss function are aligned in the scale of values of the solution.

{\bf Suggestion 6: Relative Residuals (RR)}

Based on Suggestion 5, we have implemented a more balanced weighting scheme. The residual loss term (\ref{eq6}) for every collocation point at $k$th iteration of training we are weighting by following the scheme 
\begin{eqnarray}
	&& \mathcal{L}^l_{r}\left(\bm \theta \right) = \frac{1}{N_r} \sum_{i=1}^{N_{r}} \frac{1}{\left({\nu}^{(r)}_{l;i} + \epsilon_0\right)^2} \left| \mathcal{R}^l\left[{\vec u}_{{\bm \theta}} \right] \left({x}_{i}^{(r)} \right)\right|^{2}\label{eq45}.
\end{eqnarray}
Here
\begin{eqnarray}
	&& \left({\nu}^{(r)}_{l;i}\right)_{k} =  \gamma\left({\nu}^{(r)}_{l;i}\right)_{k-1} + (1-\gamma ) \left|\rasymbol{\mathcal{N}}_l\left[ {\vec u}_{{\bm \theta}}\left({x}_{i}^{(r)} \right), {x}_{i}^{(r)}  \right] \right| \label{eq46},
\end{eqnarray}
$\gamma$ is the decay parameter ($\gamma=0.99$ in our calculations), initial values are $\left(v^{(r)}_{l;i}\right)_{0}=1$, $\epsilon_0$ is parameter that excludes division by zero ($\epsilon_0=10^{-8}$ in our calculations). Due to this weighting, the elements of the loss function are aligned in the scale of values of the solution.

\subsection{Numerical experiments}

\subsubsection{Example 12. The exponential growth equation}\label{Ex12}

In this part of our experiments we used three-stage training which consisted of 5000 epochs of optimization by Adam optimizer with a learning rate $10^{-3}$ and 15 epochs of optimization by LBFGS optimizer (maximal number of iterations per optimization step is 5000) under frozen parameters of the hidden layer of neural network, then 15 epochs of optimization by LBFGS optimizer under unfrozen parameters of the hidden layer of neural network. The initial values of lambdas were taken $\lambda^{k}_{r} = 10$ and $\lambda_{ic}^{k} = {\beta^k} \lambda^{k}_{r}$, ${\beta^k}=\dfrac{N_t}{T}$. Other parameters are $N=512$, number of collocation points $N_t=5048$, $T=10$ and $\Delta \zeta = 0.7$. A detaching approach was used.  Results are presented in the Table~{\ref{table5}}. The gradient normalization approach wasn't considered because this method did not lead to increasing accuracy. The best result is achieved when using $\delta$-causal and PSBW weighting approaches.

\begin{table*}[h!]
	\centering
	\begin{tabular}{cc|ccc}
		\toprule
		\multicolumn{2}{c}{\bf{Train Settings}} & \multicolumn{3}{c}{\bf{Performance}}\\
		\midrule
		\bf{$\delta$-causal} & \bf{PSBW}  &  \bf{Relative ${\mathbb L}_2$ errors ($u$)} & \bf{Run time (s)} & \bf{Iterations of LBFGS} \\
		\midrule
		{\bf --} & {\bf --} & {$0.88$} & {124} & {554} \\
		{\bf +} & {\bf --} & {$0.97$} & 192 & 6105 \\
		{\bf --} & + & {$0.88$} & 93 & 554 \\
		{\bf +} & {\bf +} & {${\bf 1.7{\times} 10^{-3}}$} & 941 & 53948 \\
		\bottomrule
	\end{tabular}
	\caption{The exponential growth equation: Relative ${\mathbb L}_2$ error, run time and iterations of LBFGS optimizer for illustrating the impact of individual components of the proposed methods of training}
	\label{table5}
\end{table*}

\begin{figure}[t!]
	\centering
	\begin{subfigure}[t]{0.3\textwidth}
		\centering
		\includegraphics[width=\linewidth]{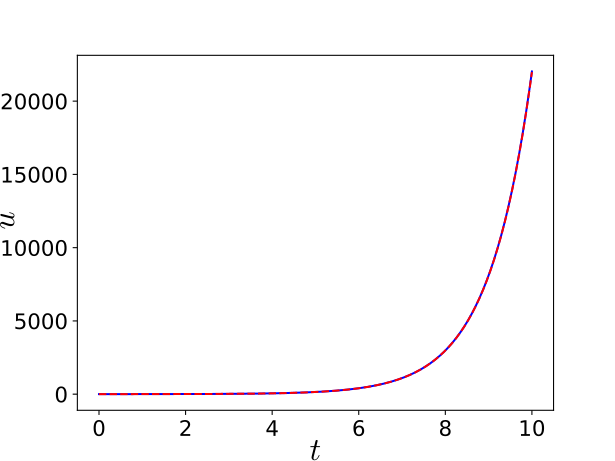}
		\caption{}
	\end{subfigure}
	\begin{subfigure}[t]{0.3\textwidth}
		\centering
		\includegraphics[width=\linewidth]{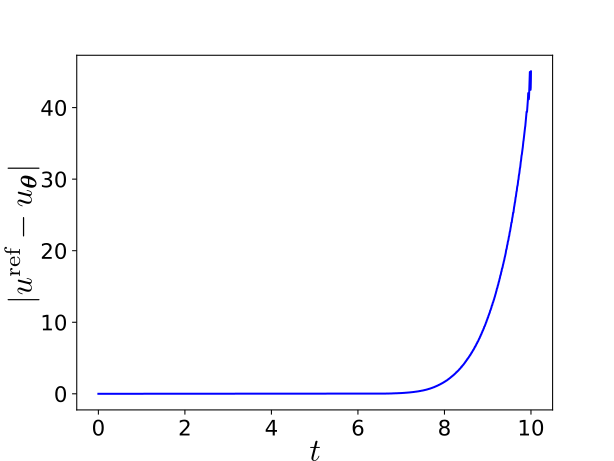}
		\caption{}
	\end{subfigure}
	\begin{subfigure}[t]{0.32\textwidth}
		\centering
		\includegraphics[width=\linewidth]{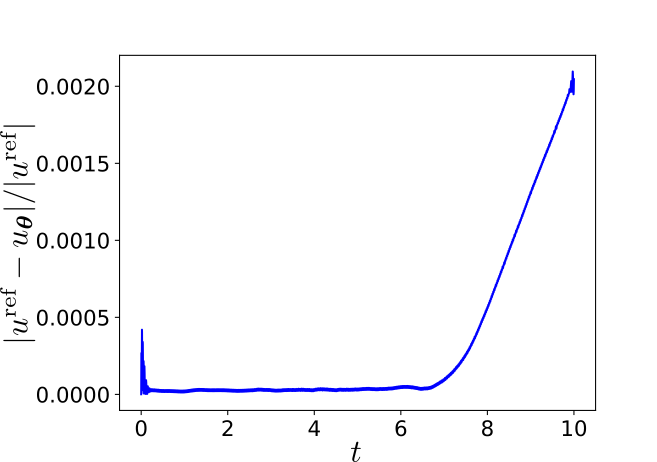}
		\caption{}
	\end{subfigure}
	\caption{The exponential growth. (a) is a comparison of the predicted (red dash lines) and reference solutions (blue solid lines) corresponding to $u$. (b) and (c) are absolute errors $|u^{\text{ref}} - u_{\bm \theta}|$ and relative error $|u^{\text{ref}} - u_{\bm \theta}|/|u^{\text{ref}}|$. Number of neurons is $N=512$, the number of collocation points $N_t = 5024$, $T=10$.  The relative ${\mathbb L}_2$ errors is $\epsilon[u_{\bm \theta}, {u^{\text{ref}}}] = 1.72{\times} 10^{-3}$.}
	\label{fig21}
\end{figure}

Results of training for a neural network with $N=512$, number of collocation points $N_t=5048$ and $t=10$ are shown on Figure~\ref{fig21}. The relative ${\mathbb L}_2$ errors is $\epsilon[u_{\bm \theta}, {u^{\text{ref}}}] = 1.72{\times} 10^{-3}$. In our calculation the initial value of causal parameter $\varepsilon$ is taken $\varepsilon = 10^{-8}$.

\subsubsection{Example 13. Relativistic slingshot in interval $t\in[9,20]$}
In this part of our experiments, we used four-stage training which is consisted of 5000 epochs of optimization by Adam optimizer with a learning rate $10^{-3}$ and 15 epochs of optimization by LBFGS optimizer (maximal number of iterations per optimization step is 5000) under frozen parameters of the hidden layer of the neural network, then 50000 epochs of optimization by Adam optimizer with initial learning rate $10^{-5}$ and step-learning-rate scheduler with decay--rate is $0.9$ for every 5000 training epochs every, and then 15 epochs of optimization by LBFGS optimizer under unfrozen parameters of the hidden layer of neural network. The parameters $(v^{(ic)}_l)_k$ and $\left(v^{(r)}_{l;i}\right)_{k}$ are recalculating by Eqns.~(\ref{eq43}) and (\ref{eq44}) every 1000 epochs. The values of lambdas were taken $\lambda^{k}_{r} = 1$ and $\lambda_{ic}^{k} = {\beta^k} \lambda^{k}_{r}$, ${\beta^k}=\dfrac{N_t}{T}$. Other parameters are $N=512$, number of collocation points $N_t=5048$, $T=11$ (initial $t=9$ and finish $t=20$), $\Delta \zeta = 0.7$ and the initial value of causality parameter is taken $\varepsilon = 10^{-8}$.

\begin{figure}[t!]
	\centering
	\begin{subfigure}[t]{0.24\textwidth}
		\centering
		\includegraphics[width=\linewidth]{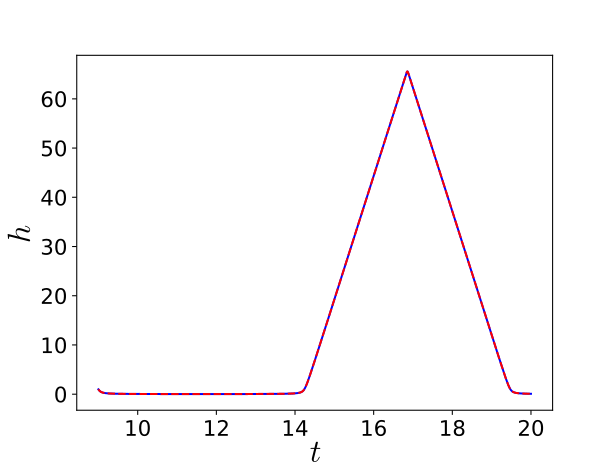}
		\caption{$h$}
	\end{subfigure}
	\begin{subfigure}[t]{0.24\textwidth}
		\centering
		\includegraphics[width=\linewidth]{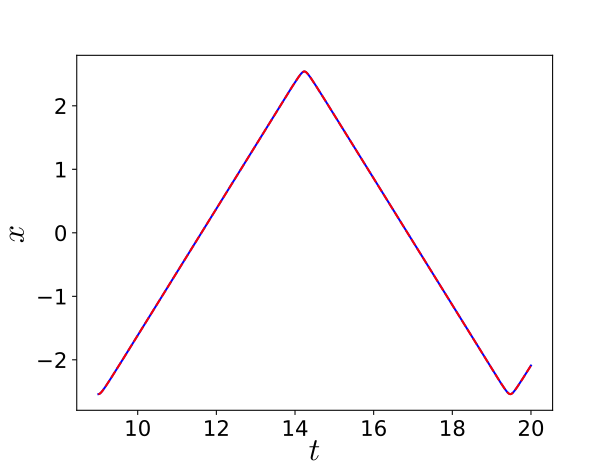}
		\caption{$x$}
	\end{subfigure}
	\begin{subfigure}[t]{0.24\textwidth}
		\centering
		\includegraphics[width=\linewidth]{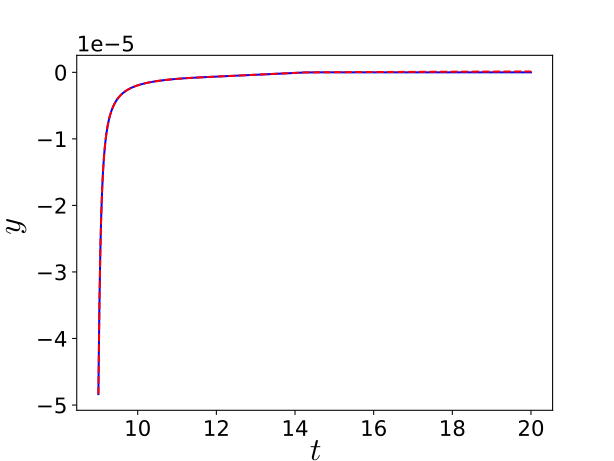}
		\caption{$y$}
	\end{subfigure}
	\begin{subfigure}[t]{0.24\textwidth}
		\centering
		\includegraphics[width=\linewidth]{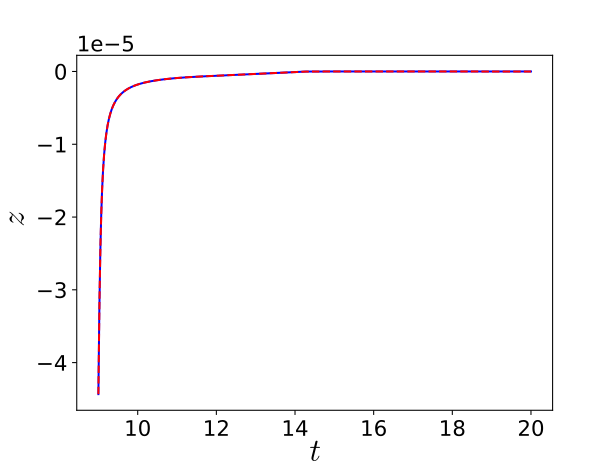}
		\caption{$z$}
	\end{subfigure}

	\begin{subfigure}[t]{0.24\textwidth}
		\centering
		\includegraphics[width=\linewidth]{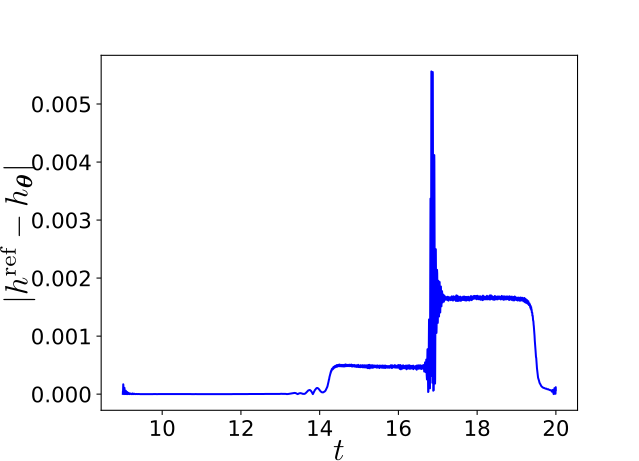}
		\caption{$h$}
	\end{subfigure}
	\begin{subfigure}[t]{0.24\textwidth}
		\centering
		\includegraphics[width=\linewidth]{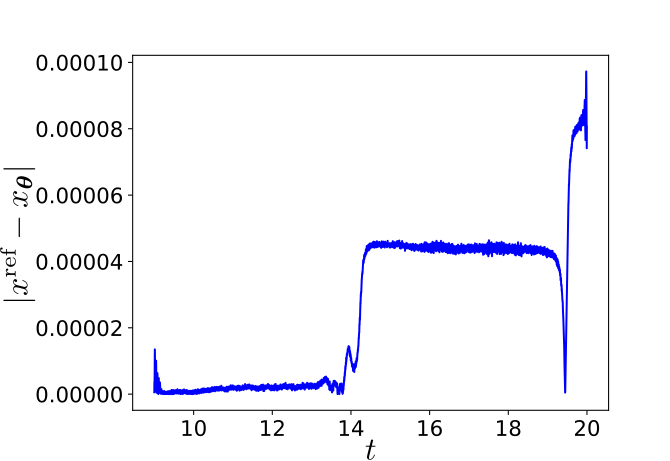}
		\caption{$x$}
	\end{subfigure}
	\begin{subfigure}[t]{0.24\textwidth}
		\centering
		\includegraphics[width=\linewidth]{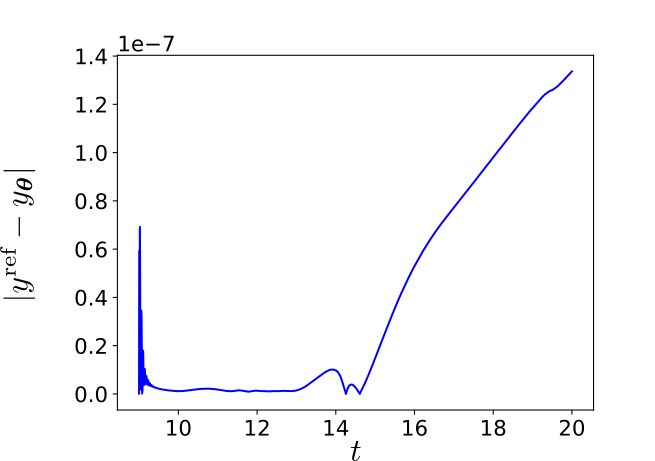}
		\caption{$y$}
	\end{subfigure}
	\begin{subfigure}[t]{0.24\textwidth}
		\centering
		\includegraphics[width=\linewidth]{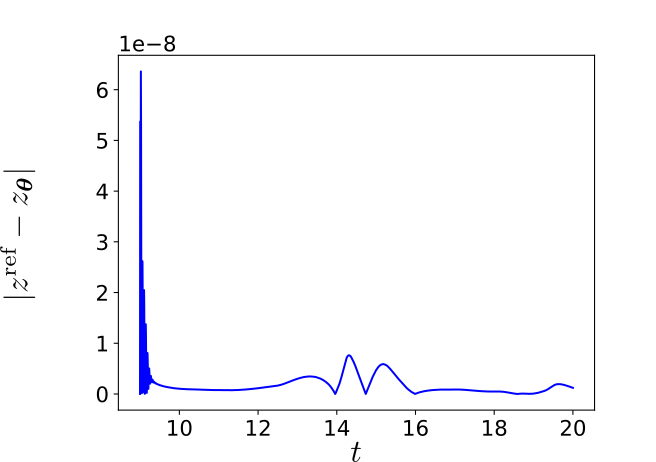}
		\caption{$z$}
	\end{subfigure}
	\caption{Relativistic slingshot. (a), (b), (c) and (d) are comparisons of the predicted (red dash lines) and reference solutions (blue solid lines) corresponding to $h$, $x$, $y$ and $z$, respectively. (e), (f), (g) and (h) are absolute errors $|h^{\text{ref}} - h_{\bm \theta}|$, $|x^{\text{ref}} - x_{\bm \theta}|$, $|y^{\text{ref}} - y_{\bm \theta}|$, and $|z^{\text{ref}} - z_{\bm \theta}|$, respectively. Number of neurons is $N=512$, the number of collocation points $N_t = 5048$. The relative ${\mathbb L}_2$ errors are $\epsilon[{h}_{\bm \theta}, {h^{\text{ref}}}] = 3.4\times 10^{-5}$, $\epsilon[{x}_{\bm \theta}, {x^{\text{ref}}}] = 2.2\times 10^{-5}$, $\epsilon[{y}_{\bm \theta}, {y^{\text{ref}}}] = 1.8\times 10^{-2}$, and $\epsilon[{z}_{\bm \theta}, {z^{\text{ref}}}] = 1.1\times 10^{-3}$.}
	\label{fig22}
\end{figure}

Results of four-stage training which consisted of 5000 epochs of optimization by Adam optimizer with learning rate $10^{-3}$ and 15 epochs of optimization by LBFGS optimizer (maximal number of iterations per optimization step is 5000) under frozen parameters of the hidden layer of neural network, then 50000 epochs of optimization by Adam optimizer with learning rate $10^{-5}$ and 15 epochs of optimization by LBFGS optimizer under unfrozen parameters of the hidden layer of neural network are shown on the Figure~\ref{fig22}. The reference solutions are obtained by using the \verb*|odeint| solver of \verb*|scipy.integrate| library. The relative ${\mathbb L}_2$ errors are $\epsilon[{h}_{\bm \theta}, {h^{\text{ref}}}] = 3.4\times 10^{-5}$, $\epsilon[{x}_{\bm \theta}, {x^{\text{ref}}}] = 2.2\times 10^{-5}$, $\epsilon[{y}_{\bm \theta}, {y^{\text{ref}}}] = 1.8\times 10^{-2}$, and $\epsilon[{z}_{\bm \theta}, {z^{\text{ref}}}] = 1.1\times 10^{-3}$. The training was carried out on GPU and took 770 iterations of LBFGS optimizer and time is $516$ seconds.

\subsubsection{Example 14. Relativistic slingshot in interval $t\in[0,20]$}

In this part of our experiments, we used two-stage training which consisted of 5000 epochs of optimization by Adam optimizer with a learning rate $10^{-3}$ under frozen parameters of the hidden layer of the neural network, then  100000 epochs of optimization by Adam optimizer with initial learning rate $10^{-5}$ and step-learning-rate scheduler with decay--rate is $0.9$ for every 5000 training epochs every under unfrozen parameters of the hidden layer of neural network. $\delta$-causal and detaching were used for the calculations. The parameters $(v^{(ic)}_l)_k$ and $\left(\nu^{(r)}_{l;i}\right)_{k}$ are recalculating by Eqns.~(\ref{eq43}) and (\ref{eq46}) at every of training epoch. The values of lambdas were taken $\lambda^{k}_{r} = 1$ and $\lambda_{ic}^{k} = {\beta^k} \lambda^{k}_{r}$, ${\beta^k}=\dfrac{N_t}{T}$. Other parameters are $N=1024$, number of collocation points $N_t=10000$, $T=20$ (initial $t=0$ and finish $t=20$), $\Delta \zeta = 0.7$ and the initial value of causality parameter is taken $\varepsilon = 10^{-24}$.

\begin{figure}[t!]
	\centering
	\begin{subfigure}[t]{0.24\textwidth}
		\centering
		\includegraphics[width=\linewidth]{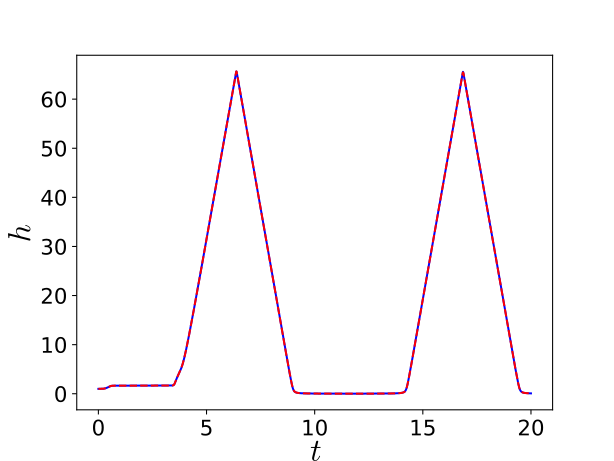}
		\caption{$h$}
	\end{subfigure}
	\begin{subfigure}[t]{0.24\textwidth}
		\centering
		\includegraphics[width=\linewidth]{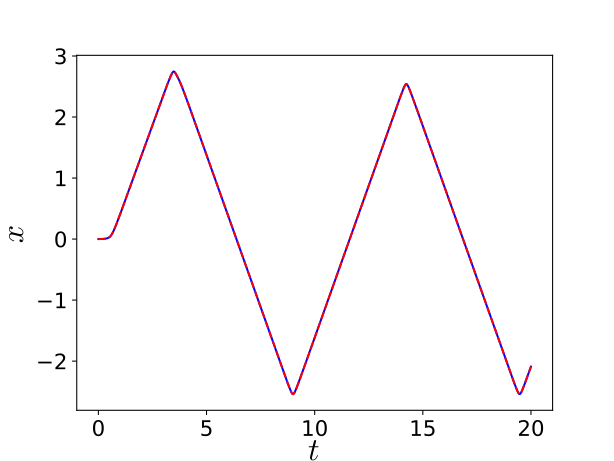}
		\caption{$x$}
	\end{subfigure}
	\begin{subfigure}[t]{0.24\textwidth}
		\centering
		\includegraphics[width=\linewidth]{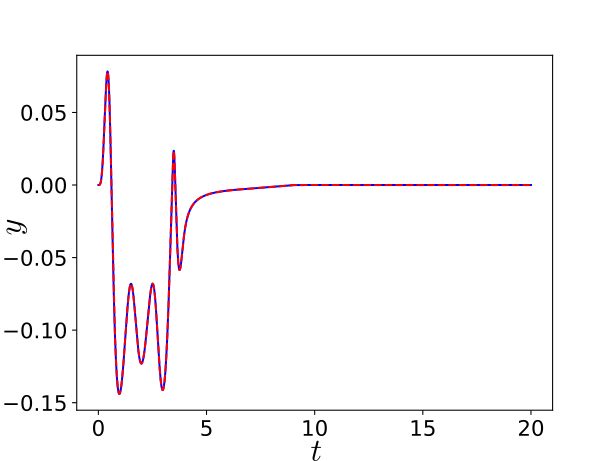}
		\caption{$y$}
	\end{subfigure}
	\begin{subfigure}[t]{0.24\textwidth}
		\centering
		\includegraphics[width=\linewidth]{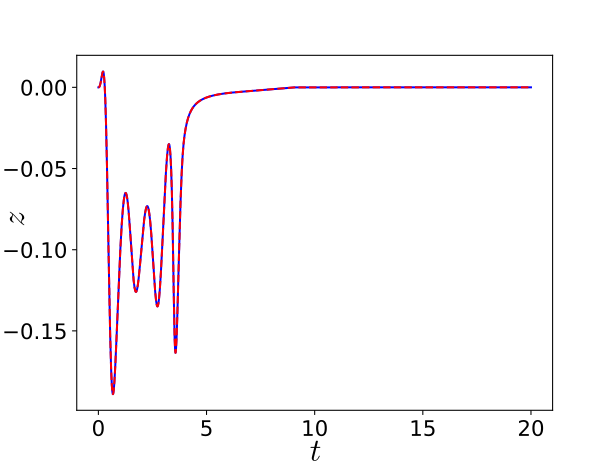}
		\caption{$z$}
	\end{subfigure}

	\begin{subfigure}[t]{0.24\textwidth}
		\centering
		\includegraphics[width=\linewidth]{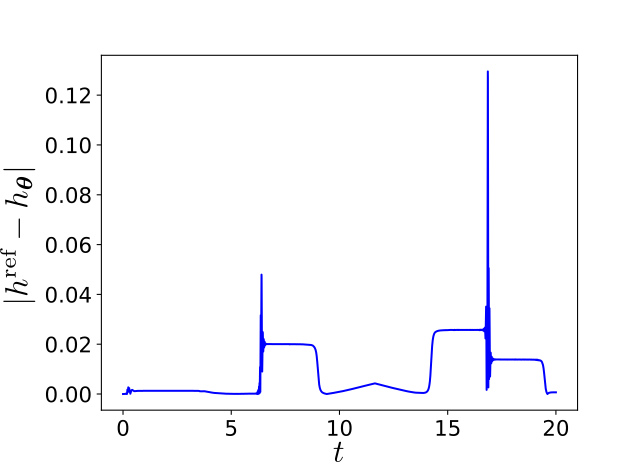}
		\caption{$h$}
	\end{subfigure}
	\begin{subfigure}[t]{0.24\textwidth}
		\centering
		\includegraphics[width=\linewidth]{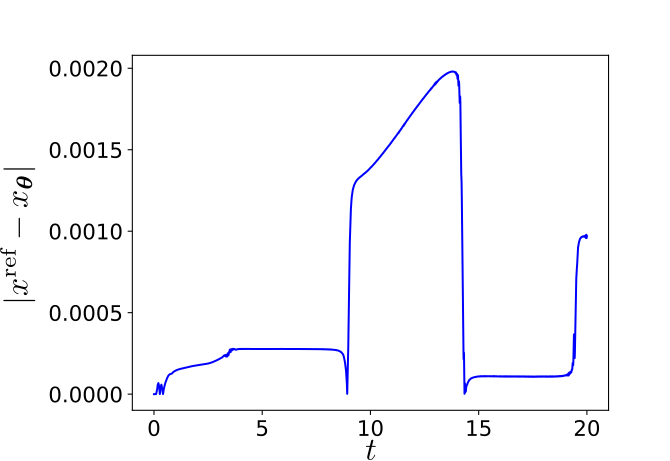}
		\caption{$x$}
	\end{subfigure}
	\begin{subfigure}[t]{0.24\textwidth}
		\centering
		\includegraphics[width=\linewidth]{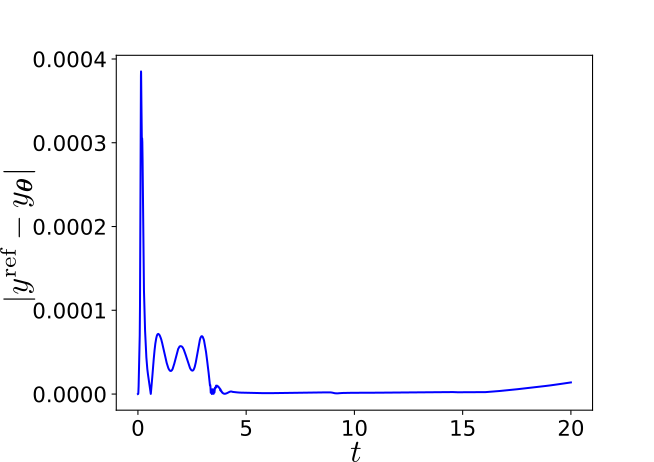}
		\caption{$y$}
	\end{subfigure}
	\begin{subfigure}[t]{0.24\textwidth}
		\centering
		\includegraphics[width=\linewidth]{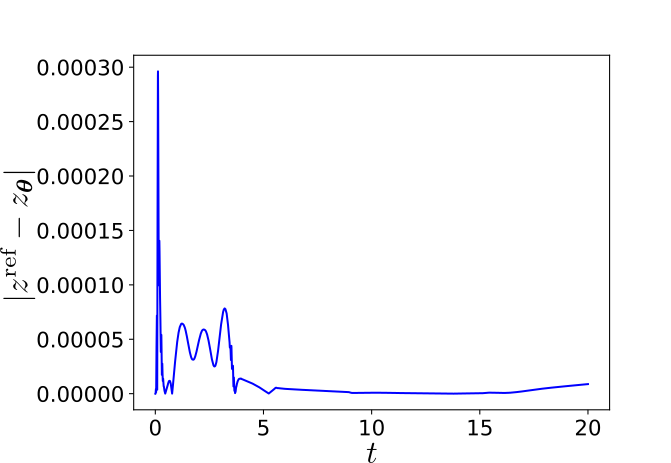}
		\caption{$z$}
	\end{subfigure}
	\caption{Relativistic slingshot. (a), (b), (c) and (d) are comparisons of the predicted (red dash lines) and reference solutions (blue solid lines) corresponding to $h$, $x$, $y$ and $z$, respectively. (e), (f), (g) and (h) are absolute errors $|h^{\text{ref}} - h_{\bm \theta}|$, $|x^{\text{ref}} - x_{\bm \theta}|$, $|y^{\text{ref}} - y_{\bm \theta}|$, and $|z^{\text{ref}} - z_{\bm \theta}|$, respectively. Number of neurons is $N=1024$, the number of collocation points $N_t = 10000$. The relative ${\mathbb L}_2$ errors are $\epsilon[{h}_{\bm \theta}, {h^{\text{ref}}}] = 4.9\times 10^{-4}$, $\epsilon[{x}_{\bm \theta}, {x^{\text{ref}}}] = 5.7\times 10^{-4}$, $\epsilon[{y}_{\bm \theta}, {y^{\text{ref}}}] = 8.1\times 10^{-4}$, and $\epsilon[{z}_{\bm \theta}, {z^{\text{ref}}}] = 5.2\times 10^{-4}$.}
	\label{fig23}
\end{figure}

Results of two-stage training which is consisted of 5000 epochs of optimization by Adam optimizer with learning rate $10^{-3}$ under frozen parameters of the hidden layer of neural network and 100000 epochs of optimization by Adam optimizer with a learning rate $10^{-5}$ under unfrozen parameters of the hidden layer of the neural network are shown on Figure~\ref{fig23}. The references solutions are obtained by using the \verb*|odeint| solver of \verb*|scipy.integrate| library.  The training was carried out on GPU and took time $1559$ seconds. The relative ${\mathbb L}_2$ errors for calculating by RR and PSBW given by Eqns.~(\ref{eq43}) and (\ref{eq46}) are presented in the table~\ref{table6}. It can be seen, the best result is achieved under RR approach of weighting of training.

\begin{table*}[h!]
	\centering
	\begin{tabular}{c|cc}
		\toprule
		{\bf{Weighting Method}} &  \bf{Relative ${\mathbb L}_2$ error (h, x, y, z)} & \bf{Run time (s)} \\
		\midrule
		{\bf RR} & (${\bf 4.9\times 10^{-4}}$, ${\bf 5.7\times 10^{-4}}$, ${\bf 8.1\times 10^{-4}}$, ${\bf 5.2\times 10^{-4}}$) & 1559 \\
		{\bf PSBW} &  ($6.1\times 10^{-3}$, $8.4\times 10^{-3}$, $5.8\times 10^{-2}$, $6.7\times 10^{-2}$) & 1743 \\		
		\bottomrule
	\end{tabular}
	\caption{Relativistic slingshot in interval $t\in[0,20]$: Relative ${\mathbb L}_2$ error and run time for RR and PSBW weighting methods}
	\label{table6}
\end{table*}

\subsection{Discussion}
The neural network training methods described above require calculating the gradient of the loss function. Despite the improvement in the learning results of neural networks for the presented problems, it would be better to have in the arsenal of training methods one that does not require calculating such a gradient and is relatively simple and intuitive. We will try to create and develop such a method in the next section of the article.

\newpage
\section{Gradient-free methods for fitting neural network parameters}
We begin this section with a kind of data-driven learning method and then extend some of its approaches to physics-informed training.

\subsection{Physics-informed data-driven initialization}
Let's assume that for each point with coordinate $x_k$ from uniformly distributed $N$ points in the interval $[0,X]$, the solution ${\vec u}_k$ of Equation (\ref{eq1}) is known. We use these coordinates and functions to obtain a trained neural network.

{\bf Suggestion 7: Physics-informed data-driven initialization.} 
Firstly, we initialize the parameters of the hidden layer of a neural network as in Algorithm~\ref{alg2}:
\begin{enumerate}
	\item Weights of the 1st layer are taken $2 \Delta \zeta / \Delta x$ (${W}^{(1)}_k = 2 \Delta \zeta / \Delta x$) and biases $b^{(1)}_k=-2 k \Delta \zeta$ for the $k$th neuron. Here, unlike algorithm~\ref{alg2}, we take the first bias of the hidden layer ($b^{(1)}_0$) corresponding to $x_0$.
	\item The parameters of the output layer are initialized with ${W}^{(2)}_{k} =  - 2 \Delta x \mathcal{N}_l [{\vec u}_k, x_k] / \Delta \zeta$.
	\item The bias of the output layer is taken equal to initial value ${u}_0$ ($b^{(2)}_0=u_0$).
\end{enumerate}
Here, to initialize ${W}^{(2)}_{k}$, the values of ${\vec u}_k$ at the points $x_k$ are taken, and not ${\vec u}_0$ for all the points, as in Algorithm~\ref{alg2}. For the derivation of weights ${W}^{(2)}_{k}$, as the derivation of Eqn.~(\ref{eq16}), it is assumed that the main contribution to the sum on the right-hand side of the equation~(\ref{eq14}) is made by the term containing $\sigma'(0)$.  However, this is not true for all values of $\Delta \zeta$. We attempt to derive a more accurate equation for these weights using the following equation
\begin{equation}\label{eq47}
	\sum\limits_{k = 0}^{N-1} \frac{2\Delta \zeta}{\Delta x}{W}^{(2)}_k\sigma'\left(\frac{2\Delta \zeta}{\Delta x} \left[x_m - x_k\right]\right) = - \mathcal{N}_l [{\vec u}_m, x_m].
\end{equation}
Note that the function $\sigma'(x)$ decreases relatively quickly as the absolute value of the argument increases. Based on this behaviour of $\sigma'(x)$ and the assumption that the weights $W^{(2)}_k$, for $k \in [m-L,m+L]$, do not change significantly from $W^{(2)}_{m}$, we can simplify equation (\ref{eq47}) to the following
\begin{equation}\label{eq48}
	\frac{2\Delta \zeta}{\Delta x}{W}^{(2)}_{m}\sum\limits_{k = m-L}^{m+L} \sigma'\left({2\Delta \zeta} \left[m - k\right]\right) = - \mathcal{N}_l [{\vec u}_m, x_m].
\end{equation}
Thus, we obtain the following equation for these weights
\begin{align}\label{eq49_1}
	& {W}^{(2)}_{m} = - \frac{\Delta x}{2\Delta \zeta} \frac{\mathcal{N}_l[{\vec u}_{m}, x_{m}]}{\kappa_m},\notag\\ 
	& \kappa_m = \sum\limits_{k = m-L}^{m+L} \sigma'\left({2\Delta \zeta} \left[m - k)\right]\right).
\end{align}
In our calculations, we limited $L$ to 10. It should be noted that the contributions of the terms $\sigma'(0)$, $\sigma'({2\Delta\zeta})$, and $\sigma'(-{2\Delta\zeta})$ to the value $\kappa_m$ are defining.

The bias of the output layer is taken equal to the initial value $u_{l;0}$ minus the output all layer at $x_0$ without this bias:
\begin{equation}\label{eq49_2}
	b^{(2)}_0 = u_{l;0} - \sum\limits_{k = 0}^{N-1} {W}^{(2)}_k\sigma\left({W}^{(1)}_k x_0 + {b}^{(1)}_k\right).
\end{equation}

We require that the higher-order derivatives in the response of the neural network at the point $x_m$ are determined mainly by the $m$th neuron. Fulfilling this requirement leads to better control of this neural network coefficient calculation scheme: neighboring neurons to the $m$ neuron lead to minimal distortion into its signal when responding at the point $x_m$. Since the response is determined mainly by this neuron and the nearest neurons, the following conditions must be satisfied
\begin{align}
	& \sigma''(-{2\Delta\zeta}) + \sigma''(0) + \sigma''({2\Delta\zeta}) = \sigma''(0),\label{eq50a}\\
	& \sigma'''(-{2\Delta\zeta}) + \sigma'''(0) + \sigma'''({2\Delta\zeta}) = \sigma'''(0). \label{eq50b}
\end{align}
Due to the oddness of the function $\sigma''(x)$, the equation (\ref{eq50a}) is executed for any $\Delta\zeta$. Due to the evenness of the function $\sigma''(x)$, the equality $\sigma'''(-{2\Delta\zeta}) = \sigma'''({2\Delta\zeta})$ is satisfied. Thus, it is necessary to solve the equation $\sigma'''({2\Delta\zeta}) = 0$, which has a following positive solution
\begin{align}
	{\Delta\zeta}= \ln\left(2 + \sqrt{3}\right) / 2. \label{eq51}
\end{align}
With such a $\Delta\zeta$ value in Eqn.~(\ref{eq49_1}) is $\kappa_k \approx 0.75934$ for $k\geq L$. Note that $\ln\left(2 +\sqrt{3}\right) /2 \approx 0.66$ is very close to the optimal value for a number of problems ${\Delta\zeta}=0.7$, which we experimentally determined in the previous section.

As a result, we have the following physics-informed data-driven (PIDD) initialization of neural networks (Algorithm~\ref{alg3}).
\begin{algorithm}[h!]
	\caption{Physics-informed data-driven initialization of neural networks}\label{alg3}
	\KwData{$\left\{{\vec u}_k\right\}_{k=0}^{N-1}$ for $\left\{x_k\right\}_{k=0}^{N-1}$ (uniform grid with a step $\Delta x \gets X / N$)}
	\KwResult{Initialized neural network ${u}_{{\bm \theta};l}(x)$ of PINN ${\vec u}_{\bm \theta}$, which consists of $N$ neurons on hidden layer}
	${\Delta\zeta} \gets \ln\left(2 + \sqrt{3}\right) / 2$\;
	\For{$k=1,\dots, N-1$}{
		${W}^{(1)}_k \gets 2 \Delta \zeta / \Delta x$\;
		${b}^{(1)}_k \gets -2 k \Delta \zeta$\;
		$\kappa_k \gets \sum\limits_{m = k-L}^{k+L} \sigma'\left({2\Delta \zeta} \left[k - m)\right]\right)$\;
		${W}^{(2)}_k \gets - \dfrac{\Delta x}{2\Delta \zeta} \dfrac{\mathcal{N}[{\vec u}_{k}, x_{k}]}{\kappa_k}$\;}
	${b}^{(2)}_0 \gets {u}_{l;0} - \sum\limits_{k = 0}^{N-1} {W}^{(2)}_k\sigma\left({W}^{(1)}_k x_0 + {b}^{(1)}_k\right)$.
\end{algorithm}

\subsection{Numerical experiments}
\subsubsection{Example 15. Harmonic Oscillator}

In our experiments, we used PIDD initialization, which is described above (Algorithm~\ref{alg3}). The reference solutions are obtained by using the \verb*|odeint| solver of \verb*|scipy.integrate| library. The data $\left\{{\vec u}_k\right\}_{k=0}^N$ for uniform distribution points $\left\{x_k\right\}_{k=0}^N$ were generated with \verb*|odeint| for $N=20000$  in ${\bf 0.004}$ seconds.

Results of PIDD initialization for neural network with $N=20000$ are following: the relative ${\mathbb L}_2$ errors are $\epsilon[{u}_{\bm \theta;1}, {u_1^{\text{ref}}}] = 5.67{\times} 10^{-5}$ and $\epsilon[{u}_{\bm \theta;2}, {u_2^{\text{ref}}}] = 6.82{\times} 10^{-4}$. Such initialization took ${\bf 0.076}$ seconds of computing time of CPU. It can be seen, it took much less time for such initialization than with standard training (compare with Example 4), while the size of the hidden layer increased significantly. 

The dependence of relative ${\mathbb L}_2$ error $\epsilon[{u}_{\bm \theta;1}, {u_1^{\text{ref}}}]$ on the number of neurons $N$ of neural network are shown on the Figure~\ref{fig25}(a) by blue solid line. The red dashed line corresponds to dependence $\alpha_0 / N$, where $\alpha_0$ related to value $\epsilon[{u}_{\bm \theta;1}, {u_1^{\text{ref}}}]$ for $N=1000$. It can be seen, the accuracy of the neural network increases with the increasing the number of neurons faster than the accuracy of the method of integrating a differential equation based on rectangle methods.

\subsubsection{Example 16. Lorentz system}
As a further example, let us consider the chaotic Lorenz system. This system of equations arises in studies of convection and instability in planetary atmospheric convection, in which variables describe convective intensity and horizontal and vertical temperature differences~\cite{Lorenz1963DeterministicNF}. This system is given by the following set of ordinary differential equations:
\begin{align}
	&\frac{dx}{d t} = \tilde{\sigma} \left(y - x\right),\notag\\
	& \frac{dy}{d t} = x \left(\rho - z\right) - y,\notag\\
	& \frac{dz}{d t} = x y - \beta z,\label{eq54}		
\end{align}
where $\rho$, $\tilde{\sigma}$ and $\beta$ are the Prandtl number, Rayleigh number, and a geometric factor, respectively.

We take the parameters $\tilde{\sigma} = 10$, $\rho = 28$ and $\beta = 8/3$, under which the calculations in~\cite{Wang2022} were carried out. (Note that the work~\cite{Wang2022} incorrectly specifies the parameter $\sigma$ at which numerical calculations were performed: $\tilde{\sigma}=3$ is wrote instead $\tilde{\sigma}=10$). The max time is T = 20, and initial conditions are $x(0)=1$, $y(0)=1$, and $z(0)= 1$. 

\begin{figure}[t!]
	\centering
	\begin{subfigure}[t]{0.3\textwidth}
		\centering
		\includegraphics[width=\linewidth]{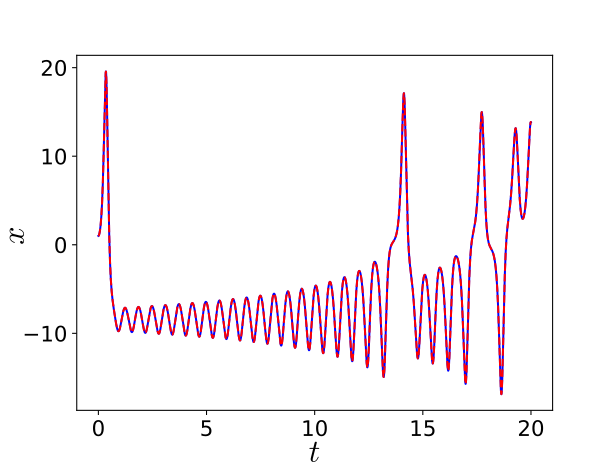}
		\caption{$x$}
	\end{subfigure}
	\begin{subfigure}[t]{0.3\textwidth}
		\centering
		\includegraphics[width=\linewidth]{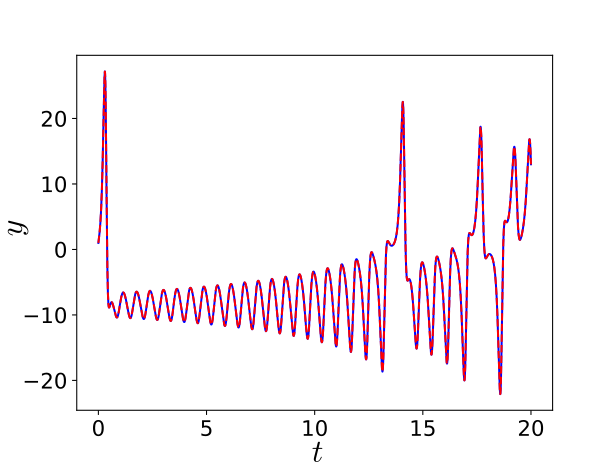}
		\caption{$y$}
	\end{subfigure}
	\begin{subfigure}[t]{0.3\textwidth}
		\centering
		\includegraphics[width=\linewidth]{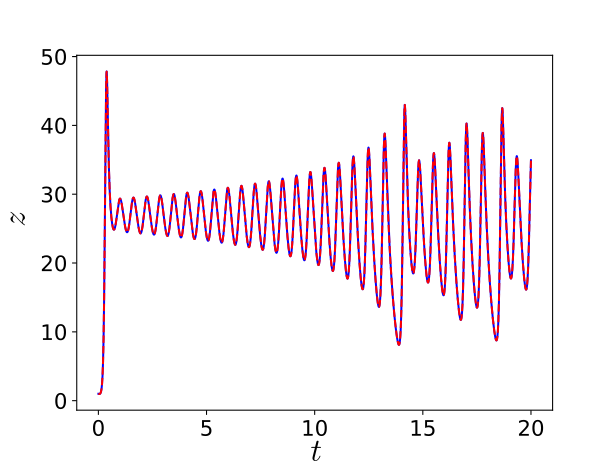}
		\caption{$z$}
	\end{subfigure}
	
	\caption{Lorentz system. (a), (b), and (c) are comparisons of the predicted (red dash lines) and reference solutions (blue solid lines) corresponding to $x$, $y$ and $z$, respectively. A number of neurons is $N=10000$, regime is PIDD initialization. The relative ${\mathbb L}_2$ errors are $\epsilon[{x}_{\bm \theta}, {x^{\text{ref}}}] = 2.0\times 10^{-4}$, $\epsilon[{y}_{\bm \theta}, {y^{\text{ref}}}] = 6.7\times 10^{-4}$, and $\epsilon[{z}_{\bm \theta}, {z^{\text{ref}}}] = 1.3\times 10^{-4}$.}
	\label{fig24}
\end{figure}

We used PIDD initialization. The reference solutions are obtained by using \verb*|scipy.integrate.odeint| with default settings. The data for initialization ($\left\{{x}_k\right\}_{k=0}^N$, $\left\{{y}_k\right\}_{k=0}^N$, $\left\{{z}_k\right\}_{k=0}^N$) for uniform distribution points $\left\{t_k\right\}_{k=0}^N$ were generated with \verb*|odeint| for $N=10000$  in ${\bf 0.014}$ seconds. Results of PIDD initialization for neural networks with $N=10000$ are following: the relative ${\mathbb L}_2$ errors are $\epsilon[{x}_{\bm \theta}, {x^{\text{ref}}}] = 2.0\times 10^{-4}$, $\epsilon[{y}_{\bm \theta}, {y^{\text{ref}}}] = 6.7\times 10^{-4}$, and $\epsilon[{z}_{\bm \theta}, {z^{\text{ref}}}] = 1.3\times 10^{-4}$. Such initialization took ${\bf 0.039}$ seconds of computing time of CPU. Results of such data-driven initialization are shown on Figure~\ref{fig24}.

\begin{figure}[t!]
	\centering
	\begin{subfigure}[t]{0.3\textwidth}
		\centering
		\includegraphics[width=\linewidth]{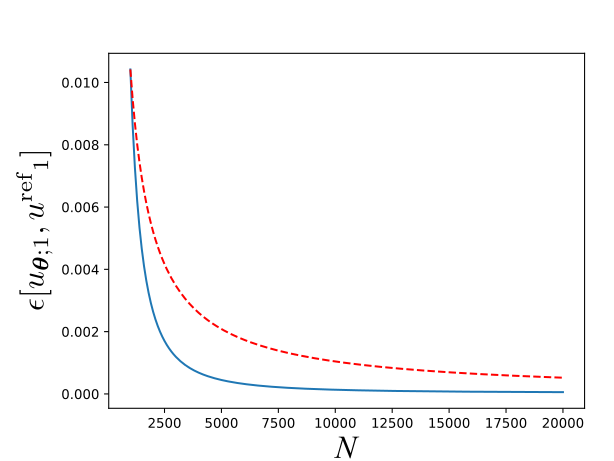}
		\caption{$x$}
	\end{subfigure}
	\begin{subfigure}[t]{0.3\textwidth}
		\centering
		\includegraphics[width=\linewidth]{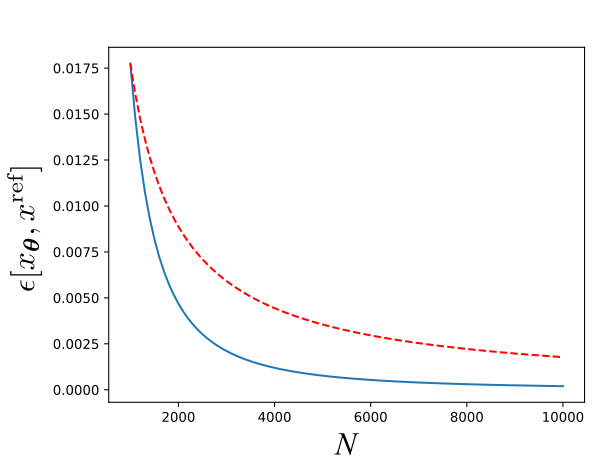}
		\caption{$y$}
	\end{subfigure}
	\caption{Blue solid lines of (a) and (b) are dependences of relative ${\mathbb L}_2$ error for $u$ of example 15 (harmonic oscillator) and $x$ for Example 16 (Lorentz system), respectively, on the number of neurons $N$ of neural network for the PIDD initialization. Red dashed lines of (a) and (b) are dependences $\alpha / N$.}
	\label{fig25}
\end{figure}

The dependence of relative ${\mathbb L}_2$ error $\epsilon[{x}_{\bm \theta}, {x^{\text{ref}}}]$ on the number of neurons $N$ of neural network are shown on the Figure~\ref{fig25}(b) by blue solid line. The red dashed line corresponds to dependence $\alpha_0 / N$, where $\alpha_0$ related to value $\epsilon[{x}_{\bm \theta}, {x^{\text{ref}}}]$ for $N=1000$. As in the case of the harmonic oscillator example, the accuracy of the neural network increases with the increasing number of neurons faster than the accuracy of the method of integrating a differential equation based on rectangle methods.


\subsection{Generalization}
In this subsection, we check how well a trained neural network, created on the basis of a neural network with one hidden layer, can predict solutions of a differential equation at parameters, at which it was not trained. I.e., we investigate the generalizing properties of such a neural network. For simplicity and clarity, we assume that differential equations depend on only one parameter. Thus, the solutions of differential equations should depend on two quantities: the coordinates and the specified parameter.

In our experiments on generalization properties of NN for the solving of ODE we used a network architecture, called separable PINN (SPINN), which facilitates significant forward-mode automatic differentiation (AD) for more efficient computation~\cite{cho2023,eskin2024SPINN}. SPINN consists of $d$ body-networks (which are usually presented by MLPs), each of which takes an individual one-dimensional coordinate component as an input from $d$-dimensional space. Each body-network ${\vec f}^{{\bm \theta}_i}:\mathbb{R}\rightarrow \mathbb{R}^{r}$ (parameterized
by ${{\bm \theta}_i}$) is a vector-valued function which transforms the coordinates of $i$th axis into a $r$-dimensional feature representation (functional basis, or set of eigenfunctions).

The predictions of SPINN are computed by basis functions merging which is defined by following
\begin{equation}
	u_{\bm \theta}({\vec x}) = \sum\limits_{j=1}^{r} \prod\limits_{i=1}^{d} {f}^{{\bm \theta}_i}_{j} (x_i),\label{eq49}
\end{equation}
where ${u}^{{\bm \theta}}:\mathbb{R}^{d}\rightarrow \mathbb{R}$ is the predicted solution function, {${\vec x} \in \mathbb{R}^{d}$ are input coordinates, $x_i$ denotes $i$th input coordinate}, and ${f}^{{\bm \theta}_i}_{j}$ is the $j$th element of basis ${\vec f}^{{\bm \theta}_i}$ $\left({\vec f}^{{\bm \theta}_i} = \left({f}^{{\bm \theta}_i}_{1},{f}^{{\bm \theta}_i}_{2},\dots,{f}^{{\bm \theta}_i}_{r}\right)\right)$. {SPINN possesses a universal approximation property, which is peculiar to neural networks. The proof of the presence of this capability in such neural networks is provided by the authors of SPINN  in the seminal work~\cite{cho2023}. This property allows us to use SPINN to solve any physical problem described by a PDE.}

{\bf Suggestion 8.} Every body-networks of SPINN can be presented by simple MLP, given by equation (\ref{eq3}). We use the following SPINN representation of the solution of ODE, depending on the coordinate $x$ and on the parameter $q$
\begin{equation}
	u_{\bm \theta}({x},q) = \sum\limits_{j=1}^{r} U_j (x) Q_j(q) .\label{eq56}
\end{equation}
In the examples discussed above (Examples 1,3, and 6), the parameter $q$ is the only one. So for such solutions, the function $Q$ can be equated to 1, and $r = 1$. Based on this, we can require that the value of $Q_j$ be kept equal to 1 in order to prepare an increase in the number of parameters considered $q$. We naturally come to the following method of constricting the generalized solution
\begin{enumerate}
	\item Let we have homogeneous set $\left\{{q}_{i}\right\}^{N_{q}}_{i=1}$ of parameter $q$ with step $\Delta q$ taken from domain $[q_{\rm min}, q_{\rm max}]$;
	\item For each parameter $q_i$, the solution $U_i(x)$ for ODE is found in domain [$0, X$];
	\item The generalized solution of ODE is
	\begin{equation}
		u_{\bm \theta}({x},q) = \sum\limits_{j=1}^{N_q} U_j (x) Q_j(q).\label{eq57}
	\end{equation}	
\end{enumerate}
Here the functions $Q_j(q)$ should be designed in such a way that $u_{\bm \theta}({x},q_n) \approx U_n (x)$. This means that $Q_n(q)$ should give 1 for $q=q_n$ and values close to zero for other $q$. According to the SPINN approach, the neural network for solving $u_{\bm\theta}({x},q)$ is two neural networks $U$ and $Q$ with $N_q$ outputs.

{\bf Suggestion 9: Initialization of functions $U$ and $Q$.} We initialize the parameters of neural network $Q$ with $N_q$ output as follows:
\begin{enumerate}
	\item Weights of the hidden layer are taken $2 \Delta \zeta / \Delta q$ (${W}^{(1)}_k = 2 \Delta \zeta / \Delta q$) and biases $b^{(1)}_k=-(2 k - 1) \Delta \zeta$ for the $k$th neuron.
	\item Activation function is the sum of two sigmoids as follows: $\sigma(x) - \sigma(x - 2 \Delta \zeta)$. The maximum of this activation function is achieved when $x = 2 \Delta \zeta$.
	\item The biases of output layer is taken zero ($b^{(2)}_j=0$).
	\item Additionally, we refine the weights by solving the following equation using the least squares method
	\[{\vec W}^{(2)}{\vec a} = {\vec I},\]
	where vector ${\vec W}^{(2)}$ is presented by weights of second layer of $Q$ (${W}^{(2)}_{kj}$),
	\[{\vec a} = \sigma\left({\vec W}^{(1){\rm T}}_k {\vec q} + {\vec b}^{(1)}_k\right) - \sigma\left({\vec W}^{(1){\rm T}}_k {\vec q} + {\vec b}^{(1)}_k  - 2 \Delta \zeta\right),\] 
	${\vec q}$ is set $\left\{{q}_{i}\right\}^{N_{q}}_{i=1}$, ${\vec I}$ is identity matrix of size $N_q$.
	\item Neural network $U$ is initialized by Algorithm~\ref{alg3}, but to obtain the weights of the output layer, one needs to change the formulas to the following:
	\begin{align}
		&{W}^{(2)}_{kj} \gets - \dfrac{\Delta x}{2\Delta \zeta} \dfrac{\mathcal{N}[{u}(x_{k},q_j), x_{k}]}{\kappa_k}\notag\\
		&{b}^{(2)}_j \gets {u}(x_0,q_j) - \sum\limits_{k = 0}^{N-1} {W}^{(2)}_{kj}\sigma\left({W}^{(1)}_k x_0 + {b}^{(1)}_k\right).\notag
	\end{align}	
\end{enumerate}
As a result, we have the following Algorithm~\ref{alg4} of physics-informed data-driven initialization of neural networks with generalization properties.
\begin{algorithm}[h!]
	\caption{Physics-informed data-driven initialization of neural networks with generalization properties}\label{alg4}
	\KwData{$\left\{{u}_{kj}\right\}$ for $\left\{x_k\right\}_{k=0}^{N_x-1}$ (uniform grid with a step $\Delta x \gets X / N_x$) and $\left\{q_j\right\}_{j=0}^{N_q-1}$ }
	\KwResult{Initialized NNs $U$ and $Q$ for solution ${u}_{\bm \theta}$. $U$ and $Q$ consist of $N_x$ and $N_q$ neurons on hidden layer and $N_q$ ouputs}
	${\Delta\zeta} \gets \ln\left(2 + \sqrt{3}\right) / 2$\;
	Initialization for $U$\;
	\For{$k=0,\dots, N_x-1$}{
		${W}^{(1)}_k \gets 2 \Delta \zeta / \Delta x$\;
		${b}^{(1)}_k \gets -2 k \Delta \zeta$\;
		$\kappa_k \gets \sum\limits_{m = k-L}^{k+L} \sigma'\left({2\Delta \zeta} \left[k - m)\right]\right)$\;
		${W}^{(2)}_{kj} \gets - \dfrac{\Delta x}{2\Delta \zeta} \dfrac{\mathcal{N}[{u}_{kj}, x_{k}]}{\kappa_k}$\;}
	${b}^{(2)}_j \gets {u}_{0j} - \sum\limits_{k = 0}^{N_x-1} {W}^{(2)}_{kj}\sigma\left({W}^{(1)}_k x_0 + {b}^{(1)}_k\right)$\;
	Initialization for $Q$\;
	\For{$k=0,\dots, N_q-1$}{
		${b}^{(1)}_k \gets - (2k - 1) \Delta \zeta$\;
		${W}^{(1)}_k \gets 2 \Delta \zeta / \Delta q$\;}
	Weights ${W}^{(2)}_{kj}$ are determined by solving the following equation using the least squares method
	\begin{align}
		&{\vec W}^{(2)}{\vec a} = {\vec I},\notag\\
		&{\vec a} = \sigma\left({\vec W}^{(1){\rm T}}_k {\vec q} + {\vec b}^{(1)}_k\right) - \sigma\left({\vec W}^{(1){\rm T}}_k {\vec q} + {\vec b}^{(1)}_k  - 2 \Delta \zeta\right);\notag
	\end{align}
	${b}^{(2)}_j \gets 0$.
\end{algorithm}

For the evaluation of generalization properties we use the measure of generalization $\mu$ which is given in appendix~\ref{appB}.

\subsection{Numerical experiments}

We investigate the generalization properties of the described neural network using the example of linear and nonlinear ODEs.

\subsubsection{Example 17. Harmonic Oscillator}

Consider the harmonic oscillator, which is governed by a system of two hidden-order equations for the $t\in [0, T]$
\begin{eqnarray}
	&& \frac{d u_1}{d t} = u_2,\quad \frac{d u_2}{d t} = - \omega^2 u_1, \quad t \in [0,T]\label{eq58}\\
	&& u_1(0) = 1,\quad u_2(0) = 0 \label{eq59},
\end{eqnarray}
where $\omega$ is the frequency of the considered system. We used the following parameters: $T=12$, $\omega \in [1,10]$. The exact analytical solution of this problem is $u_1^{\text{ref}} = \cos\left(\omega t\right)$, $u_2^{\text{ref}} = - \sin\left(\omega t\right)$. Consider how the generalization parameters of the solution provided by a neural network depend on the number of neurons $N_q$ in the hidden layer for the $Q$ network.

In our experiments, we used physics-informed data-driven initialization of neural networks with generalization properties, which is described above (Algorithm~\ref{alg4}). The data $\left\{{u}_{kj}\right\}$ for uniform distribution points $\left\{x_k\right\}_{k=0}^{N_x-1}$ and $\left\{\omega_j\right\}_{j=0}^{N_q-1}$.

\begin{figure}[t!]
	\centering
	\begin{subfigure}[t]{0.6\textwidth}
		\centering
		\includegraphics[width=\linewidth]{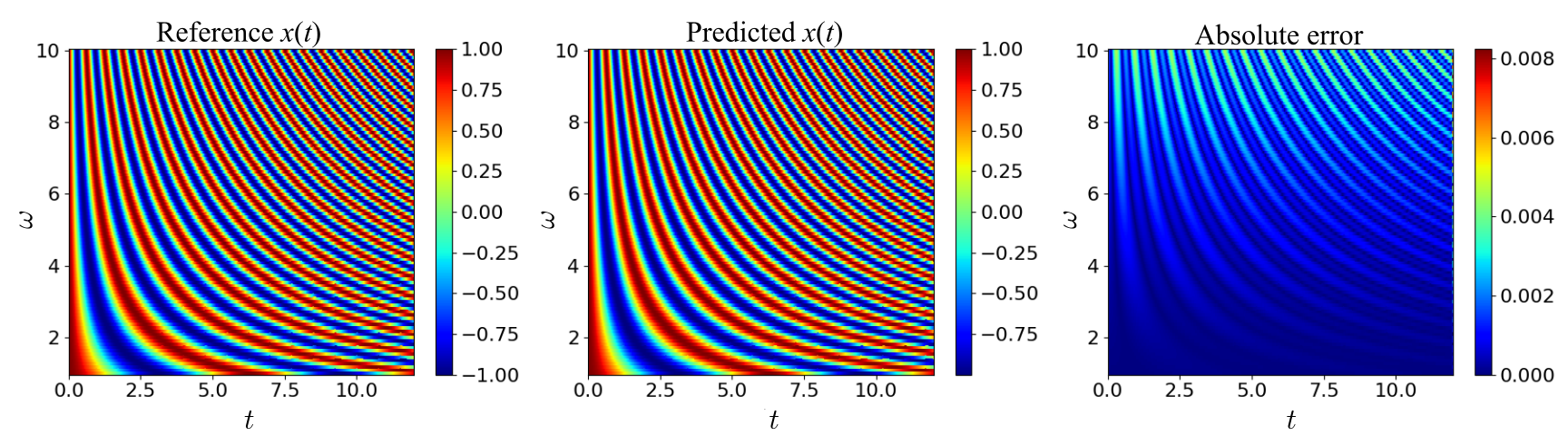}
		\caption{}
	\end{subfigure}
	\begin{subfigure}[t]{0.6\textwidth}
		\centering
		\includegraphics[width=\linewidth]{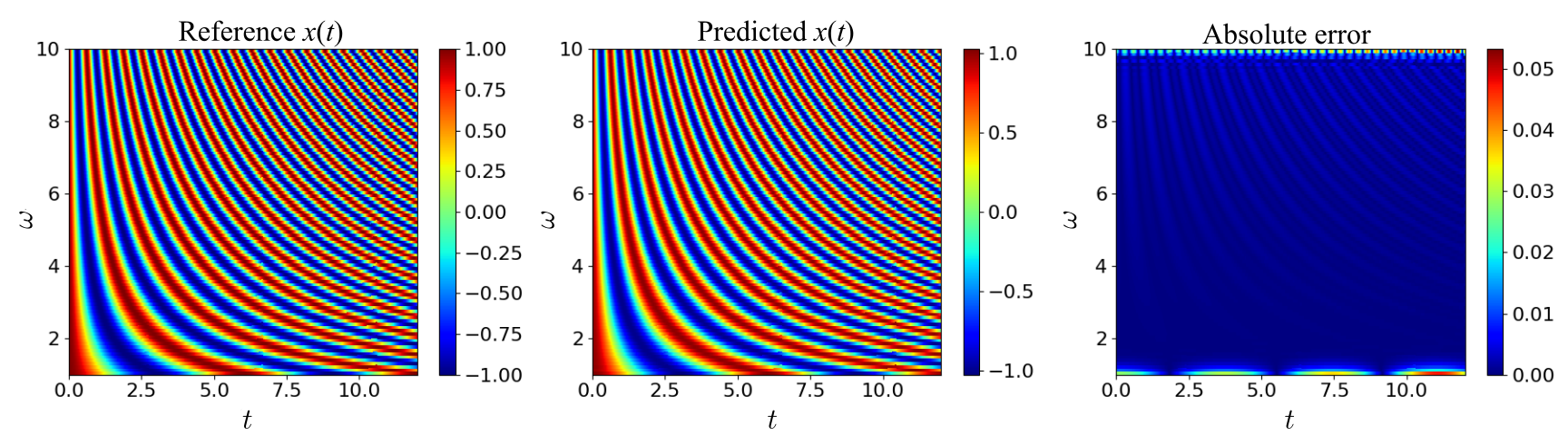}
		\caption{}
	\end{subfigure}
	\caption{Harmonic Oscillator. (a) the left panel is the reference solution, (a) the middle panel is a prediction of a trained physics-informed neural network based on algorithm~\ref{alg4}, (a) the right panel is absolute difference between the reference solution and the predicted solution for the frequencies $\left\{\omega_j\right\}_{j=0}^{N_q}$. (b) are same for the frequencies $\left\{\omega_{j+1/2}\right\}_{j=0}^{N_q-1}$ at which training was not carried out. The relative error $\epsilon_1[{u}_{{\bm \theta};1}, {u^{\text{ref}}_1}] = 5.14\times 10^{-3}$ for (a) and $\epsilon_{1/2}[{u}_{{\bm \theta};1}, {u^{\text{ref}}_1}] = 2.42\times 10^{-2}$ for (b). $N_q = 100$ and $N_x = 2000$. The measure of generalization $\mu$ is $0.21$.}
	\label{fig26_1}
\end{figure}

\begin{figure}[t!]
	\centering
	\begin{subfigure}[t]{0.3\textwidth}
		\centering
		\includegraphics[width=\linewidth]{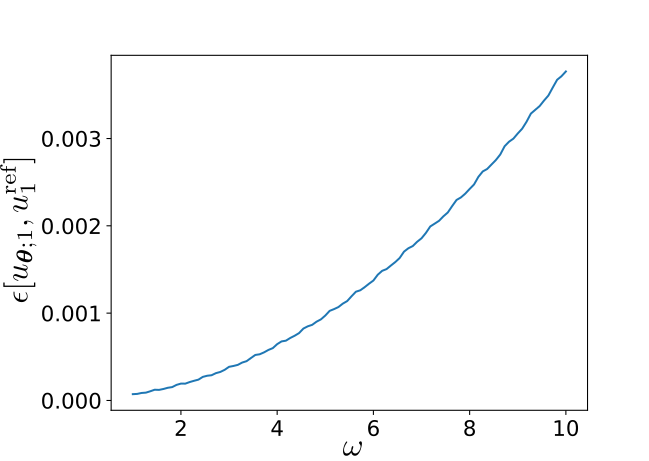}
		\caption{}
	\end{subfigure}
	\begin{subfigure}[t]{0.3\textwidth}
		\centering
		\includegraphics[width=\linewidth]{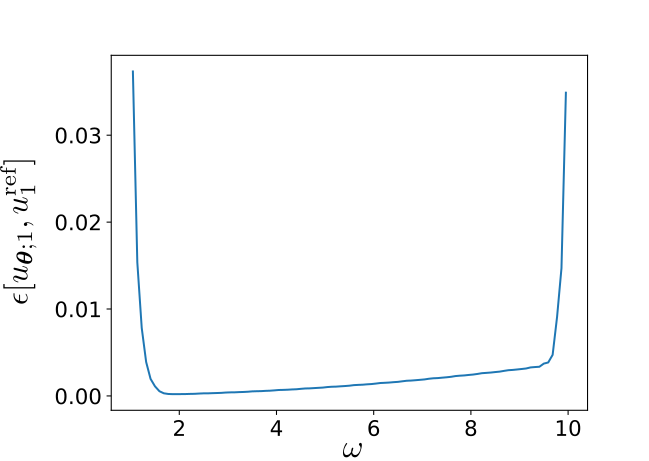}
		\caption{}
	\end{subfigure}
	\caption{Harmonic Oscillator. (a) and (b) are dependences of relative ${\mathbb L}_2$ error for $u_1$ on frequencies for sets of frequencies $\left\{\omega_j\right\}_{j=0}^{N_q}$ and $\left\{\omega_{j+1/2}\right\}_{j=0}^{N_q-1}$.}
	\label{fig27_1}
\end{figure}

\begin{figure}[t!]
	\centering
	\includegraphics[width=0.3\linewidth]{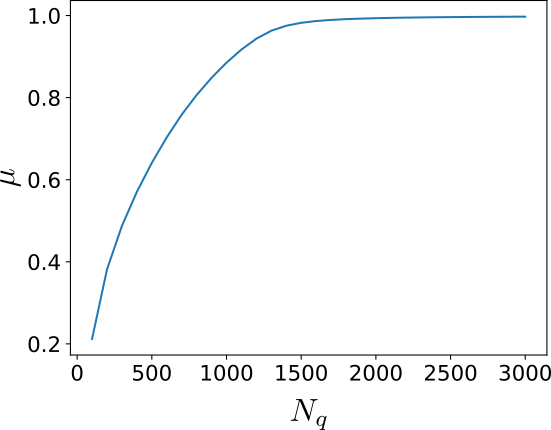}
	\caption{Harmonic Oscillator. Dependence of measure of generalization properties $\mu$ on size of set of frequencies $N_q$.}
	\label{fig28}
\end{figure}


Results of initialization by algorithm~\ref{alg4} for neural network with $N_x=2000$ and $N_q=100$ are presented on Figure~\ref{fig26_1} for the sets of frequencies $\left\{\omega_j\right\}_{j=0}^{N_q}$ of training domain and $\left\{\omega_{j+1/2}\right\}_{j=0}^{N_q-1}$ at which training was not carried out. The relative errors $\epsilon_1[{u}_{{\bm \theta};1}, {u^{\text{ref}}_1}] = 5.14\times 10^{-3}$, $\epsilon_{1/2}[{u}_{{\bm \theta};1}, {u^{\text{ref}}_1}] = 2.42\times 10^{-2}$ and the measure of generalization $\mu = 0.21$. Figure~\ref{fig27_1} shows the dependences of relative ${\mathbb L}_2$ error for $u_1$ on frequency for these sets. It can be seen, the main discrepancies for these two sets of frequencies are observed at the boundary of the frequency domain. Figure\ref{fig28} shows the dependence of the measure of generalizing properties $\mu$ on the size of the frequency set $N_q$, from which it can be seen that with increasing of $N_q$, the generalization measure quickly reaches values close to 1.

\subsubsection{Example 18. Lorentz system}

Consider behaviour of generalization properties of neural networks under PIDD initialization for the system of equations~(\ref{eq54}) depending on the Prandtl number $\rho$.

\begin{figure}[t!]
	\centering
	\begin{subfigure}[t]{0.6\textwidth}
		\centering
		\includegraphics[width=\linewidth]{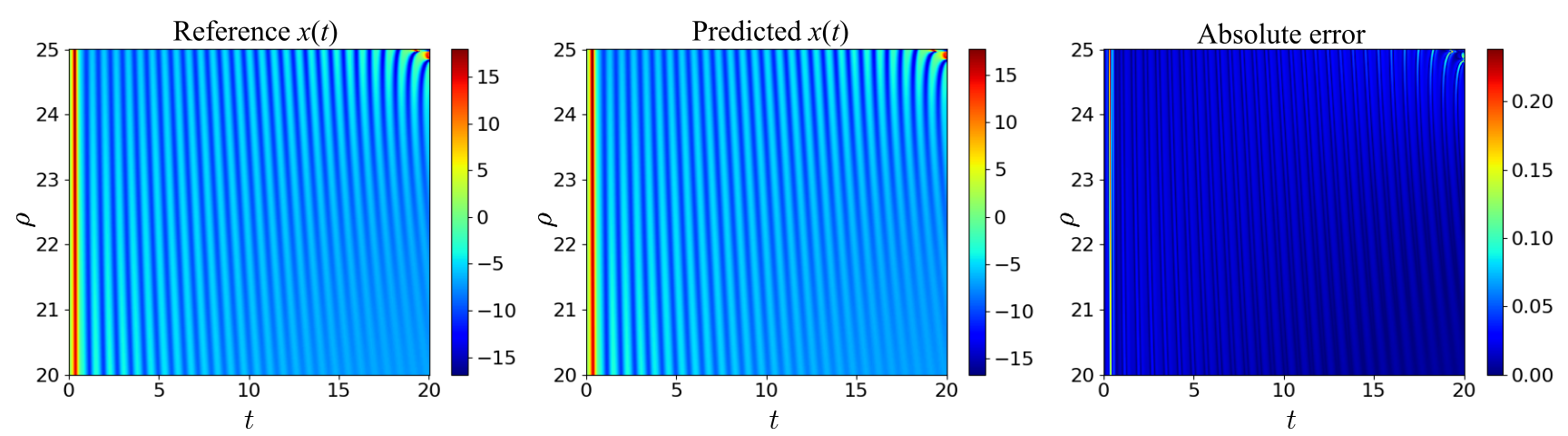}
		\caption{}
	\end{subfigure}
	\begin{subfigure}[t]{0.6\textwidth}
		\centering
		\includegraphics[width=\linewidth]{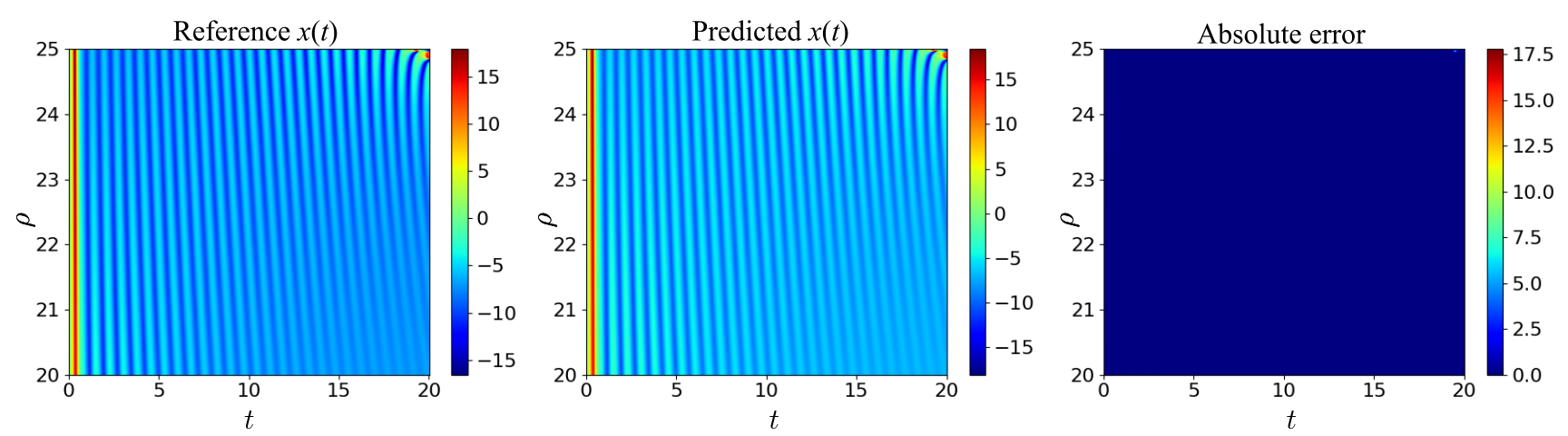}
		\caption{}
	\end{subfigure}
	\begin{subfigure}[t]{0.6\textwidth}
		\centering
		\includegraphics[width=\linewidth]{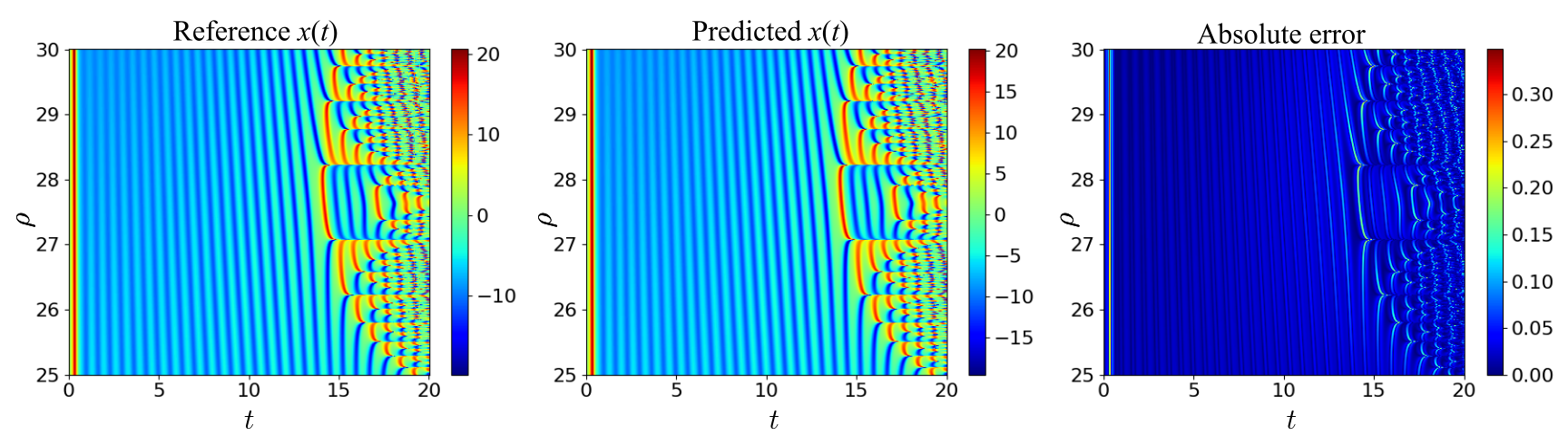}
		\caption{}
	\end{subfigure}
	\begin{subfigure}[t]{0.6\textwidth}
		\centering
		\includegraphics[width=\linewidth]{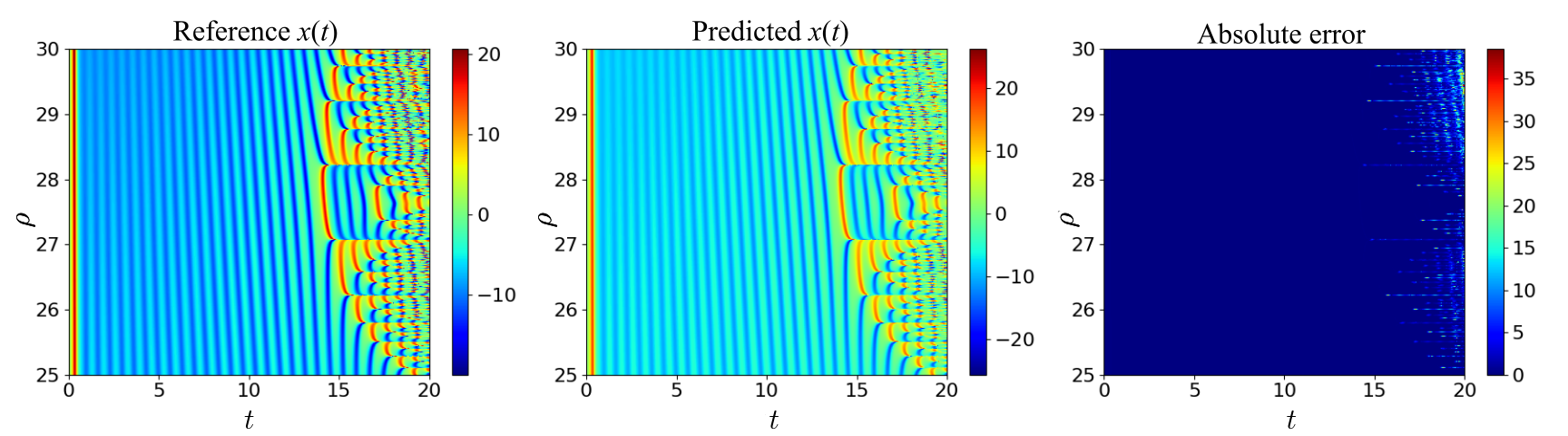}
		\caption{}
	\end{subfigure}
	\caption{Lorentz system. (a) the left panel is the reference solution, (a) the middle panel is a prediction of a trained physics-informed neural network based on algorithm~\ref{alg4}, (a) the right panel is the absolute difference between the reference solution and the predicted solution for the parameters $\left\{\rho_j\right\}_{j=0}^{N_q}$ in domain $[20,25]$. (b) are the same for the Lorentz system $\left\{\rho_j{j+1/2}\right\}_{j=0}^{N_q-1}$ at which training was not carried out in domain $[20,25]$. Panels (c) and (d) are same at the (a) and (b), respectively, for the parameters $\rho \in [25,30]$. The relative errors are $\epsilon_1[{x}_{\bm \theta}, {x^{\text{ref}}}] = 2.01\times 10^{-3}$ (a), $\epsilon_{1/2}[{x}_{\bm \theta}, {x^{\text{ref}}}] = 6.9\times 10^{-3}$ (the measure of generalization $\mu = 0.29$) for  $\rho\in [20,25]$, and $\epsilon_1[{x}_{\bm \theta}, {x^{\text{ref}}}] = 3.02\times 10^{-3}$ (a), $\epsilon_{1/2}[{x}_{\bm \theta}, {x^{\text{ref}}}] = 1.27\times 10^{-1}$ (the measure of generalization $\mu = 0.023$) for  $\rho\in [25,30]$. $N_q = 1000$ and $N_t = 2000$.}
	\label{fig29}
\end{figure}

\begin{figure}[t!]
	\centering
	\begin{subfigure}[t]{0.3\textwidth}
		\centering
		\includegraphics[width=\linewidth]{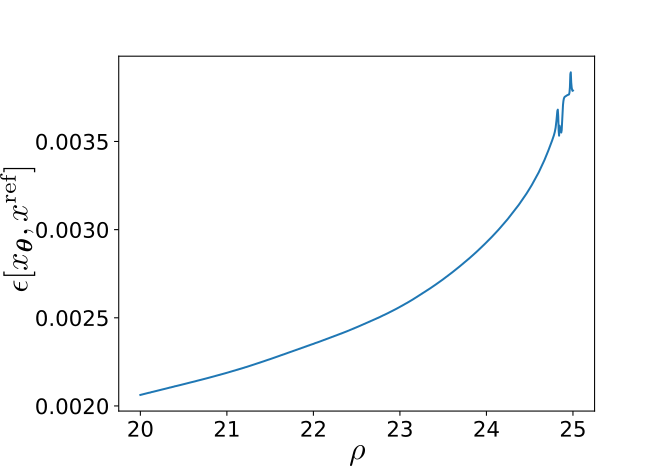}
		\caption{}
	\end{subfigure}
	\begin{subfigure}[t]{0.3\textwidth}
		\centering
		\includegraphics[width=\linewidth]{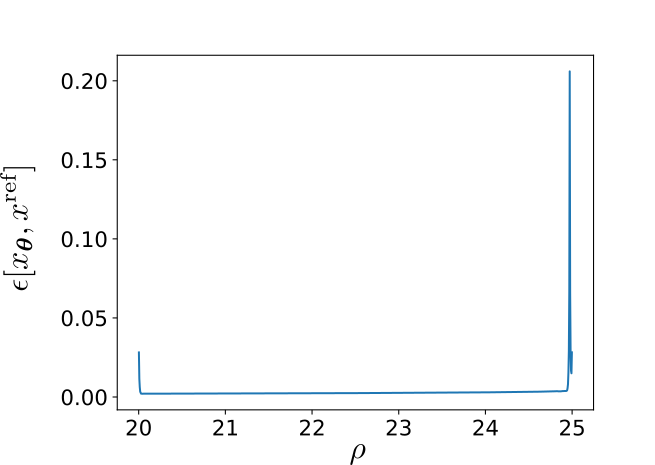}
		\caption{}
	\end{subfigure}\\
	\begin{subfigure}[t]{0.3\textwidth}
		\centering
		\includegraphics[width=\linewidth]{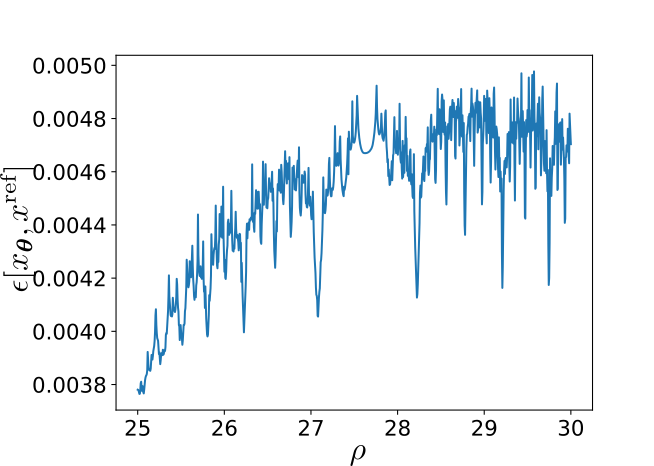}
		\caption{}
	\end{subfigure}
	\begin{subfigure}[t]{0.3\textwidth}
		\centering
		\includegraphics[width=\linewidth]{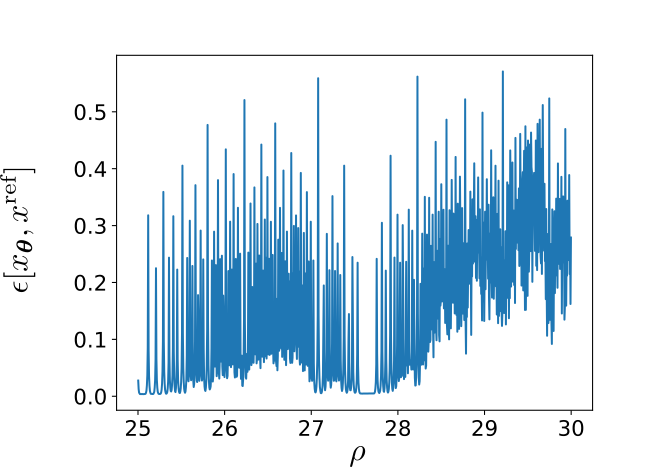}
		\caption{}
	\end{subfigure}
	\caption{Lorentz system. (a) and (b) are dependences of relative ${\mathbb L}_2$ error for $x(t)$ on $\rho$ for sets of parameters $\left\{\rho_j\right\}_{j=0}^{N_q}$ and $\left\{\rho_j{j+1/2}\right\}_{j=0}^{N_q-1}$ in domain $\rho\in[20,25]$. (c) and (d) are the same in domain $\rho\in[25,30]$.}
	\label{fig30}
\end{figure}

\begin{figure}[t!]
	\centering
	\begin{subfigure}[t]{0.3\textwidth}
		\centering
		\includegraphics[width=\linewidth]{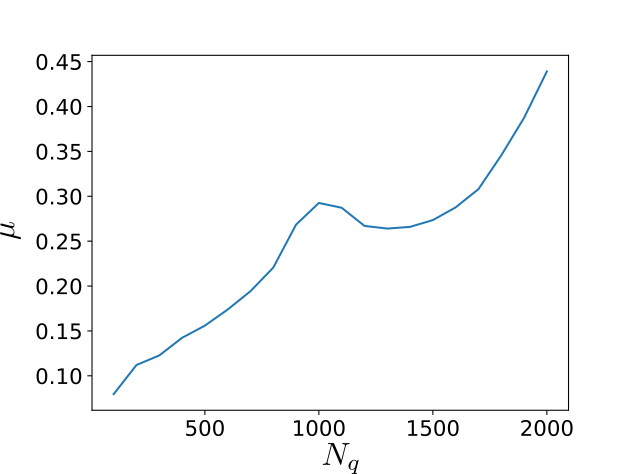}
		\caption{}
	\end{subfigure}
	\begin{subfigure}[t]{0.3\textwidth}
		\centering
		\includegraphics[width=\linewidth]{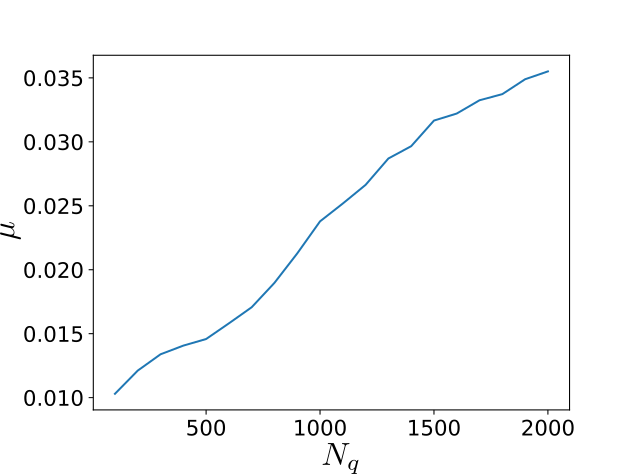}
		\caption{}
	\end{subfigure}
	\caption{Lorentz system. Dependences of measure of generalization properties $\mu$ on size of sets of the Prandtl numbers $\rho$ ($N_q$). (a) is for $\rho\in[20,25]$ and (b) is for $\rho\in[25,30]$.}
	\label{fig31}
\end{figure}

Results of initialization by algorithm~\ref{alg4} for neural network with $N_t=2000$ and $N_q=1000$ are presented on the Figure~\ref{fig29} for the sets of parameter $\left\{\rho_j\right\}_{j=0}^{N_q}$ of training domain and $\left\{\rho_{j+1/2}\right\}_{j=0}^{N_q-1}$ at which training was not carried out in domains: $\rho\in [20,25]$ for (a) and (b), $\rho \in [25,30]$  for (c) and (d). The relative errors are $\epsilon_1[{x}_{\bm \theta}, {x^{\text{ref}}}] = 2.01\times 10^{-3}$ (a), $\epsilon_{1/2}[{x}_{\bm \theta}, {x^{\text{ref}}}] = 6.9\times 10^{-3}$ (the measure of generalization $\mu = 0.29$) for  $\rho\in [20,25]$, and $\epsilon_1[{x}_{\bm \theta}, {x^{\text{ref}}}] = 3.02\times 10^{-3}$ (a), $\epsilon_{1/2}[{x}_{\bm \theta}, {x^{\text{ref}}}] = 1.27\times 10^{-1}$ (the measure of generalization $\mu = 0.023$) for  $\rho\in [25,30]$. The Figure~\ref{fig30} shows the dependences of relative ${\mathbb L}_2$ error for $x(t)$ on $\rho$ for this sets. It can be seen, the main discrepancies for these two sets of $\rho$ are observed at the boundary of the $\rho$ domain for $\rho\in [20,25]$ and multiple spikes for $\rho\in [25,30]$. Such instability for $\rho\in [25,30]$ is due to the fact that at these parameters, the system of equations under consideration begins to exhibit chaotic behaviour. This chaotic behaviour is shown in Figures~\ref{fig29} (c) and (d) in the form of fractal structures for $t>15$. Figure~\ref{fig31} shows the dependences of the measure of generalizing properties $\mu$ on $N_q$, the size of the sets, from which it can be seen that with increasing of $N_q$. Note that despite the increasing of the generalizing abilities of the network for  $\rho\in [25,30]$, the value of the generalization measure remains relatively small. I.e., it can be concluded that in the areas of chaotic behaviour of the system of differential equations, the proposed neural network demonstrates relatively weak generalizing abilities. This was to be expected, because no neural network at the moment can predict the detailed fractal structure of those functions that possess this structure.

\subsubsection{Discussion}
In the previous two experiments, it was demonstrated that the proposed neural network, in the case of PIDD initialization, has pronounced generalizing properties from such parameters of a system of differential equations in which this system remains regular (not chaotic). Note that this initialization method is easily parallelized. We can train-initialize functions $U$ for different values of parameters $q$ independently of each other. Despite the power of PIDD initialization, its operation requires solutions of the ODE system obtained numerically or based on an analytical solution. It would be better to exclude the leverage of third-party solvers in the construction of the solution, i.e. to remain within the framework of the PINN approach.

\subsubsection{Neuron-by-neuron training}

{\bf Suggestion 10: Neuron-by-neuron training (NbN).} A neural network initialized using Algorithm~\ref{alg2} can be iteratively improved neuron by neuron based on differential equations and predictions of this neural network. Such training is based on the following steps:
\begin{enumerate}
	\item Initialize NNs ${\vec u}_{\bm \theta}$ with Algorithm~\ref{alg2};
	\item Additionally, we iteratively change the weights ${W}^{(2)}_{k}$ using the Euler integration method based on the rectangle method: $x_{k} = k \Delta x$, $\hat{\vec u}_{0} = {\vec u}_0$, $\hat{\vec u}_{k} = \hat{\vec u}_{k-1} + \mathcal{N} [\hat{\vec u}_{k-1}, x_{k-1}]\Delta x$,  ${W}^{(2)}_{k} =  - 2 \Delta x \mathcal{N} [\hat{\vec u}_{k}, x_{k}] / \Delta \zeta$. This step is optional, but it additionally accelerates the convergence of the method.
	\item Next, the weights ${W}^{(2)}_{k}$ change based on the predictions of the neural network. I.e., the neural network is self-improving. ${W}^{(2)}_{k} =  - 2 \Delta x \mathcal{N}_l [{\vec u}_{\bm \theta}(x_{k}), x_{k}] / \Delta \zeta$.
\end{enumerate}
As a result, we have the following neural network training Algorithm~\ref{alg5}.
\begin{algorithm}[h!]
	\caption{Neuron-by-neuron training}\label{alg5}
	\KwData{---}
	\KwResult{Trained for the $E$ epoch NNs ${\vec u}_{\bm \theta}$, each of which consists of $N$ neurons on hidden layer}
	$\Delta x \gets X / N$\;
	${\Delta\zeta} \gets \ln\left(2 + \sqrt{3}\right) / 2$\;
	Initialization of ${\vec u}_{\bm \theta}$ with following steps\;
	\For{$k=0,\dots, N-1$}{
		${W}^{(1)}_k \gets 2 \Delta \zeta / \Delta x$\;
		${b}^{(1)}_k \gets -2 k \Delta \zeta$\;
		${W}^{(2)}_k \gets - 2 \Delta x \mathcal{N}_l [{\vec u}(0), x_{k}] / \Delta \zeta$\;}
	${b}^{(2)}_0 \gets {u_l}(0)$\;
	Iteratively change the weights ${W}^{(2)}_{k}$ of ${\vec u}_{\bm \theta}$ using the Euler integration method with the following steps (optional procedure)\;
	$\hat{\vec u}_{0} = {\vec u}_0$\;
	\For{$k=0,\dots, N-1$}{
		$\hat{\vec u}_{k} = \hat{\vec u}_{k-1} + \mathcal{N} [\hat{\vec u}_{k-1}, x_{k-1}]\Delta x$\;
		${W}^{(2)}_{k} =  - 2 \Delta x \mathcal{N}_l [\hat{\vec u}_{k}, x_{k}] / \Delta \zeta$\;}
	Neuron-by-neuron improving weights ${W}^{(2)}_{k}$ with following steps\;
	\For{$e=0,\dots, E-1$}{
		\For{$k=0,\dots, N-1$}{
			${W}^{(2)}_{k} =  - 2 \Delta x \mathcal{N}_l [{\vec u}_{\bm \theta}(x_{k}), x_{k}] / \Delta \zeta$.}
	}
	In our experiments, it was enough to use 3 epochs ($E=3$) of training to obtain high-precision results.
\end{algorithm}

\subsection{Numerical experiments}

We investigate the generalization properties of the described neural network using the example of linear and nonlinear ODEs.

\subsubsection{Example 19. Harmonic Oscillator}
Consider system of equations~(\ref{eq9})--(\ref{eq11}) with same parameters. The whole time domain $[0;100]$ was split into 10 disjoint equivalent time windows of size $\Delta t = 10$. In our experiments for each
time window we took 2 neural networks (which correspond to $u_1$, and $u_2$ values) with 10000 neurons on the hidden layer. The training was carried out on GPU and took 3 iterations of NbN training and $208.8$ seconds of computational time.

\begin{figure}[t!]
	\centering
	\begin{subfigure}[t]{0.3\textwidth}
		\centering
		\includegraphics[width=\linewidth]{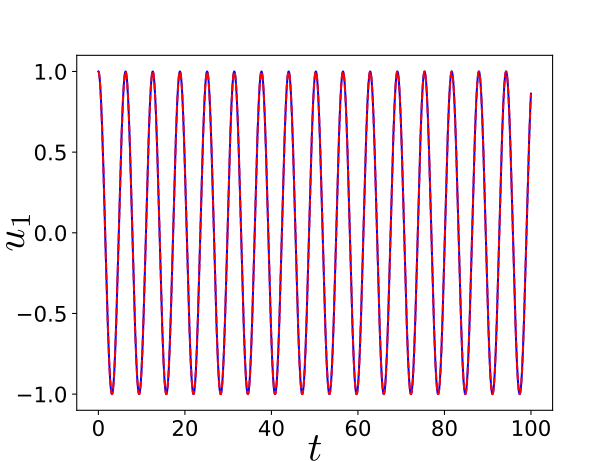}
		\caption{}
	\end{subfigure}
	\begin{subfigure}[t]{0.34\textwidth}
		\centering
		\includegraphics[width=\linewidth]{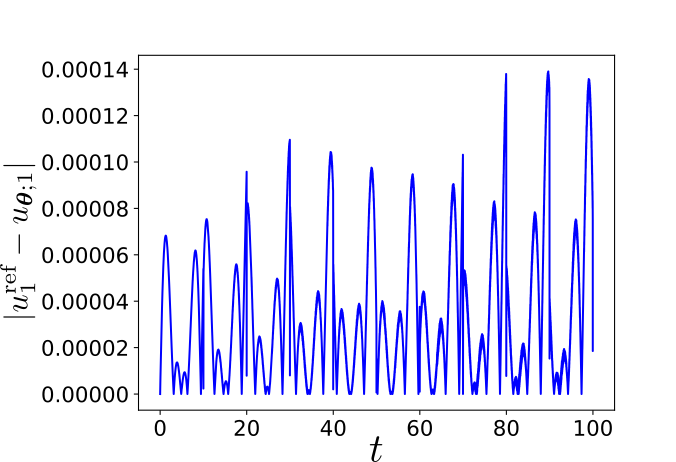}
		\caption{}
	\end{subfigure}
	
	\caption{Harmonic oscillator. (a) is a comparison of the predicted (red dash lines) and reference solutions (blue solid lines) corresponding to $u_1$. (b) is absolute error $|u^{\text{ref}}_1 - u_{{\bm \theta};1}|$. The relative ${\mathbb L}_2$ errors are $\epsilon[{u}_{{\bm \theta};1}, {u^{\text{ref}}_1}] = 6.56\times 10^{-5}$, and $\epsilon[{u}_{{\bm \theta};2}, {u^{\text{ref}}_2}] = 6.37\times 10^{-5}$.}
	\label{fig32}
\end{figure}

Figure~\ref{fig32} shows the predicted $u_1$ and the reference $u_1$ was obtained by \verb*|odeint| solver of \verb*|scipy.integrate| library. The relative ${\mathbb L}_2$ errors are $\epsilon[{u}_{{\bm \theta};1}, {u^{\text{ref}}_1}] = 6.56\times 10^{-5}$, and $\epsilon[{u}_{{\bm \theta};2}, {u^{\text{ref}}_2}] = 6.37\times 10^{-5}$. 

\subsubsection{Example 20. Lorentz system}
Consider a system of equations~(\ref{eq54}) with the same parameters. The whole time domain $[0;20]$ was split into 20 disjoint equivalent time windows of size $\Delta t = 1$. In our experiments for each
time window we took 3 neural networks (which correspond to $x$, $y$, and $z$ values) with 10000 neurons on the hidden layer.  The training was carried out on GPU and took 3 iterations of NbN training and $545.02$ seconds of computational time.

\begin{figure}[t!]
	\centering
	\begin{subfigure}[t]{0.3\textwidth}
		\centering
		\includegraphics[width=\linewidth]{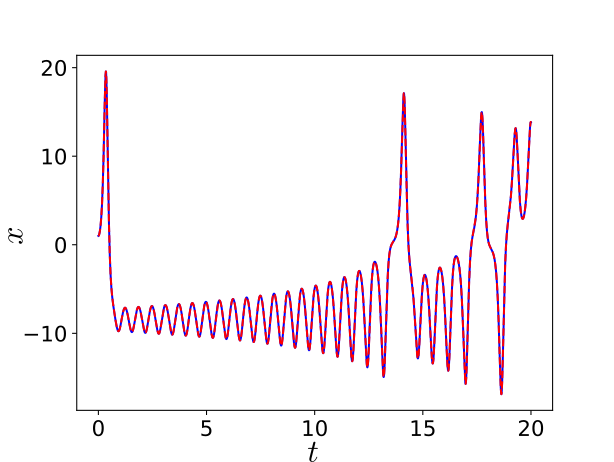}
		\caption{}
	\end{subfigure}
	\begin{subfigure}[t]{0.3\textwidth}
		\centering
		\includegraphics[width=\linewidth]{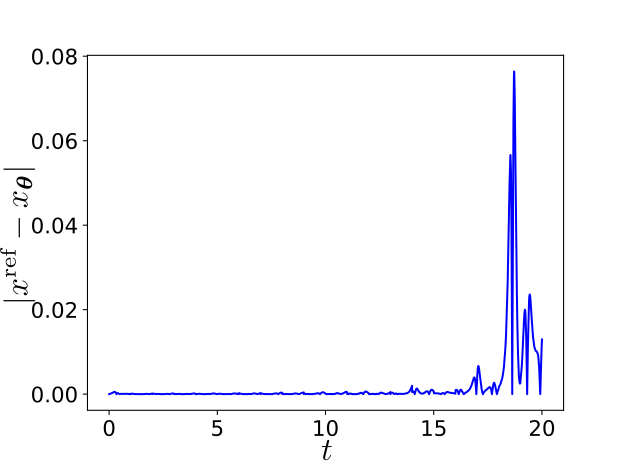}
		\caption{}
	\end{subfigure}
	
	\caption{Lorentz system. (a) is a comparison of the predicted (red dash lines) and reference solutions (blue solid lines) corresponding to $x$. (b) is absolute error $|x^{\text{ref}} - x_{\bm \theta}|$. The relative ${\mathbb L}_2$ errors are $\epsilon[{x}_{\bm \theta}, {x^{\text{ref}}}] = 9.8\times 10^{-4}$, $\epsilon[{y}_{\bm \theta}, {y^{\text{ref}}}] = 1.4\times 10^{-3}$, and $\epsilon[{z}_{\bm \theta}, {z^{\text{ref}}}] = 6.0\times 10^{-4}$.}
	\label{fig33}
\end{figure}

Figure~\ref{fig33} shows the predicted $x$ and the reference $x$ of the trajectory obtained by \verb*|odeint| solver of \verb*|scipy.integrate| library. It can be seen, prediction errors increase significantly in the area of $t>15$. The relative L2 errors are $\epsilon[{x}_{\bm \theta}, {x^{\text{ref}}}] = 9.8\times 10^{-4}$, $\epsilon[{y}_{\bm \theta}, {y^{\text{ref}}}] = 1.4\times 10^{-3}$, and $\epsilon[{z}_{\bm \theta}, {z^{\text{ref}}}] = 6.0\times 10^{-4}$. Table~{\ref{table7}} demonstrates the accuracy of solutions which are given by different approaches for the Lorenz system of differential equations.

\subsection{Discussion}
The presented results clearly indicate that a neural network trained with the help of neuron-by-neuron approach can demonstrate excellent accuracy in solving problems with a relatively short training time. 

In the next section, we extend the approaches described above to tasks with multiple coordinate variables.

\begin{table}[!h]
	\begin{center}
		\begin{tabular}{l|c}
			\hline
			Method \rule[-1ex]{0pt}{3.5ex}  & Relative ${\mathbb L}_2$ errors ($x, y, z$)  \\
			\hline
			Causal training~\cite{Wang2022} \rule[-1ex]{0pt}{3.5ex}  & $1.14\times 10^{-2}$, $1.66\times 10^{-2}$, $7.04\times 10^{-3}$ \\
			{\bf PIDD initialization} & $\bf{2.0\times 10^{-4}}$, $\bf{6.7\times 10^{-4}}$, $\bf{1.3\times 10^{-4}}$\\	
			{\bf  NbN training} & 
			$\bf{9.8\times 10^{-4}}$, $\bf{1.4\times 10^{-3}}$, $\bf{6.0\times 10^{-4}}$ \\
			\hline
		\end{tabular}
		\caption{Lorentz system: Relative ${\mathbb L}_2$ errors obtained by different approaches}\label{table7}
	\end{center}
\end{table}

\newpage
\section{Training of strictly deterministic initialized neural networks for solving PDE}

In our experiments for the solving of PDE we used the above-mentioned SPINN architecture. We discuss solving only 2D problems. The solution of the problem for greater dimensions can be naturally extended based on the solution of 1D and 2D problems.

\subsection{2D problems}

\begin{figure}[th!]
	\centering
	\includegraphics[width=0.3\textwidth]{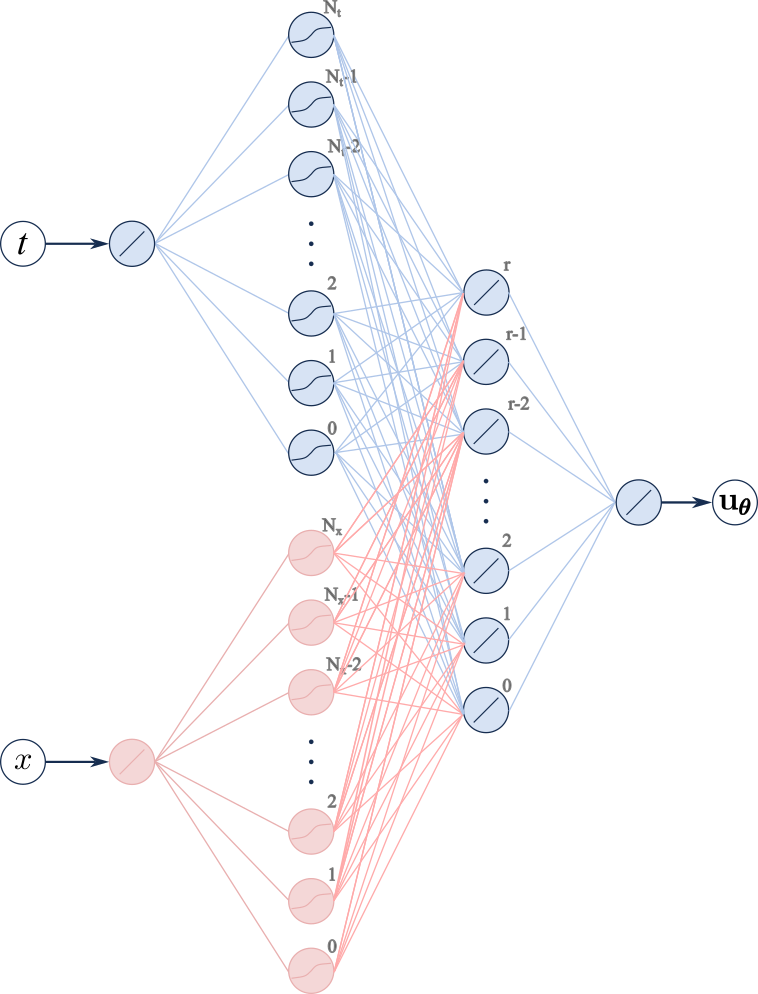}
	\caption{Schematic diagram of the neural network for the SPINN approach at 2D problem.}
	\label{fig26}
\end{figure}

{\bf Suggestion 11:} According SPINN approach for a 2D problem with coordinates $t$ and $x$ the solution can be written as
\begin{equation}
	u_{\bm \theta}(t, x) = \sum\limits_{j=1}^{r} {v}_{j} (t) {w}_{j} (x),\label{eq60}
\end{equation}
where
\begin{eqnarray}\label{eq61}
	&& {v}_{j} (t) = \sum\limits_{k = 0}^{N_t-1} {W}^{(2)}_{t;j,k}\sigma\left({W}^{(1)}_{t;j,k} t + {b}^{(1)}_{t;j,k}\right) + {b}^{(2)}_{t;0},\notag\\
	&& {w}_{j} (x) = \sum\limits_{k = 0}^{N_x-1} {W}^{(2)}_{x;j,k}\sigma\left({W}^{(1)}_{x;j,k} x + {b}^{(1)}_{x;j,k}\right) + {b}^{(2)}_{x;0},
\end{eqnarray}
where ${W}^{(i)}_{t;j,k}$ is $k$th weight of the $i$th layer of $j$th basis-function ${v}_{j}$, ${W}^{(i)}_{x;j,k}$ is $k$th weight of the $i$th layer of $j$th basis-function ${w}_{j}$, $N_t$ is number neurons in the hidden layer of ${v}_{j}$ and  $N_x$ is number neurons in the hidden layer of ${w}_{j}$. A schematic diagram of this neural network is shown in Figure~\ref{fig26}.

\subsection{Numerical experiment. Example 21. Allen--Cahn equation}

Consider a classical phase field modeled by an one-dimensional Allen--Cahn equation with periodic boundary condition
\begin{align}
	& \frac{\partial u}{\partial t} - 0.0001 \frac{\partial^2 u}{\partial x^2} + 5 u^3 - 5u = 0, \quad t \in [0,1],\quad x\in[-1,1],\label{eq62}\\
	& u(0,x) = x^2 \cos(\pi x), \quad u(t,-1) = u(t,1), \quad \left.\frac{\partial u(t,x)}{\partial x}\right|_{x=-1} = \left.\frac{\partial u(t,x)}{\partial x}\right|_{x=1} \notag.
\end{align}
We represent the latent variable $u$ by a SPINN neural network $u_{\bm \theta}$ with $N_t = 1025$ (neurons on hidden layer of ${\vec v}(t)$), $N_x = 513$ (neurons on hidden layer of ${\vec w}(x)$) and number of basis functions is $N_r = 256$. The number of collocation points along the $t$-axis is 8041 and along the $x$-axis is 513. We did not use the Fourier expansion of the input coordinates for our approach. The reference solution can be generated using the Chebfun package~\cite{driscollChebfunGuide2014}, or a dataset from supplementary materials for ~\cite{Wang2022} can be used. The neural networks ${\vec v}(t)$ and ${\vec w}(x)$ were initialized with using Algorithm~\ref{alg2} and initial distribution of $u$: 
\[{W}^{(2)}_{t;j,k} = 2 \dfrac{\Delta t}{\Delta \zeta} \left.\left(0.0001 \dfrac{\partial^2 u}{\partial x^2} - 5 u^3 + 5u\right)\right|_{t=0},\] 
\[{W}^{(2)}_{x;j,k} = 2 \dfrac{\Delta x}{\Delta \zeta} \left.\left(\dfrac{\partial u}{\partial x} \right)\right|_{t=0},\] for every $j$th output.

\begin{figure}[t!]
	\centering
	\includegraphics[width=0.8\textwidth]{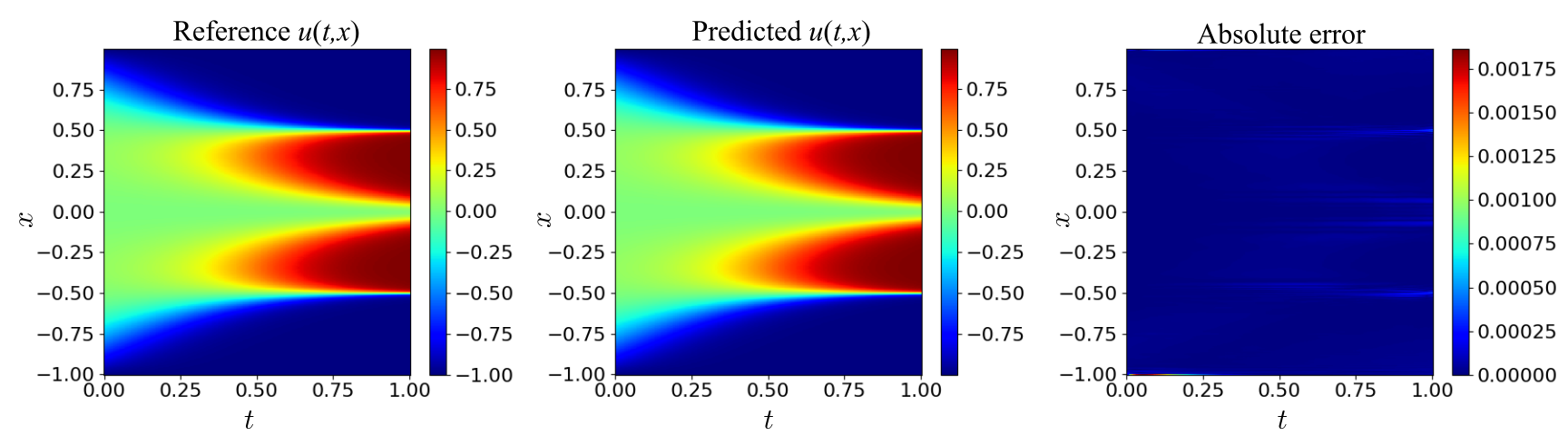}
	\caption{Allen--Cahn equation. (a) is reference solution, (b) is a prediction of a trained physics-informed neural network based on the current approach, (c) is the absolute difference between the reference solution and the predicted solution. The relative error $\epsilon$ is $4.40\times 10^{-5}$.}
	\label{fig35}
\end{figure}

Results of four-stage training which consisted of 5000 epochs of optimization by Adam optimizer with learning rate $10^{-3}$ and 15 epochs of optimization by LBFGS optimizer (maximal number of iterations per optimization step is 25000) under frozen biases of the hidden layer of neural network, then 300000 epochs of optimization by Adam optimizer with learning rate $10^{-5}$ and 15 epochs of optimization by LBFGS optimizer under unfrozen parameters of the hidden layer of neural network, and then 100000 epochs of optimization by Adam optimizer with learning rate $10^{-7}$ and 15 epochs of optimization by LBFGS optimizer under unfrozen parameters of the hidden layer of neural network are shown on the Figure~\ref{fig35}. The following approaches were used simultaneously during the training: detaching, $\delta$-causality (on $t$ and $x$ coordinates (see~\cite{eskin2023optimal}), the initial value is $\varepsilon = 10^{-24}$), additional weighting based on second derivative on $t$ and $x$ coordinates, and relative residuals. To avoid getting stuck in the learning process of a neural network in the area of an incorrect solution, it is proposed to periodically resetting the RR weights to the initial values. Strangely enough, the best convergence is achieved by reset of RR weights at each step of the training. According to method~\cite{eskin2023optimal} the weights $\lambda_{ic}$ and $\lambda_{r}$ are related to each other by means of the relation (\ref{eq26}), where $\beta = N_t / T$, $T=1$. The parameter $\Delta \zeta$ was taken  $ 0.7$. The training was carried out on GPU and took $9875$ seconds of calculation time. The achieved relative error $\epsilon$ is $4.40\times 10^{-5}$.

In table~{8} we also report the accuracy for Allen--Cahn problem achieved using other methods from the literature~\cite{Raissi2019,Wight2020,McClenny2020,Mattey2022,Wang2022,Wang2022_2,Zhang2023}.

\begin{table}[!h]
	\begin{center}
		\begin{tabular}{l|c}
			\hline
			Method \rule[-1ex]{0pt}{3.5ex}  & Relative ${\mathbb L}_2$ error  \\
			\hline
			Original formulation by Raissi et al.~\cite{Raissi2019} \rule[-1ex]{0pt}{3.5ex}  & $4.98\times 10^{-1}$ \\
			Adaptive time sampling~\cite{Wight2020} \rule[-1ex]{0pt}{3.5ex}  & $2.33\times 10^{-2}$ \\
			Causal training (MLP)~\cite{Wang2022} \rule[-1ex]{0pt}{3.5ex}  & $1.43\times 10^{-3}$ \\
			Modified MLP with conservation laws~\cite{Wang2022_2} \rule[-1ex]{0pt}{3.5ex}  & $5.42\times 10^{-4}$ \\
			Causal training (modified MLP)~\cite{Wang2022} \rule[-1ex]{0pt}{3.5ex}  & $1.39\times 10^{-4}$ \\
			{$\delta$-causal training}~\cite{eskin2023optimal} & ${6.29\times 10^{-5}}$ \\			
			{\bf current approach} & $\bf{4.40\times 10^{-5}}$ \\			
			PirateNet~\cite{wang2024} & ${2.24\times 10^{-5}}$ \\
			\hline
		\end{tabular}
		\caption{Allen--Cahn equation: Relative ${\mathbb L}_2$ errors obtained by different approaches}\label{table8}
	\end{center}
\end{table}

\subsection{Physics-informed data-driven initialization}
{\bf Suggestion 11:} If there is a numerical solution to the problem, we can use the neural network (\ref{eq57}) instead of (\ref{eq60}). The functions $U_j$ and $Q_j$ must dependence on $t$ and $x$, respectively.
\subsection{Numerical experiment.}
In our experiments, we used physics-informed data-driven initialization of neural networks, which is described above (Algorithm~\ref{alg4}). 

\subsubsection{Example 22. Allen--Cahn equation}
Consider equation~(\ref{eq62}) with same parameters.

The reference solution was generated using the Chebfun package~\cite{driscollChebfunGuide2014} for $N_t=201$ and $N_x = 512$. The data $\left\{{u}_{kj}\right\}$ for distribution points $\left\{t_k\right\}_{k=0}^{N_t}$ and $\left\{x_j\right\}_{j=0}^{N_x}$. The ${W}^{(2)}_{kj}$ for functions $U_j$ are calculated as follows
\[{W}^{(2)}_{kj} = \dfrac{\Delta x}{2\kappa_k\Delta \zeta} \left(0.0001 \left.\dfrac{\partial^2 u(0,x)}{\partial x^2}\right|_{x_k} - 5 u_{jk}^3 + 5u_{jk}\right).\]
Here $u(0,x)$ is known initial conditions (see (\ref{eq62})).
Results of PIDD initialization for neural network are the following: the relative ${\mathbb L}_2$ errors are $\epsilon[{u}_{\bm \theta}, {u^{\text{ref}}}] = {\bf 6.43{\times} 10^{-3}}$. Such initialization took ${\bf 0.30}$ seconds of computing time of CPU.

\subsubsection{Example 23. Advection equation}
\begin{figure}[t!]
	\centering
	\includegraphics[width=0.8\textwidth]{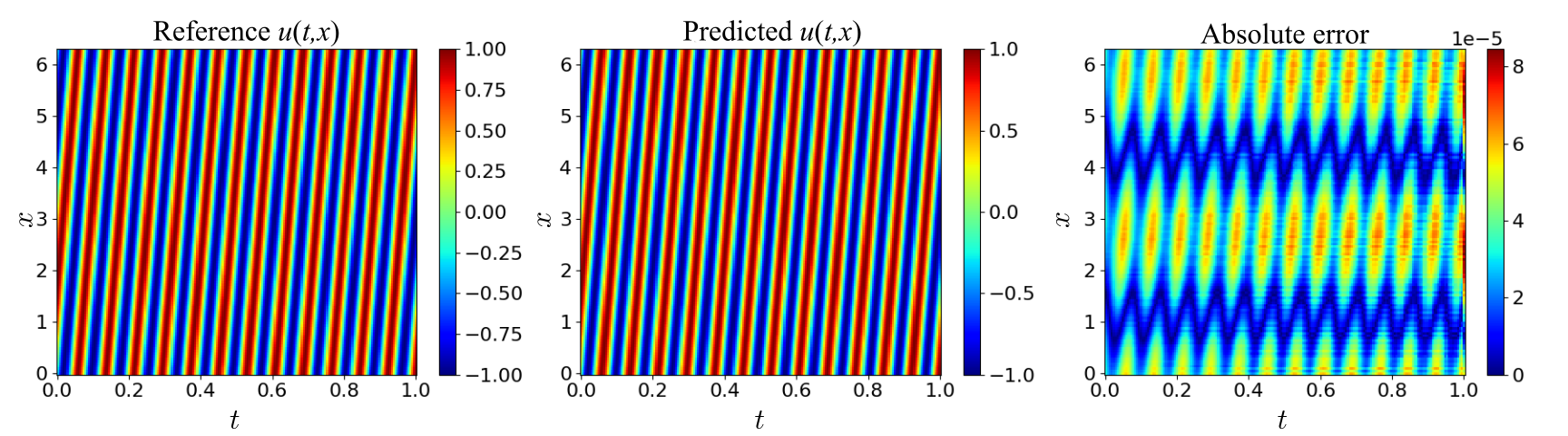}
	\caption{Advection equation. (a) is reference solution, (b) is a prediction of PIDD initialized neural network, (c) is the absolute difference between reference solution and predicted solution. The relative error $\epsilon$ is ${\bf 6.48{\times} 10^{-5}}$.}
	\label{fig36}
\end{figure}

Consider linear hyperbolic equation commonly used to model transport phenomena in the following form
\begin{align}
	& \frac{\partial u}{\partial t} + c \frac{\partial u}{\partial x} = 0, \quad t \in [0,1],\quad x\in[0,2\pi],\label{eq63}\\
	& u(0,x) = \sin(x),\quad x\in[0,2\pi], \quad u(t,0) = u(t,2\pi)\notag.
\end{align}
In our experiments, we consider the challenging setting of
$c = 80$. The reference solution is exact solution $u(t,x) = \sin\left(x - c t\right)$. Results of initialization by algorithm~\ref{alg4} for neural network with $N_t=20000$ and $N_x = 128$ are presented on the Figure~\ref{fig36}. The relative ${\mathbb L}_2$ error is $\epsilon[{u}_{\bm \theta}, {u^{\text{ref}}}] = {\bf 6.48{\times} 10^{-5}}$. Such initialization took ${\bf 0.033}$ seconds of computing time of GPU.

\section{Conclusions}

In this paper, a comprehensive study has been conducted on the development and application of various methods and techniques for the initialization and training of neural networks with one hidden layer, as well as for separable physics-informed neural networks consisting of neural networks with a single hidden layer has been explored, as a means to solve physical problems described by ordinary and partial differential equations.

We have proposed a method for strictly deterministic initialization of a neural network with one hidden layer for solving physical problems described by an ODE. The major advantages of such initialization are the elimination of the need to select a method for random initialization of weights and biases depending on the activation function, controlled network accuracy depending on the number of neurons in the hidden layer, and improved accuracy of the neural network approximations for multi-scale functions. In our opinion, strictly deterministic initialization of a neural network can be effectively using together with Extreme Learning Machine~\cite{HUANG2006489,DWIVEDI202096,SCHIASSI2021334,Huang2011ExtremeLM,lee2024nonoverlappingdomaindecompositionmethod} for increasing accuracy of solution of some problems. In addition, the method of initializing weights and biases of the hidden layer can be used to initialize the hidden layer of a multilayer neural network with several inputs corresponding to its coordinates. With this initialization, the weights and biases of each input are set in accordance with the assignment of weights and offsets of the hidden layer given in Algorithms~\ref{alg1} and~\ref{alg2}.

We have also modified the existing methods of weighting the loss function, such as $\delta$-causal training~\cite{eskin2023optimal} and gradient normalization~\cite{Wang2021}, which improve the accuracy of trained neural networks to predict solutions of considered linear and nonlinear systems of differential equations. New methods for training of strictly deterministic-initialized neural networks to solve ODEs and PDEs, such as detaching, additional weighting based on the second derivative, predicted solution based weighting (PSBW), and relative residuals (RR) have been developed. It is shown that detaching improves the convergence of the learning procedure and increases the accuracy of the neural network trained in this way for nonlinear problems. For linear problems, detaching does not lead to noticeable improvements in the accuracy of the trained model. It is established that the use of additional weighting based on the second derivative leads to significant improvements in the accuracy of the neural network for the problem in which there is a relatively rapid change in the first derivative of the desired solution from the coordinate. Numerical experiments have demonstrated that weighing the terms of the loss function using $\delta$-causal together with PSBW and RR approaches significantly increases the accuracy of a neural network for multi-scale problems. Such an improvement is associated with bringing the elements of the loss function to a uniform scale by such weightings. Note that the PSBW and RR approaches to weighing are inspired by the work~\cite{anagnostopoulos2023residualbased} in which authors proposed an efficient, gradient-less weighting scheme based on evolving cumulative residual. In all numerical experiments using the approaches described in this paragraph, the relations between $\lambda_{r}$ and $\lambda_{ic}$ loss functions were taken based on the method proposed in~\cite{eskin2023optimal}, which gives an unambiguous relationship between these values.

Physics-informed data-driven (PIDD) initialization algorithm for a neural network with one hidden layer was proposed, which allows us to obtain a neural network that provides a solution to the problem with high accuracy, without a tedious data-driven training procedure. A neural network with pronounced generalizing properties was presented, which consists of two neural networks with one hidden layer combined using the SPINN architecture. It was shown the generalizing properties of such neural network can be precisely controlled by changing network parameters. Note that it is relatively simple to construct a neural network with pronounced generalizing properties in the case of parameterization of differential equations by several parameters. To do this, it is enough to create a neural network with one hidden layer for each parameter and combine all neural networks with a SPINN architecture (see Eqn.~\ref{eq57}). A measure of generalization of a neural network has been proposed, which allows us to give a numerical estimation of the generalizing properties of a neural network by the parametrization parameter of differential equations. Such a measure of neural network generalization can be useful in evaluating the generalizing properties of neural networks that differ in architecture~\cite{belbuteperes2021,cho2024} from the neural networks considered in this paper.

A gradient-free neuron-by-neuron (NbN) fitting method has been developed to adjust the parameters of a neural network with one hidden layer, which does not require an optimizer or solver in its application. Using this approach, it was possible to obtain a trained neural network that demonstrates the state-of-the-art results on the Lorentz system.

Experiments on physical problems, such as solving various ODEs and PDEs, have demonstrated that these methods for initializing and training neural networks with one or two hidden layers (SPINN) achieve competitive accuracy and, in some cases, state-of-the-art results.

We hope that our research will unlock the full potential of neural networks with one or two hidden layers and develop more accurate, effective and reliable methods to solve complex physical problems for them.

\newpage

\appendix

\section{Measure of accuracy}\label{appA}
To evaluate the accuracy of the approximate solution obtained with the help of PINN method, the values of the solution of (\ref{eq1}) predicted by the neural network at given points are compared with the values calculated on the basis of classical high-precision numerical methods. As a measure of accuracy, the relative total ${\mathbb L}_2$ error of prediction is taken, which can be expressed with the following relation
\begin{equation}\label{eq61}
	\epsilon[{u}_{\bm \theta}, {u}] = \left\{\frac{1}{N_e} \sum^{N_e}_{i=1} \left[{u}_{\bm \theta}({\vec x}_i) - {u}({\vec x}_i)\right]^2 \right\}^{1/2} {\times} \left\{\frac{1}{N_e} \sum^{N_e}_{i=1} \left[{u}({\vec x}_i)\right]^2 \right\}^{-1/2},
\end{equation}
where $\left\{{\vec x}_{i}\right\}^{N_{e}}_{i=1}$ is the set of evaluation points taken from the domain $\Omega$, ${u}_{\bm \theta}$ and ${u}$ are the predicted and reference solutions respectively.

\section{Measure of generalization properties}\label{appB}
To evaluate the generalization properties of the approximate solution obtained using the PINN method, two sets of solution values (\ref{eq1}) predicted by a neural network on homogeneous grids of this parameter are required for the parameter of differential equations $q$. The first set of solution values is calculated for the grid of the parameter $q$ in steps of $\Delta q$ for those values of $q$ at which the training was carried out. We denote this set of parameters as $\left\{{\vec q}_{i}\right\}^{N_{q}}_{i=1}$. The second set of solution values is calculated for the calculated grid of the parameter $q$ in steps of $\Delta q$, offset by $\Delta q/ 2$, i.e. at points as far away as possible from the nearest points of $q$ where the training was carried out, but not out of the training domain. We will denote this set of parameters as $\left\{{\vec q}_{i+1/2}\right\}^{N_{q}-1}_{i=1}$. As a measure of generalization properties, the ration of the relative total ${\mathbb L}_2$ error of prediction of set $\left\{{\vec q}_{i+1/2}\right\}^{N_{q}-1}_{i=1}$ to the ration of the relative total ${\mathbb L}_2$ error of prediction of set $\left\{{\vec q}_{i}\right\}^{N_{q}}_{i=1}$  is taken, which can be expressed with the following relation
\begin{align}\label{eq65}
	& \mu[{u}_{\bm \theta}, {u}] = \frac{\epsilon_1[{u}_{\bm \theta}, {u}]}{\epsilon_{1/2}[{u}_{\bm \theta}, {u}]}.
\end{align}
Here
\begin{align}
	& \epsilon_{1/2}[{u}_{\bm \theta}, {u}] = \left\{\sum^{N_q -1}_{j=1}\sum^{N_e}_{i=1} \left[{u}_{\bm \theta}({\vec x}_i,q_{j+1/2}) - {u}({\vec x}_i,q_{j+1/2})\right]^2 \right\}^{1/2} {\times} \left\{ \sum^{N_q -1}_{j=1}\sum^{N_e}_{i=1} \left[{u}({\vec x}_i,q_{j+1/2})\right]^2 \right\}^{-1/2},\notag\\
	&  \epsilon_{1}[{u}_{\bm \theta}, {u}] = \left\{\sum^{N_q}_{j=1}\sum^{N_e}_{i=1} \left[{u}_{\bm \theta}({\vec x}_i,q_{j}) - {u}({\vec x}_i,q_{j})\right]^2 \right\}^{1/2} {\times} \left\{ \sum^{N_q}_{j=1}\sum^{N_e}_{i=1} \left[{u}({\vec x}_i,q_{j})\right]^2 \right\}^{-1/2}, \notag
\end{align}
where $\left\{{\vec x}_{i}\right\}^{N_{e}}_{i=1}$ is the set of evaluation points taken from the domain $\Omega$, $\left\{{q}_{i}\right\}^{N_{q}}_{i=1}$ and $\left\{{q}_{i+1/2}\right\}^{N_{q}}_{i=1}$ are sets taken from domain $[q_{\rm min}, q_{\rm max}]$, ${u}_{\bm \theta}$ and ${u}$ are the predicted and reference solutions, respectively, second arguments of function ${u}_{\bm \theta}$ and ${u}$ means parameter of differential equation at which this functions had been calculated.

Thus, the closer the value of the measure of generalization properties $\mu$ is to one, the neural network outputs solutions to the problem at the parameters of the ODE at which it was not trained, closer in accuracy to its solutions to the problem at the parameters of the ODE at which it was trained.

\newpage
\section{Nomenclature}\label{appE}

The summarizes the main notations, abbreviations and symbols are given in table~9.
\begin{table}[h!]
	\centering
	\begin{tabular}{ll}
		\toprule
		Notation & Description  \\
		\midrule
		GN1& Gradient Normalization method on Eqn.~(\ref{eq32})\\
		GN2& Gradient Normalization method on Eqn.~(\ref{eq34})\\
		MLP & Multilayer perceptron  \\
		NbN & Neuron-by-neuron training\\
		NN & neural network\\
		ODE & Ordinary differential equations \\
		PDE & Partial differential equation\\
		PIDD & Physics-informed data-driven\\
		PINN & Physics-informed neural network \\
		PSBW & Predicted solution based weighting\\
		RR & Relative Residuals\\		
		SPINN & Separable physics-informed neural network \\
		${\vec u}$ & solution of PDE	\\
		${\mathcal{N}}$ & nonlinear differential operator\\
		$\mathcal{R}$ & PDE residual \\
		${\vec u}_{\bm \theta}$ & neural network representation of the ODE or PDE solution \\
		${\bm \theta}$ & vector of the trainable parameters of the neural network \\
		$w_i$ & residual weight at the time $t_i$\\
		$\varepsilon$ & causality parameter\\
		$\rasymbol{f}$ & Ra symbol over value $f$ means calculated value $f$ is detached from calculation graph.\\
		& This operation corresponding to \verb*|detach| in PyTorch, and \verb*|lax.stop_gradient| in JAX \\		
		${\mathcal{L}_r^{(\rm t)}\left(t,{\bm \theta}\right)}$ & temporal residual loss\\
		${\mathcal{L}(\bm \theta)}$ & aggregate training loss\\
		$\delta$-causal & Dirac delta function causal training~\cite{eskin2023optimal}\\
		$D^2$ & weighting based on second derivative(Eqns. (\ref{eq36_2}) and (\ref{eq36_3}))\\
		\bottomrule
	\end{tabular}\label{tableAppA}
	\caption{Nomenclature}
\end{table}

\newpage
\bibliographystyle{IEEEtran}
\bibliography{Eskin_tinyPINN}

\end{document}